\documentclass{amsart}

\usepackage{xy}
\usepackage{amssymb}
\DeclareMathAlphabet{\mathscr}{U}{rsfs}{m}{n}

\title[Localization in $THH$ and $TC$]%
{Localization theorems in topological Hochschild homology and
topological cyclic homology}

\author{Andrew J. Blumberg}
\address{Department of Mathematics, The University of Texas,
Austin, TX \ 78712}
\email{blumberg@math.utexas.edu}
\thanks{The first author was supported in part by an NSF postdoctoral
fellowship, NSF grant DMS-0111298, and a Clay Liftoff Fellowship}
\author{Michael A. Mandell}
\address{Department of Mathematics, Indiana University,
Bloomington, IN \ 47405}
\email{mmandell@indiana.edu}
\thanks{The second author was supported in part by NSF grants
DMS-0504069, DMS-0804272}


\newcommand{\ssdot}{\bullet}
\newcommand{\supdot}{^\ssdot}
\newcommand{\subdot}{_\ssdot}

\newcommand{\Sdot}[1][{\ssdot}]{S_{#1}}
\newcommand{\barw}{\bar w}

\newcommand{\orth}{\aI\sS}
\newcommand{\symm}{\Sigma\sS}

\newcommand{\cuph}{\mathbin{\mathop{\cup}\limits^{h}}}

\newcommand{\parf}{\mathrm{parf}}
\newcommand{\DSparf}[1][{X}]{\aD^{S}_{\parf}({#1})}
\newcommand{\DSparflat}[1][{X}]{\aD^{S}_{\parf}({#1})_{\flat}}
\newcommand{\DGparf}[1][{X}]{\aD^{DG}_{\parf}({#1})}
\newcommand{\DGparflat}[1][{X}]{\aD^{DG}_{\parf}({#1})_{\flat}}

\newcommand{\RMod}[1]{\mathfrak{Mod}_{#1}}

\newcommand{\THM}{N^{\cy}}
\newcommand{\CTHH}{CTHH}
\newcommand{\CTHM}{C\THM}
\newcommand{\WTHH}{WTHH}
\newcommand{\WT}{WT}
\newcommand{\WTC}{WTC}
\newcommand{\TB}{B}

\newcommand{\Cell}{\mathrm{Cell}}

\newcommand{\Fix}[1]{\mathrm{Fix}^{C_{#1}}}

\let\iso\cong
\let\sma\wedge
\newcommand{\htp}{\simeq}
\renewcommand{\to}{\mathchoice{\longrightarrow}{\rightarrow}{\rightarrow}{\rightarrow}}
\newcommand{\from}{\mathchoice{\longleftarrow}{\leftarrow}{\leftarrow}{\leftarrow}}

\newcommand{\bm}{\mathbf{m}}
\newcommand{\bn}{\mathbf{n}}
\newcommand{\bp}{\mathbf{p}}
\newcommand{\bz}{\mathbf{0}}

\let\catsymbfont\mathcal
\newcommand{\aA}{{\catsymbfont{A}}}
\newcommand{\aB}{{\catsymbfont{B}}}
\newcommand{\aC}{{\catsymbfont{C}}}
\newcommand{\aD}{{\catsymbfont{D}}}
\newcommand{\aE}{{\catsymbfont{E}}}
\newcommand{\aF}{{\catsymbfont{F}}}
\newcommand{\aG}{{\catsymbfont{G}}}
\newcommand{\aI}{{\catsymbfont{I}}}
\newcommand{\aK}{{\catsymbfont{K}}}
\newcommand{\aM}{{\catsymbfont{M}}}
\newcommand{\aN}{{\catsymbfont{N}}}
\newcommand{\aP}{{\catsymbfont{P}}}
\newcommand{\aQ}{{\catsymbfont{Q}}}
\newcommand{\aS}{{\catsymbfont{S}}}

\newcommand{\aU}{{\catsymbfont{U}}}

\newcommand{\aZ}{{\catsymbfont{Z}}}

\newcommand{\bC}{{\mathbb{C}}}
\newcommand{\bP}{{\mathbb{P}}}

\newcommand{\bZ}{{\mathbb{Z}}}

\let\sheafsymbfont\mathcal
\newcommand{\sE}{{\sheafsymbfont{E}}}
\newcommand{\sO}{{\sheafsymbfont{O}}}
\newcommand{\sJ}{{\mathscr{J}}}
\newcommand{\sS}{{\mathscr{S}}}

\let\fusymbfont\mathcal
\newcommand{\fL}{{\fusymbfont{L}}}
\newcommand{\fM}{{\fusymbfont{M}}}
\newcommand{\fN}{{\fusymbfont{N}}}
\newcommand{\fP}{{\fusymbfont{P}}}
\newcommand{\fQ}{{\fusymbfont{Q}}}

\newcommand{\fLL}[1][{}]{\fL^{\aB_{#1}}_{\aA_{#1}}}
\newcommand{\fQQ}[1][{}]{\fQ^{\aB_{#1}}_{\aA_{#1}}}
\newcommand{\fLM}{\fL_{\aA_{c}}^{\aM_{c}}}
\newcommand{\fQM}{\fQ_{\aA_{c}}^{\aM_{c}}}

\def\quickop#1{\expandafter\DeclareMathOperator\csname
#1\endcsname{#1}}
\quickop{id}\quickop{Id}\quickop{colim}\quickop{hocolim}\quickop{holim}\quickop{op}
\quickop{Ar}\quickop{sing}\quickop{Hom}\quickop{w}\quickop{Ho}
\quickop{ob}\quickop{cy}\quickop{Spec}\quickop{Lan}\quickop{diag}
\DeclareMathOperator{\Coker}{Coker}
\quickop{sd}

\DeclareMathOperator*{\lcolim}{colim}

\newcommand{\on}{\mathbin{\mathrm{on}}}

\numberwithin{equation}{section}

\newtheorem{thm}{Theorem}
\numberwithin{thm}{section}

\newtheorem*{thm*}{Theorem}
\newtheorem{cor}[thm]{Corollary}
\newtheorem{lem}[thm]{Lemma}
\newtheorem{prop}[thm]{Proposition}
\theoremstyle{definition}
\newtheorem{defn}[thm]{Definition}
\newtheorem{notn}[thm]{Notation}
\newtheorem*{conv}{Convention}
\theoremstyle{remark}
\newtheorem{rem}[thm]{Remark}

\makeatletter\let\c@equation\c@thm\makeatother

\xyoption{arrow}
\xyoption{curve}
\xyoption{matrix}
\xyoption{cmtip}
\SelectTips{cm}{}

\newcommand{\term}[1]{\textit{#1}}

\newdir{ >}{{}*!/-5pt/\dir{>}}

\begin{document}

\begin{abstract}
We construct localization cofibration sequences for the topological
Hochschild homology ($THH$) and topological cyclic homology ($TC$) of
small spectral categories.  Using a global construction of the $THH$ and
$TC$ of a scheme in terms of the perfect complexes in a spectrally
enriched version of the category of unbounded complexes, the sequences
specialize to localization cofibration sequences associated to the
inclusion of an open subscheme.  These are the targets of the
cyclotomic trace from the localization sequence of Thomason-Trobaugh
in $K$-theory.  We also deduce versions of Thomason's blow-up formula
and the projective bundle formula for $THH$ and $TC$.
\end{abstract}

\maketitle

\tableofcontents

\section{Introduction}

Algebraic $K$-theory provides a powerful and subtle invariant of
schemes.  The $K$-theory of a scheme encodes many of its arithmetic
and algebraic properties, captures information about its geometry and
singularities, and is closely connected to its \'etale and motivic
cohomology.  One of the fundamental underpinnings of the subject is
the localization theorem of Thomason and Trobaugh
\cite[7.4]{ThomasonTrobaugh}, which for a quasi-separated quasi-compact
scheme $X$ provides a cofibration sequence of (non-connective) $K$-theory
spectra 
\[ K(X \on (X-U)) \to K(X) \to K(U) \to \Sigma K(X \on (X-U)),\]
for $U$ a quasi-compact open subscheme contained in $X$.  Here $K(X \on (X-U))$
denotes the $K$-theory of the category of perfect complexes on 
$X$ which are supported on the complement of
$U$ in $X$.  This localization sequence and the
closely related Mayer-Vietoris sequence for $K$-theory
allow global assembly of local information.

Keller \cite{KellerCyclic} constructed the analogue of the
Thomason-Trobaugh localization sequence for Hochschild homology ($HH$)
and for the variants of cyclic homology, including negative cyclic
homology ($HC^{-}$).  The Dennis trace (or Chern character) connects
the localization sequence in $K$-theory to the localization sequence
in $HC^{-}$.  Using this, together with generalizations to blow-ups
along regular sequences and Hironaka's resolution of singularities,
Corti\~nas, Haesemeyer, Schlichting, and Weibel 
\cite{Cortinas,CHSW,CHW} recently resolved Weibel's conjecture
bounding below the negative $K$-groups and Vorst's conjecture that
$K_{d+1}$-regularity implies regularity, for finite-type schemes of
dimension $d$ over a field of characteristic zero.

The purpose of this paper is to generalize Keller's localization
sequences to topological Hochschild homology ($THH$) and topological
cyclic homology ($TC$).  Over the course of the last two decades,
$THH$ and $TC$ have revolutionized $K$-theory computations.  Roughly,
topological Hochschild homology for a ring is obtained by promoting the
ring to a ring spectrum and substituting the smash product of spectra
for the tensor product of rings in the Hochschild complex
\cite{Bokstedt}.  The $THH$ spectrum comes with a 
``cyclotomic'' structure (which involves an $S^{1}$-action and extra
structure maps), and for each prime $p$, topological cyclic homology
is then defined as a certain homotopy limit over the fixed point
spectra.  The Dennis trace map lifts to a ``cyclotomic trace'' map
from $K$-theory to $TC$ \cite{BHM}, and McCarthy \cite{McCarthy}
showed that this captures all the relative information at $p$ for
surjections with nilpotent kernel, just as $HC^{-}$ does rationally
\cite{GoodwillieHN}.  Starting from Quillen's computation of the
$K$-theory of finite fields, Hesselholt and Madsen have used $TC$ to
make extensive computations in $K$-theory \cite{HM1,HM2,HM3}.
Moreover, because of the close relationship between $K$-theory and
$TC$ (and analogy with $HC^{-}$), this paper provides the key
ingredients needed to generalize the work of Corti\~nas, Haesemeyer,
Schlichting, and Weibel \cite{CHSW, CHW} to cases in characteristic $p$
where resolution of singularities holds.  Geisser and Hesselholt have
already started applying the results of this paper in this direction
\cite{GHWeib}.

Between $TC$ and $THH$ is an intermediate theory called
$TR$, whose homotopy groups have the structure of a Witt complex (the
structure whose universal example is the de\,Rham-Witt complex of
Bloch-Deligne-Illusie).  The Hesselholt-Madsen computations proceed by
studying this structure on $TR$. Hesselholt has observed that in all
known examples, the de\,Rham-Witt 
complex has the same relationship to $TR$ that Milnor $K$-theory has
to algebraic $K$-theory.  This led Geisser and Hesselholt to conjecture
an ``additive'' motivic spectral sequence converging to a modified
version of $TR$ with edge homomorphism the universal map from the
de\,Rham-Witt complex.  Recent work of Levine \cite{LevineHCT}
axiomatizes the role of localization and Mayer-Vietoris theorems in
the construction of the 
motivic spectral sequence \cite{BlochLichtenbaum,FriedlanderSuslin},
and such theorems for $TR$ should provide key input to the
construction of this conjectural ``additive'' motivic spectral
sequence.  We prove the following results in this direction.

\begin{thm}\label{inttt}
Let $X$ be a quasi-compact and semi-separated scheme.  For a
quasi-compact open subscheme 
$U$, there are homotopy cofibration sequences 
\begin{gather*}
THH(X \on (X-U)) \to THH(X) \to THH(U)\to \Sigma THH(X \on (X-U))\\[1ex]
TR(X \on (X-U)) \to TR(X) \to TR(U) \to \Sigma TR(X \on (X-U)) \\[1ex]
TC(X \on (X-U)) \to TC(X) \to TC(U) \to \Sigma TC(X \on (X-U)) 
\end{gather*}
where $THH(X \on (X-U))$ denotes the $THH$ of the spectral category of perfect
complexes on $X$ which are supported on $X-U$.

For quasi-compact open subschemes $U,V$ with $X=U\cup V$, the squares
\[ \scriptstyle 
\xymatrix@C=1em{
THH(X) \ar[r] \ar[d] & THH(U) \ar[d] 
&TR(X) \ar[r] \ar[d] & TR(U) \ar[d]
&TC(X) \ar[r] \ar[d] & TC(U) \ar[d]\\
THH(V) \ar[r] & THH(U \cap V) 
&TR(V) \ar[r] & TR(U \cap V) 
&TC(V) \ar[r] & TC(U \cap V) 
}
\]
are homotopy cocartesian.
\end{thm}

In the statement, a scheme is \term{semi-separated}
\cite[B.7]{ThomasonTrobaugh} when it has a basis of affine open
subsets whose intersections are also affine.  Semi-separated is a
slightly stronger condition than \term{quasi-separated} which means
that it has a basis of affine open subsets whose intersection is a
finite union of affine open subsets.  If a scheme has an ample family
of line bundles then it is semi-separated
\cite[B.7]{ThomasonTrobaugh}.

Geisser and Hesselholt \cite{GeisserHesselholt} proved the second
statement in Theorem~\ref{inttt} for $THH$ of rings and used it to
define $THH$ of quasi-compact quasi-separated schemes in 
terms of Thomason's hypercohomology
construction \cite[1.33]{ThomasonEtale}.  The relative term $THH(X \on
(X-U))$ does not have an intrinsic description in the context of the
Geisser-Hesselholt definition of $THH$.  Here we 
describe it in terms of a construction of $THH$ for \term{spectral
categories}, i.e., categories enriched over symmetric spectra, the
stable homotopy theory refinement of DG-categories.

Dundas and McCarthy \cite{DundasMcCarthy} generalized B\"okstedt's
construction of $THH$ to small spectral categories.  We build on the
foundations there and study more general invariance properties; see in
particular Theorems~\ref{thmdkequiv} and~\ref{thmthick} below.
We use these invariance properties to generalize the localization
theorem of Keller to the setting of spectral categories.
Roughly, we show that the $THH$ of a triangulated quotient is the
cofiber on $THH$; Theorems~\ref{thmgenone} and~\ref{thmgentwo} provide
precise statements.  Although we work in the context of spectral categories,
our localization theorem specializes to the setting of DG-categories,
as small DG-categories may be functorially converted to small spectral categories
with the same objects and spectral refinements of the Hom complexes;
see for example \cite[\S6]{SSMonoidalEq} or
\cite[App.~A]{DuggerShipleyEnriched}, among others. 
Just as $THH$ of a ring captures much more torsion information than
$HH$ of the ring, $THH$ provides a richer invariant of a
DG-category than $HH$.  Using an appropriate functor from
small DG-categories to small
spectral categories and DG-categories, we make the following
observation at the end of Section~\ref{secthh}.

\begin{thm}\label{intthhfunct}
The constructions of $THH$, $TR$, and $TC$ as defined in
Section~\ref{secthh} are functors from the category of small
DG-categories and DG-functors to the stable category.
\end{thm}

We define $THH$ of a scheme in terms of a spectral category refinement
$\DSparf$ of the DG-category quotient $\DGparf$ modeling the derived
category of perfect complexes.  In Section~\ref{secappl}, we prove the
following consistency theorem that compares this definition to the
definition of Geisser-Hesselholt.

\begin{thm}\label{intgh}
Let $X$ be a quasi-compact and semi-separated scheme, and $\DSparf$ a
spectral category refinement of $\DGparf$.  Then $THH(\DSparf)$ is
equivalent to the Thomason hypercohomology of the presheaf of
symmetric spectra $U \mapsto
THH(\sO_{U})$ on the small Zariski site of $X$.
\end{thm}

This theorem in particular constructs a trace map from the $K$-theory
of the scheme to $THH(\DSparf)$ and $TC(\DSparf)$.  In
Section~\ref{secnoncon}, we show that the trace map factors through
Thomason-Trobaugh's Bass' non-connective $K$-theory spectrum using
their spectral version of Bass' fundamental theorem.  In
Section~\ref{apptrace}, we give a direct
construction of the trace for $\DSparf$ that does not use the
hypercohomology construction.

In addition to Theorem~\ref{inttt}, we also establish $THH$ and $TC$
versions of two classical geometric calculations in algebraic
$K$-theory using our general localization machinery.  First, we prove
the following formula for blow-ups along regular sequences, which
already has been applied by Geisser and Hesselholt \cite{GHWeib} to
prove the characteristic $p$ analogue of Weibel's conjecture.  We
state the theorem using the notation of \cite[\S 1]{CHSW}, and prove
it in Section~\ref{secappl}.

\begin{thm}\label{intbl}
Let $X$ be a quasi-compact and semi-separated scheme.
Let $i\colon Y\subset X$ be a regular embedding of a closed subscheme, 
$p\colon X\to X'$ the blowup along $Y$, $j\colon Y'\subset X'$ the
exceptional divisor, and write $q$ for the map $Y'\to Y$.  Then the squares
\[
\xymatrix{%
THH(X)\ar[r]^{Lp^{*}}\ar[d]_{Li^{*}}&THH(X')\ar[d]^{Lj^{*}}
&TC(X)\ar[r]^{Lp^{*}}\ar[d]_{Li^{*}}&TC(X')\ar[d]^{Lj^{*}}\\
THH(Y)\ar[r]_{Lq^{*}} &THH(Y')
&TC(Y)\ar[r]_{Lq^{*}} &TC(Y')
}
\]
are homotopy cocartesian.
\end{thm}

We also prove a projective bundle
theorem \cite[4.1,7.3]{ThomasonTrobaugh} in Section~\ref{secappl}.

\begin{thm}\label{intproj}
Let $X$ be a quasi-compact and semi-separated scheme.  Let $\sE$ be
an algebraic vector bundle of rank $r$ over $X$, and let
$\pi\colon \bP \sE_X \to X$ be the associated projective 
bundle.  Then a spectral lift of the derived functor
\[
\bigoplus_{i=0}^{r-1} \sO_{\bP \sE_X} (-i) \otimes L\pi^*(-)
\]
induces weak equivalences
\[
\prod_{i=0}^{r-1} THH(X) \to THH(\bP \sE_X) \quad \text{and}\quad \prod_{i=0}^{r-1} TC(X) \to TC(\bP \sE_X).
\]
\end{thm}

The proof of the previous two theorems require the machinery of the
construction of $THH$ of small spectral categories that we develop in this
paper; it is not known how to prove them using just a hypercohomology
construction like that of Geisser-Hesselholt \cite{GeisserHesselholt}.

The paper is organized as follows.  In Section~\ref{secspecc}, we
review the basic definitions for spectral categories (categories
enriched in symmetric spectra).  As indicated above, this is the
appropriate setting for studying $THH$, $TR$, and $TC$, and is a
stable homotopy theory generalization of the setting of DG-categories.
In Section~\ref{secthh}, we review the definition of $THH$ of small spectral
categories due to B\"okstedt \cite{Bokstedt} and Dundas-McCarthy
\cite{DundasMcCarthy}.  Because of the work of Shipley
\cite{ShipleyD}, the technical hypotheses of connectivity and
convergence on the input symmetric spectra that seemed necessary for
the last 20 years may now be omitted.  We take the viewpoint, first
articulated by Dwyer and Kan, that enriched mapping spaces (or
spectra) encode the ``higher homotopy theory'' of a category, and we
view $THH$, $TR$, and $TC$ as invariants of the higher homotopy theory
of the category, as is $K$-theory \cite{ToenVezzosi,BlumbergMandell2}.
Section~\ref{secdefcyc} spells out in detail the definition of the
point-set category of cyclotomic spectra of orthogonal spectra. 
In Section~\ref{secdkequiv}, we list several invariance theorems for
$THH$ in this context.  Section~\ref{secdwm} reviews an elementary
tilting argument for $THH$, Proposition~\ref{propcoremorita},
originally due to Dennis and Waldhausen \cite[p.~391]{WaldhausenA2}.
We demonstrate how to apply the tilting argument to prove powerful
comparison theorems.  Using these techniques, in
Section~\ref{secgenloc} we prove the general localization
theorems~\ref{thmgenone} and~\ref{thmgentwo}, which we apply in
Section~\ref{secappl} to prove Theorems~\ref{inttt}, \ref{intgh},
\ref{intbl}, and \ref{intproj} above.  In Section~\ref{secnoncon}, we
extend the cyclotomic trace over Bass' non-connective $K$-theory, using
Thomason and Trobaugh's spectral version of Bass' fundamental theorem.
In order to simplify the discussion, we use the ad hoc version of the
cyclotomic trace for schemes in \cite{GeisserHesselholt} induced by
the cyclotomic trace for rings; Section~\ref{apptrace} constructs a
canonical version of the cyclotomic trace map for arbitrary
DG-Waldhausen categories. This requires a review of Waldhausen's
$S\subdot$ construction and the construction of algebraic $K$-theory.
The last section gives a version of Theorem~\ref{thmgenone} that is
more useful in the context of spectral model categories.

\bigskip

The authors would like to thank the Department of Mathematics and the
Mathematics Research Center at Stanford University and the Institute
for Advanced Study for their hospitality and support while some of
this work was being done.  The authors would like to thank Lars
Hesselholt for asking motivating questions and for sharing his ideas
in this direction, as well as Christian Haesemeyer, John Rognes, Marco
Schlichting, Brooke Shipley, and Charles Weibel for interesting and
useful conversations.

\section{Review of spectral categories}\label{secspecc}

Modern constructions of the stable category with point-set level smash
products allow easy generalization of the concepts of simplicial
category or DG-category to the context of spectra.  Symmetric
spectra in particular often arise naturally as the refinement of
mapping sets.  In fact, symmetric ring spectra (the analogue of
DG-rings) and categories enriched in symmetric spectra (the analogue
of DG-categories) predated Smith's insight that the homotopy theory of
symmetric spectra models the stable category.  In older $K$-theory
literature, they were called FSPs (or FSPs defined on spheres) and
FSPs with many objects, respectively, and treatments generally
included hypotheses on connectivity or convergence.  A modern approach
to $THH$ and $TC$, taking advantage of \cite{HSS} and especially
\cite{ShipleyD} obviates the need for any such connectivity or
convergence hypotheses.  In this section, we review the definition of
spectral categories, and modules and bimodules over spectral categories
in terms of enriched category theory.

\begin{defn}
A \term{spectral category} is a category enriched over symmetric
spectra in simplicial sets.  Specifically, a spectral category $\aC$
consists of: 
\begin{enumerate}
\item A collection of objects $\ob\aC$ (which need not be a small set),
\item A symmetric spectrum $\aC(a,b)$ for each pair of objects $a,b\in \ob\aC$,
\item A unit map $S\to \aC(a,a)$ for each object $a \in \ob\aC$, and
\item A composition map $\aC(b,c)\sma \aC(a,b)\to \aC(a,c)$ for each
triple of objects $a,b,c\in \ob\aC$,
\end{enumerate}
satisfying the usual associativity and unit properties.  We say that a
spectral category is small when the objects $\ob \aC$ form a small set.
\end{defn}

We emphasize that the data in (iii) and (iv) consist of point-set
maps (rather than maps in the stable category) and that ``$\sma$''
denotes the point-set smash product of symmetric spectra.  
The definition of spectral functor between spectral categories is the
usual definition of an enriched functor:

\begin{defn}
Let $\aC$ and $\aD$ be spectral categories.  A \term{spectral functor}
$F\colon \aC\to \aD$ is an enriched functor.  Specifically, a spectral
functor consists of:
\begin{enumerate}
\item A function on objects $F\colon \ob\aC\to \ob\aD$, and
\item A map of symmetric spectra $F_{a,b}\colon \aC(a,b)\to
\aD(Fa,Fb)$ for each pair of objects $a,b\in \ob\aC$,
\end{enumerate}
which is compatible with the units and the compositions in the obvious sense.
\end{defn}

Again, we emphasize that the compatibility condition holds in the
point-set category of symmetric spectra rather than in the stable
category.  We use the term \term{weak equivalence} to mean a
spectral functor that is a bijection on objects and a weak equivalence
(stable equivalence of symmetric spectra) on all mapping spectra.  See
Definition~\ref{defdkequiv} for a more general kind of equivalence.

We have the evident concepts of module and bimodule over
spectral categories: 

\begin{defn}
Let $\aC$ and $\aD$ be spectral categories.  A left $\aC$-module
is a spectral functor from $\aC$ to symmetric spectra.  A right $\aD$-module
is a spectral functor from $\aD^{\op}$ to symmetric
spectra.  A $(\aD,\aC)$-bimodule is a spectral functor from
$\aD^{\op}\sma \aC$ to symmetric spectra.
\end{defn}

Here $\aD^{\op}$ denotes the spectral category with the same objects
and mapping spectra as $\aD$ but the opposite composition map.  The
spectral category $\aD^{\op}\sma \aC$ has as its objects the cartesian
product of the objects, 
\[
\ob (\aD^{\op}\sma \aC)=\ob\aD^{\op}\times
\ob\aC,
\]
and as its mapping spectra the smash product of the mapping
spectra
\[
(\aD^{\op}\sma \aC)((d,c),(d',c'))=\aD^{op}(d,d')\sma\aC(c,c'),
\]
with unit maps the smash product of the unit maps and composition maps
the smash product of the composition maps for $\aD^{\op}$ and $\aC$.
Explicitly, a $(\aD,\aC)$-bimodule $\fM$ consists of a choice of
symmetric spectrum $\fM(d,c)$ for each $d$ in $\ob \aD$ and $c$ in
$\ob \aC$, together with maps
\[
\aC(c, c') \sma \fM(d,c) \sma \aD(d',d) \to \fM(d',c')
\]
for each $d'$ in $\ob \aD$ and $c'$ in
$\ob \aC$, making the obvious unit and associativity diagrams
commute.  In particular, for any spectral category $\aC$, the mapping
spectra $\aC(-,-)$ define a $(\aC,\aC)$-bimodule. (This example
motivates the convention of listing the right module structure first.)

The work of \cite{SSMonoidalEq} provides the category of
$(\aD,\aC)$-bimodules with a closed model structure.

\begin{prop}\textup{(\cite[6.1]{SSMonoidalEq})}
The category of $(\aD,\aC)$-bimodules forms a closed model category
where the fibrations are the objectwise fibrations and the weak
equivalences are the objectwise weak equivalences in the stable model
structure on symmetric spectra.
\end{prop}

Older $K$-theory literature required ``convergence'' hypotheses on spectral
categories and bimodules, asking for the homotopy groups of the
constituent spaces in each mapping spectrum to stabilize.  These
hypotheses appeared necessary at the time to analyze the homotopy
colimits arising in B\"okstedt's construction of $THH$.  It was thought
that these homotopy colimits could be wrong
for a non-convergent symmetric spectrum because the homotopy groups
they computed generally differed from the homotopy groups expected
from the underlying prespectrum.  Because of \cite{HSS,ShipleyD}, we
now understand that it is the homotopy groups of the underlying
prespectrum that may be wrong: The homotopy groups of the prespectrum
underlying a symmetric spectrum $X$ do not necessarily agree with the
homotopy groups of the object represented by $X$ in the stable
category.  In general, every symmetric spectrum $X$ admits a 
weak equivalence $X\to \tilde X$ to a symmetric $\Omega$-spectrum
$\tilde X$, i.e., one whose underlying prespectrum is an
$\Omega$-spectrum (level fibrant with adjoint structure maps $\tilde
X_{n}\to \Omega \tilde X_{n+1}$ weak equivalences). The correct
homotopy groups of $X$ are the homotopy groups of the underlying
prespectrum of $\tilde X$; when these agree under the comparison map
with the homotopy groups of the underlying prespectrum of $X$, then
$X$ is said to be \term{semistable}.  In particular, symmetric
$\Omega$-spectra and (more generally) convergent symmetric spectra are
semistable.  Since we do not include convergence or even semistability
hypotheses, for brevity and clarity we adhere to the following
convention.

\begin{conv}
The \term{homotopy groups} of a symmetric spectrum $X$ will always mean the
homotopy groups of $X$ as an object of 
the stable category, i.e., the abelian groups of maps in the stable
category from $S^q$ to $X$ (for $q\in \bZ$), and we will denote these
as $\pi_{q}X$. 
In the rare cases when we need to refer to the homotopy groups
of the underlying prespectrum of $X$, we will call them the homotopy
groups of the underlying prespectrum, and we introduce no notation for these.
By \term{weak equivalence} of symmetric spectra we shall always mean a
weak equivalence in the stable model structure.  
A weak equivalence is precisely a map that
induces an isomorphism on homotopy groups; it does not necessarily
induce an isomorphism of the homotopy groups of the underlying prespectra.
\end{conv}

Although we do not require convergence hypotheses, they tend to hold
for examples of interest.  In fact, we can replace an arbitrary
small spectral category with a weakly equivalent spectral category that has
the same objects but has mapping spectra that are symmetric
$\Omega$-spectra.  One way of doing this arises 
from the cofibrantly generated Quillen model
category structure on the category of small enriched
categories with a fixed set of objects described in \cite[\S
6]{SSMonoidalEq}.  The maps in this category are the spectral functors
that are the identity on object sets, the fibrations are the maps
$\aC\to \aD$ that restrict to fibrations of symmetric spectra
$\aC(x,y)\to \aD(x,y)$ for all $x,y$ and the weak equivalences are the
maps that restrict to weak equivalences $\aC(x,y)\to \aD(x,y)$ for all
$x,y$.  We use the following terminology from \cite[\S 6]{SSMonoidalEq}. 

\begin{defn}
A small spectral category $\aC$ is said to be \term{pointwise fibrant}
if $\aC(x,y)$ is a fibrant symmetric spectrum (in the stable model
structure) for every pair of objects $x,y$.  Likewise, $\aC$ is said
to be \term{pointwise cofibrant} if $\aC(x,y)$ is a cofibrant
symmetric spectrum for every pair of objects $x,y$.  For a spectral
functor of small spectral categories $F\colon \aC\to \aD$ that is the
identity on the object sets, we say that $F$ is a \term{pointwise weak
equivalence} or \term{pointwise level equivalence} if for every pair
of objects $x,y$, the map $F\colon \aC(x,y)\to \aD(x,y)$ is a weak
equivalence or level equivalence, respectively, of symmetric spectra.
\end{defn}

Fibrant replacement in the model structures of \cite[\S
6]{SSMonoidalEq} then gives most of the following proposition.  The
rest follows from the easy observation that the factorization functors
constructed by the small objects argument on the category of small
spectral categories with a fixed
set of objects still behave well with respect to
spectral functors that are not the identity on object sets.

\begin{prop}\label{propfibrep}\textup{(\cite[6.3]{SSMonoidalEq})}
Given a small spectral category $\aC$, there exists a small spectral
category $\aC^{\Omega}$ and a spectral functor $R\colon \aC\to \aC^{\Omega}$
such that:
\begin{enumerate}
\item $\aC^{\Omega}$ has the same objects as $\aC$ and $R$ is the
identity map on objects,
\item $\aC^{\Omega}$ is pointwise fibrant, and
\item $R$ is a pointwise weak equivalence.
\end{enumerate}
Moreover, $(-)^{\Omega}$ and $R$ may be constructed as an endofunctor and
natural transformation on the category of small spectral categories.
\end{prop}

Applying cofibrant replacement in the model structure of \cite[\S
6]{SSMonoidalEq}, we obtain the following complementary proposition.

\begin{prop}\label{propcofrep}\textup{(\cite[6.3]{SSMonoidalEq})}
Given a small spectral category $\aC$, there exists a small spectral
category $\aC^{\Cell}$ and a spectral functor $Q\colon \aC^{\Cell}\to \aC$
such that:
\begin{enumerate}
\item $\aC^{\Cell}$ has the same objects as $\aC$ and $Q$ is the
identity map on objects,
\item $\aC^{\Cell}(x,y)$ is pointwise cofibrant, and 
\item $Q$ is a pointwise level equivalence.
\end{enumerate}
Moreover, $(-)^{\Cell}$ and $Q$ may be constructed as an endofunctor and
natural transformation on the category of small spectral categories.
\end{prop}

We also use an analogous proposition in the setting of bimodules.

\begin{prop}\label{propmodcof}
Assume that $\aC$ and $\aD$ are pointwise cofibrant small spectral categories.
If $\fM$ is a cofibrant $(\aD,\aC)$-bimodule, then $\fM$ is objectwise
cofibrant, i.e., 
$\fM(d,c)$ is a cofibrant symmetric spectrum for every $(d,c)$ in
$\aD^{\op}\sma \aC$.
\end{prop}

In addition to providing the formal technical results above, the model
theory of enriched categories also explains the relationship of
small spectral categories to small DG-categories.  Sharp statements involve
categories enriched over $H\bZ$-modules (in symmetric spectra of
simplicial sets) or Quillen equivalently, categories enriched over
symmetric spectra of simplicial abelian groups.  For brevity, we will
call these \term{$H\bZ$-categories} and \term{Ab-spectral categories},
respectively.  Note that the category of $H\bZ$-modules is symmetric
monoidal under $\sma_{H\bZ}$ and its derived category is symmetric
monoidally equivalent to the derived category of $\bZ$ (in particular,
$\sma_{H\bZ}$ is more like $\otimes_{\bZ}$ than like $\sma$).  Shipley
\cite[\S 2.2]{ShipleyHZ} produces a zigzag of ``weak monoidal Quillen
equivalences'' relating $H\bZ$-modules to symmetric spectra of
simplicial abelian groups to symmetric spectra of non-negatively
graded and integer graded chain complexes to chain complexes.  A
slight modification gives a comparison for algebras or more generally
categories with a fixed object set $O$: Proposition~6.4 of
\cite{SSMonoidalEq} (or \cite[A.3]{DuggerShipleyEnriched}) gives a
zigzag of Quillen equivalences between the model categories of
DG-categories with object set $O$, Ab-spectral categories with object
set $O$, and $H\bZ$-categories with object set $O$.

\begin{defn}\label{defassoc}
Given a small DG-category, the \term{associated Ab-spectral category
model} or \term{associated $H\bZ$-category model} is the Ab-spectral
category or $H\bZ$-category (respectively) with the same object set
constructed from the zigzag of
Quillen equivalences outlined above using the fibrant or
cofibrant replacement functor (as needed) at every stage.
\end{defn}

By neglect of structure, an Ab-spectral category or $H\bZ$-category is
in particular a spectral category.  We then get the \term{associated
spectral category model} from the associated Ab-spectral or
$H\bZ$-category model.  Because the cofibrant and fibrant
replacement functors in spectral categories with fixed object sets
also behave well with respect to spectral functors that are not the
identity on object sets, the construction of associated spectral
categories in fact produces a functor from small DG-categories to
small spectral categories.

\begin{prop}\label{propDGtoSpec}
The zigzags of Quillen equivalences and cofibrant/fibrant replacement
functors above 
assemble into a functor from the category of small DG-categories to
the category of small spectral categories.
\end{prop}

Note that we do not assert that the associated spectral category
functor provides a 2-functor; it does not preserve composition of
natural transformations.  Using the usual reformulation of natural
transformations and composites of natural transformations as spectral
functors from related spectral categories, we do see that the
associated spectral category functor preserves natural transformations
and their compositions in some coherent homotopy sense.  Since we do
not need this theory here, we leave a rigorous formulation to the
interested reader.

\section{Review of $THH$, $TR$, and $TC$}\label{secthh}

In this section, we review the definition of $THH$, $TR$, and $TC$ of
small spectral categories.  We begin with a review of the cyclic bar
construction for small spectral categories and the variant defined by
B\"okstedt \cite{Bokstedt} and Dundas-McCarthy \cite{DundasMcCarthy}
necessary for the construction of 
$TC$.  We finish with a brief review of the definition of cyclotomic
spectra and the construction of $TR$ and $TC$.

The following cyclic bar construction gives the ``topological''
analogue of the Hoch\-schild-Mitchell complex.

\begin{defn}
For a small spectral category $\aC$ and $(\aC,\aC)$-bimodule $\fM$,
let
\[
\THM_{q}(\aC;\fM)=\bigvee \aC(c_{q-1},c_{q}) \sma \dotsb \sma
\aC(c_{0},c_{1}) \sma \fM(c_{q},c_{0}),
\]
where the sum is over the $(q+1)$-tuples $(c_{0},\dotsc,c_{q})$ of
objects of $\aC$.  This becomes a simplicial object in symmetric
spectra using the usual
cyclic bar construction face and degeneracy maps: The unit
maps of $\aC$ induce the degeneracy maps, and the two action maps on
$\fM$ (for $d_{0}$ and $d_{q}$) and the composition maps in $\aC$ (for
$d_{1},\dotsc,d_{q-1}$) induce the face maps.  We denote the diagonal
(geometric realization) symmetric spectrum as $\THM(\aC;\fM)$ and
write $\THM(\aC)$ for 
$\THM(\aC;\aC)$. 
\end{defn}

The previous construction turns out to be slightly inconvenient to use
as the definition of the topological Hochschild homology of a small spectral
category.  This construction typically only has the correct homotopy
type when the smash products that comprise the terms of the sum
represent the derived smash product.  The analogous problem arises in
the context of Hochschild homology of DG-categories, where the
tensor product may fail to have the right quasi-isomorphism type when
the mapping complexes are not DG-flat.  Just as in that context, this
problem can be overcome using resolutions, such as the ones in
Proposition~\ref{propcofrep} and~\ref{propmodcof}.  There is a further
more subtle difficulty with this construction, however.  While
$\THM(\aC)$ obtains an $S^{1}$-action by virtue of being the geometric
realization of a cyclic complex, the resulting equivariant spectrum
does not have the necessary additional structure to define $TC$
(a well-known problem with this kind of cyclic bar
construction definition of $THH$ in modern categories of spectra).
The correct definition, due to B\"okstedt \cite{Bokstedt} for
symmetric ring spectra and generalized by Dundas-McCarthy
\cite{DundasMcCarthy} to small spectral categories, does not suffer from
either of these deficiencies.

We give a revisionist explanation of the B\"okstedt-Dundas-McCarthy
construction, taking advantage of later results of Shipley
\cite{ShipleyD} on the derived smash product of symmetric spectra.
Let $\aI$ be the category with objects the finite sets
$\bn=\{1,\ldots,n\}$ (including $\bz=\{\}$), and with morphisms the
injective maps.  For a symmetric spectrum $A$, write $A_{n}$ for the
$n$-th space.  The association $\bn\mapsto \Omega^{n}|A_{n}|$ extends
to a functor from $\aI$ to spaces, where $|{-}|$ denotes geometric
realization.  More generally, given symmetric spectra
$A^{0},\dotsc,A^{q}$ and a space $X$, we obtain a functor from
$\aI^{q+1}$ to spaces that sends $\vec\bn=(\bn_{0},\dotsc,\bn_{q})$ to
\[
\Omega^{n_{0}+\dotsb+n_{q}}(|A^{q}_{n_{q}} \sma \dotsb \sma
A^{0}_{n_{0}}|\sma X),
\]
which is also natural in $X$.  Restricting to the case when $X$ is a
sphere $S^{n}$, we form this into a symmetric spectrum following
\cite{ShipleyD} (but using different notation).

\begin{defn}\relax\textup{(\cite[4.2.1]{ShipleyD})}\label{defD}
Let $D(A^{q},\dotsc,A^{0})$ be the symmetric spectrum (of topological
spaces) with $n$-th space
\[
D_{n}(A^{q},\dotsc,A^{0})=\hocolim_{\vec \bn\in\aI^{q+1}}
\Omega^{n_{0}+\dotsb+n_{q}}(|A^{q}_{n_{q}} \sma \dotsb \sma
A^{0}_{n_{0}}|\sma S^{n}),
\]
and the evident structure maps.
\end{defn}

The following is the main lemma of \cite{ShipleyD}.

\begin{prop}\relax\textup{(\cite[4.2.3]{ShipleyD})}
$D(A^{q},\dotsc,A^{0})$ is canonically isomorphic in the stable
category to the derived smash product of the $A^{i}$.
\end{prop}

This motivates the following definition.

\begin{defn}
Given a small spectral category $\aC$, a $(\aC,\aC)$-bimodule $\fM$,
and a space $X$, let $\aG(\aC;\fM;X)_{\vec\bn}$ be the functor from
$\aI^{q+1}$ to spaces defined on $\vec\bn=(\bn_{0},\dotsc,\bn_{q})$ by
\[
\aG(\aC;\fM;X)_{\vec\bn}=
\Omega^{n_{0}+\dotsb+n_{q}}
(\bigvee |\aC(c_{q-1},c_{q})_{n_{q}} \sma \dotsb \sma
\aC(c_{0},c_{1})_{n_{1}}\sma \fM(c_{q},c_{0})_{n_{0}}|\sma X),
\]
and let 
\[
THH_{q}(\aC;\fM)(X)=\hocolim_{\vec\bn\in\aI^{q+1}} \aG(\aC;\fM;X)_{\vec\bn}.
\]
This assembles into a simplicial space, functorially in $X$, as
follows.  The degeneracy maps are induced by
the unit maps $S^{0}\to \aC(c_{i},c_{i})_{0}$ and the functor 
\[
(\bn_{0},\dotsc,\bn_{q})\mapsto (\bn_{0},\dotsc,\bz,\dotsc,\bn_{q})
\]
from $I^{q+1}$ to $I^{q+2}$.  The face maps are induced by the
two action maps on
$\fM$ (for $d_{0}$ and $d_{q}$) and the composition maps in $\aC$ 
(for $d_{1},\dotsc,d_{q-1}$) together with a functor $\aI^{q+1}\to
\aI^{q}$ induced by the appropriate disjoint union isomorphism
$(\bn_{i},\bn_{i+1})\mapsto \bn$ or $(\bn_{q},\bn_{0})\mapsto \bn$ for
$n=n_{i}+n_{i+1}$ or $n=n_{q}+n_{0}$. We write $THH(\aC;\fM)(X)$ for
the geometric realization.
\end{defn}

$THH(\aC;\fM)(X)$ is a continuous functor in the variable $X$, and so
by restriction to the spheres $S^{n}$ specifies a symmetric spectrum
which we denote $THH(\aC;\fM)$ or $THH(\aC)$ for $\fM=\aC$.    The
fact that the symmetric spectrum $THH$ is the 
restriction of a continuous functor implies that it is semistable
\cite[8.7]{MMSS}
and so the object that it represents in the stable category agrees
with its underlying prespectrum.  With additional
hypotheses of ``convergence'' and ``connectivity'', $THH$ is often an
$\Omega$-spectrum; see, for example, Proposition~2.4
of~\cite{HM2}. 

The following propositions, which are essentially the ``many objects''
versions of \cite[4.2.8-9]{ShipleyD} and an easy consequence of the
theory developed in \cite{ShipleyD}, show that in the stable category
$THH(\aC)$ is simply a homotopically well-behaved model
of the Hochschild-Mitchell complex.

\begin{prop}\label{propTHHvsTHM}
There is a natural map in the stable category from $THH(\aC;\fM)$ to
$\THM(\aC;\fM)$ that is an isomorphism when $\aC$ is pointwise
cofibrant.
\end{prop}

\begin{prop}\label{propTHHwe}
Let $F\colon \aC\to \aC'$ be a weak equivalence of small spectral
categories, $\fM'$ a $(\aC',\aC')$-bimodule, $F^{*}\fM'$ the
$(\aC,\aC)$-bimodule obtained by restriction of scalars, and
$\fM\to F^{*}\fM'$ a weak equivalence of $(\aC,\aC)$-bimodules.  Then
the induced map $THH(\aC;\fM)\to THH(\aC';\fM')$ is a weak equivalence.
\end{prop}

As a consequence of the previous propositions, $THH(\aC;\fM)$,
and $THH(\aC)$ always have the correct homotopy type even
when $\THM(\aC;\fM)$ or $\THM(\aC)$ does not.  We also note that
Proposition~\ref{propTHHvsTHM} does not require the bimodule $\fM$ to
be objectwise cofibrant.

We now list the usual bimodule properties of $THH$ that we require in
this paper.  Proofs of these properties appear in the literature
\cite{DundasMcCarthy} under more restrictive hypotheses (i.e.,
connectivity and convergence).

\begin{prop}\label{propnibus}
Let $\aC$ be a small spectral category.
\begin{enumerate}
\item A weak equivalence of $(\aC,\aC)$-bimodules $\fM\to \fM'$ induces a weak
equivalence  $THH(\aC;\fM)\to THH(\aC;\fM')$.
\item A cofibration sequence of $(\aC,\aC)$-bimodules $\fM\to \fM'\to
\fM''\to \Sigma \fM$ induces a homotopy cofiber sequence on $THH$.
\item A fibration sequence of level fibrant $(\aC,\aC)$-bimodules $\Omega
\fM''\to \fM\to \fM'\to \fM''$ induces a homotopy fibration sequence on $THH$.
\end{enumerate}
\end{prop}

\begin{proof}
The first statement is a special case of Proposition~\ref{propTHHwe}.
For the second statement, we can identify
$THH\subdot$ levelwise as the homotopy colimit (over
$\aI^{\ssdot+1}$) of the symmetric spectra
$\aG(\aC;\fM;S^{(-)})_{\vec\bn}$.  The second
statement now follows
from the observation that $\aG$ preserves homotopy cofibration sequences
in the $\fM$ variable and that homotopy colimits and geometric
realization preserve homotopy cofibration sequences.  The third statement
follows from the second since homotopy fibration sequences and homotopy
cofibration sequences agree up to sign.
\end{proof}

We now give a minimal review of the definition of $TR$ and $TC$; we
refer the reader interested in more details to the excellent
discussions of $TR$ and $TC$ in \cite{HM2,HM3}.  For an $S^{1}$-space $X$,
the space $THH(\aC)(X)$ has two $S^{1}$-actions, one coming from $X$
and the other coming from the cyclic structure.  Using the diagonal
action and restricting to representation spheres $S^{V}$ makes
$THH(\aC)(-)$ into an equivariant orthogonal spectrum \cite[\S
II.2]{MMSS}; however, $THH(\aC)$ has even more
structure, that of a \term{cyclotomic spectrum} \cite[\S1.1]{HM3},
\cite[Def.~2.2]{HM2}.  We review the definition of a cyclotomic
spectrum in detail in the next section,
but in brief the structure on $THH$ derives from the
fundamental fixed point map 
\[
(THH(\aC)(X))^{H} \to THH(\aC)(X^{H})
\]
for $S^{1}$-spaces $X$ and finite subgroups $H$ of $S^{1}$.
This induces maps in the equivariant stable category
\[
r_{H}\colon \rho^{\#}_{H}\Phi^{H} THH(\aC)\to THH(\aC)
\]
that are non-equivariant weak equivalences.  Here $\Phi^{H}$ denotes
the (derived) geometric fixed point spectrum, and when $H$ is the
subgroup with $n$ elements, $\rho_{H}$ is the $n$-th root isomorphism
$S^{1}\iso S^{1}/H$; $\rho^{\#}_{H}$ converts the $S^{1}/H$-spectrum
$\Phi^{H}THH(\aC)$ back to an $S^{1}$-spectrum via the isomorphism
$\rho$.  Essentially, a cyclotomic spectrum consists of an
$S^{1}$-equivariant spectrum indexed on a complete universe together
with weak equivalences $r_{H}$ of the form above, called
\term{cyclotomic structure maps}, satisfying certain
coherence properties \cite[Def.~2.2]{HM2}, \cite[\S1.1]{HM3}.  We give
a precise formulation of the point-set category of cyclotomic spectra
we use here in Definition~\ref{defcyclo} in the next section; for now,
the only detail we need is that $THH$ defines a functor from small spectral
categories to the point-set category of cyclotomic spectra
(Theorem~\ref{thmthhcyc}).

For a fixed prime $p$ and each $n$, let $C_{p^{n}}\subset S^{1}$
denote the cyclic subgroup of order $p^{n}$.  We then have maps in the
(non-equivariant) stable category
\[
F,R\colon THH(\aC)^{C_{p^{n}}}\to THH(\aC)^{C_{p^{n-1}}}
\]
where $F$ is the inclusion of the fixed points and $R$ is the map
induced by the composite of the map from the fixed point spectrum to
the geometric fixed point spectrum $THH(\aC)^{C_{p}}\to
\Phi^{C_{p}}THH(\aC)$ and the cyclotomic structure map
$r_{C_{p}}\colon \Phi^{C_{p}}THH(\aC)\to THH(\aC)$; see
\cite[\S1.1]{HM3}, \cite[\S2.2]{HM2}, or Section~\ref{secdefcyc}
below.  We need functorial point-set versions of these maps to
construct $TC$ as a functor on small spectral categories.  In
\cite{HM3}, the connectivity and convergence hypotheses used there
imply that $THH(\aC)$ is an equivariant $\Omega$-spectrum relative to
the family of finite subsets of $S^{1}$; the point-set maps $F,R$ in
\cite{HM3} are then constructed using the point-set fixed point
spectra as models for the derived fixed point spectra.  In our
context, we need to use an $\Omega$-spectrum replacement functor in
the category of cyclotomic spectra; see Definition~\ref{deffibap} and
Theorem~\ref{thmQRF} in the next section.  For such a functor $Q$, we
get appropriate point-set maps
\[
F,R\colon Q(T)^{C_{p^{n}}}\to Q(T)^{C_{p^{n-1}}}.
\]
which are functorial in the cyclotomic spectrum $T$.

\begin{defn}\label{defTC}
Let $Q$ be an $\Omega$-spectrum replacement functor in the category of
cyclotomic spectra and write $T(\aC)$ for $Q(THH(\aC))$. Then 
$TR\supdot(\aC)$ is the pro-spectrum $\{T(\aC)^{C_{p^{n}}}\}$ under the
maps $R$, and $TR(\aC)$ is the homotopy limit.  $TC(\aC)$ and
$TC\supdot(\aC)$ are the spectrum and 
pro-spectrum obtained from $TR(\aC)$ and $TR\supdot(\aC)$ as the homotopy
equalizer of the maps $F$ and $R$. 
\end{defn}

Note that a map in the $S^{1}$-equivariant stable category induces a
(non-equivariant) weak equivalence on fixed point spectra for all
finite subgroups of $S^{1}$ if and only if it induces a
(non-equivariant) weak equivalence on geometric fixed point spectra
for all finite subgroups \cite[XVI.6.4]{MayAlaska}.  It follows that a
cyclotomic map of cyclotomic spectra induces a weak equivalence of
fixed point spectra for all finite subgroups of $S^{1}$ if and only if
it is a non-equivariant weak equivalence.  In particular, we obtain
the following proposition.

\begin{prop}\label{propTCequiv}
A spectral functor of small spectral categories $\aC\to \aD$ that induces a weak
equivalence on $THH$ induces a weak equivalence on $TR$ and $TC$.
\end{prop}

Likewise, using the same principle on the cofiber of a map of
cyclotomic spectra, we obtain the following proposition.  Applying
this proposition in examples when $THH(\aC)$ is contractible,
localization cofibration sequences on $TR$ and $TC$ follow from ones on
$THH$.

\begin{prop}\label{proplazy}
For a strictly commuting square of small spectral categories
\[
\xymatrix{%
\aA\ar[r]\ar[d]&\aB\ar[d]\\
\aC\ar[r]&\aD,
}
\]
if the induced square on $THH$ is homotopy cocartesian, then so are the
induced squares on $TR$ and $TC$.
\end{prop}

Finally, we turn to DG-categories.  For a DG-category $\aC^{DG}$, we
can consider $THH$ of the associated spectral category $\aC^{S}$.
Propositions~\ref{propDGtoSpec}
and~\ref{propTCequiv} show that defining $THH$, $TR$, and $TC$ of
$\aC^{DG}$ in terms of $THH(\aC^{S})(-)$ constructs $THH$, $TR$, and
$TC$ as functors from the category of DG-categories and DG-functors to
the stable category; this is Theorem~\ref{intthhfunct}.

\section{Details of the category of cyclotomic spectra}
\label{secdefcyc}

Starting in Section~\ref{secgenloc}, we will want a construction of
$THH$ that provides a point-set functor from the category of small
spectral categories to a point-set category of cyclotomic spectra.  We
describe such a category and functor in this section.  

In contrast to the canonical definition of cyclotomic spectra
in \cite{HM2} (which is in terms of Lewis-May spectra), we use
orthogonal spectra.  To make this section easier to compare with the
definition in \cite{HM2}, we follow the notation there as much as
possible (but with some abbreviation).  The following summarizes the
notation and terminology regarding the circle group, its subgroups,
and its representations.

\begin{notn}\label{notHM2}
Let $G$ denote the circle group of unit complex numbers.  For
$n=1,2,3,\dotsc$, let $C_{n}$ (or just $C$) denotes the subgroup with
$n$ elements.  Write $\rho_{n}$ for the $n$-th root isomorphism $G\to
G/C_{n}$ and for a $G/C_{n}$-space $X$ (e.g., $X=Y^{C_{n}}$ for some
$G$-space $Y$), write $\rho^{*}_{n}X$ for the $G$-space obtained via
this isomorphism.  Throughout this section \term{orthogonal $G$-representation}
will mean a finite dimensional real $G$-inner product space. Let
$\bC(0)=\bC$ denote the complex numbers with 
trivial $G$-action, $\bC(1)$ the complex numbers with the standard
$G$-action, and $\bC(n)=\bC(1)\otimes_{\bC}\dotsb \otimes_{\bC}\bC(1)$
($n$-factors).  Let
\[
\aU=\bigoplus_{n=0}^{\infty}\bigoplus_{r=1}^{\infty}\bC(n)
\]
(a direct sum of infinitely many copies of each $\bC(n)$) with the
standard inner product; this is a \term{complete $G$-universe}, a real
countable-dimensional $G$-inner product space containing an isomorphic
copy of every orthogonal $G$-representation. When we write $V<\aU$, we
will always understand $V$ to be a finite dimensional $G$-invariant
vector subspace (and thus, an orthogonal $G$-representation).  For $V<
W <\aU$, let $W-V$ denote the orthogonal complement of $V$ in $W$
\end{notn}

We now give a very brief review of the definition of orthogonal
$G$-spectra; see 
\cite[\S\S II.2,II.4]{MM} for a complete treatment.  For orthogonal
$G$-representations $V$ and $W$, let $\sJ_{G}(V,W)$ denote the Thom
space of the orthogonal complement $G$-bundle of (non-equivariant)
linear isometries from $V$ to $W$ \cite[II.4.1]{MM} (an element of the
complement $G$-bundle consists of a linear isometry $V\to W$ and a
point in the orthogonal complement of the image).  Composition of
isometries and addition in the codomain vector space induces
composition maps
\[
\sJ_{G}(W,Z)\sma \sJ_{G}(V,W)\to \sJ_{G}(V,Z),
\]
which together with the obvious identity elements
make $\sJ_{G}$ a category enriched in based $G$-spaces (with objects
the orthogonal $G$-representations).  An orthogonal $G$-spectrum is a
$\sJ_{G}$-space \cite[II.4.3]{MM},
an enriched functor from (a skeleton of) $\sJ_{G}$ to
based $G$-spaces.  That is, an orthogonal $G$-spectrum $T$ consists of
a based $G$-space $T(V)$ for each orthogonal $G$-representation $V$
together with based $G$-maps 
\[
\sJ_{G}(V,W) \sma T(V) \to T(W)
\]
satisfying the obvious unit and associativity properties. 
We recall that the homotopy groups of $T$ are defined by
\[
\pi_{q}^{H}T = 
\begin{cases}
\quad\displaystyle 
\lcolim_{V< \aU} \pi_{q}((\Omega^{V}T(V))^{H})&q\geq 0\\[1em]
\quad\displaystyle 
\lcolim_{\bC^{-q}<V< \aU} \pi_{-q}((\Omega^{V-\bC^{-q}}T(V))^{H})&q< 0\\
\end{cases}
\]
for $q\in \bZ$ and $H$ a closed subgroup of $G$ \cite[III.3.2]{MM}.
We also have a similar formula for the homotopy groups of the
geometric fixed point spectrum of the underlying object in the
equivariant stable category; these groups were denoted as
$\rho^{H}_{q}$ in \cite[\S V.4]{MM}, but we will denote them as
$\pi^{\Phi H}_{q}$ to avoid confusion with the isomorphisms $\rho_{n}$
in Notation~\ref{notHM2}.  Then
\[
\pi_{q}^{\Phi H}T = 
\begin{cases}
\quad\displaystyle 
\lcolim_{V< \aU} \pi_{q}(\Omega^{V^{H}}(T(V)^{H}))&q\geq 0\\[1em]
\quad\displaystyle 
\lcolim_{\bC^{-q}<V< \aU} \pi_{-q}(\Omega^{V^{H}-\bC^{-q}}(T(V)^{H}))&q< 0\\
\end{cases}
\]
for $q\in \bZ$ and $H$ a closed subgroup of $G$
\cite[V.4.8.(iii),V.4.12]{MM}.  

Throughout this paper, we use the following precise definition for
cyclotomic spectra.

\begin{defn}\label{defcyclo}
A \term{cyclotomic spectrum} is an orthogonal
$G$-spectrum $T$ together with $G$-equivariant maps 
\[
r_{n,V}\colon \rho^{*}_{n}(T(V)^{C_{n}}) \to T(\rho^{*}_{n}(V^{C_{n}})),
\]
for $V$ an orthogonal $G$-representation
and $n=1,2,3,\dotsc$, satisfying the following conditions.
\begin{enumerate}
\item $r_{1,V}$ is the identity for all $V$.

\item For any $V,W,n$, the diagram
\[
\xymatrix{%
\rho^{*}_{n}(\sJ_{G}(V,W)^{C_{n}})\sma \rho^{*}_{n}(T(V)^{C_{n}})
\ar[r]^-{\phi \sma r_{n,V}}\ar[d]
&\sJ_{G}(\rho^{*}_{n}(V^{C_{n}}),\rho^{*}_{n}(W^{C_{n}}))\sma 
T(\rho^{*}_{n}(V^{C_{n}}))\ar[d]\\
\rho^{*}_{n}(T(W)^{C_{n}})\ar[r]_{r_{n,W}}&T(\rho^{*}_{n}(W^{C_{n}}))
}
\]
commutes. Here 
vertical maps are induced by the $\sJ_{G}$-space structure maps.
In the top horizontal map, the map 
\[
\phi \colon \rho^{*}_{n}(\sJ_{G}(V,W)^{C_{n}})\to
\rho^{*}_{n}\sJ_{G/C_{n}}(V^{C_{n}},W^{C_{n}})
=\sJ_{G}(\rho^{*}_{n}(V^{C_{n}}),\rho^{*}_{n}(W^{C_{n}}))
\]
is induced by the fixed point functor $\phi \colon
\sJ_{G}(V,W)^{C_{n}}\to \sJ_{G/C_{n}}(V^{C_{n}},W^{C_{n}})$ described
in \cite[\S V.4]{MMSS}.

\item For all $V,m,n$, the diagram
\[
\xymatrix@C+1pc{%
\rho^{*}_{n}((\rho^{*}_{m}(T(V)^{C_{m}}))^{C_{n}})
\ar[r]^-{\rho^{*}_{n}(r_{m,V}{}^{C_{n}})}
\ar[d]_-{=}
&\rho^{*}_{n}(T(\rho^{*}_{m}(V^{C_{m}}))^{C_{n}})
\ar[d]^-{r_{n,\rho^{*}_{m}(V^{C_{m}})}}\\
\rho^{*}_{mn}(T(V)^{C_{mn}})\ar[r]_-{r_{mn,V}}
&T(\rho^{*}_{mn}(V^{C_{mn}}))
}
\]
commutes (using 
$\rho^{*}_{mn}(X^{C_{mn}})=\rho^{*}_{n}((\rho^{*}_{m}(X^{C_{m}}))^{C_{n}})$
for any $G$-space $X$).

\item For all $q\in \bZ$ and all $n$, the map
\begin{align*}
\lcolim_{V< \aU} \pi_{q}\Omega^{V^{C_{n}}}(T(V)^{C_{n}})
&\to \lcolim_{W< \aU^{C_{n}}} \pi_{q}\Omega^{W}(T(\rho^{*}_{n}W)),
&\text{if }q\geq 0\\
\lcolim_{\bC^{-q}< V< \aU} 
  \pi_{-q}\Omega^{V^{C_{n}}-\bC^{-q}}(T(V)^{C_{n}})
&\to \lcolim_{\bC^{-q}<W<\aU^{C_{n}}} 
  \pi_{-q}\Omega^{W-\bC^{-q}}(T(\rho^{*}_{n}W)),
&\text{if }q< 0
\end{align*}
is an isomorphism. 
\end{enumerate}
A map of cyclotomic spectra $(T,r)\to (T',r')$ is a map of orthogonal
spectra $f\colon T\to T'$ making the diagram
\[
\xymatrix{%
\rho^{*}_{n}T(V)^{C_{n}}\ar[r]^{r_{n}}\ar[d]_{\rho^{*}_{n}f^{C_{n}}}
&T(\rho^{*}_{n}V^{C_{n}})\ar[d]^{f}\\
\rho^{*}_{n}T'(V)^{C_{n}}\ar[r]_{r'_{n}}
&T'(\rho^{*}_{n}V^{C_{n}})
}
\]
commute for all $n$, $V$.
\end{defn}

The explanation of~(iv) is that the colimit on the left is
$\pi^{\Phi C_{n}}_{q}T=\pi_{q}\Phi^{C_{n}}T$ and the colimit on the
right is $\pi_{q}T$.  Condition~(iv) then insures that the induced
map in the stable category $\Phi^{C_{n}}T\to T$ constructed as in
\cite[Lem.~2.2]{HM2} is an isomorphism.  We can rephrase the previous
definition in terms of the point-set model $\Fix{n}T$ for the
geometric fixed point spectrum $\Phi^{C_{n}}T$ from \cite[\S V.4]{MM}.

\begin{rem}
In \cite[\S V.4]{MM}, a version of the geometric fixed point spectrum
called $\Fix{n}T$ is constructed as a $\sJ_{E_{n}}$-space (for $E_{n}$
the extension $C_{n}\rightarrowtail G\twoheadrightarrow G/C_{n}$),
where $\sJ_{E_{n}}=(\sJ_{G})^{C_{n}}$.  We write $i_{n}$ for the
inclusion of $\sJ_{E_{n}}$ into $\sJ_{G}$ and $\phi_{n} \colon
\sJ_{E_{n}}\to \sJ_{G/C_{n}}$ for the functor that sends $V$ to
$V^{C_{n}}$; then $i_{n}$ is an enriched functor of categories enriched in
based $G$-sets and $\phi_{n}$ is an enriched functor of categories
enriched in based $G/C_{n}$-sets.  The functor $\phi_{n}$ allows us to
regard an orthogonal $G/C_{n}$-spectrum as a $\sJ_{E_{n}}$-space, and
we will use $\rho^{*}_{n}\phi_{n}$ to regard a
$\rho^{*}_{n}\sJ_{G/C_{n}}$-space as a  
$\rho^{*}_{n}\sJ_{E_{n}}$-space.  Of course, $\rho_{n}$ induces an
isomorphism of enriched categories $\rho^{*}_{n}\sJ_{G/C_{n}}\to
\sJ_{G}$, and we write
$\phi^{!}_{n}$ for the composite enriched functor
$\rho^{*}_{n}\sJ_{E_{n}}\to \sJ_{G}$.
Recall that $\Fix{n}T$ is the
functor $(T\circ i_{n})^{C_{n}}$.  
Definition~\ref{defcyclo} is equivalent to asking for a map of
$\rho^{*}_{n}\sJ_{E_{n}}$-spaces
\[
r_{n}\colon \rho^{*}_{n}\Fix{n}T \to T\circ \phi^{!}_{n}
\]
for each $n$ that is a non-equivariant weak equivalence, is the
identity for $n=1$, and that makes the following diagram of
$\rho^{*}_{mn}\sJ_{E_{mn}}$-spaces commute for all $mn$:
\[
\xymatrix@C-1.375pc@R-.5pc{%
&\rho^{*}_{mn}\Fix{mn} T\ar[r]^{r_{mn}}
&T\circ \phi^{!}_{mn}\ar[dr]^{=}\\
\rho^{*}_{n}(\rho^{*}_{m}(\Fix{m}T)\circ \rho^{*}_{m}i_{mn}^{n})^{C_{n}}
\ar[dr]_{r_{m}}\ar[ur]^{=}\hspace{-5pc}
&&&\hspace{-3pc}
T\circ \phi^{!}_{n}\circ \rho^{*}_{n}\phi^{!}_{m;mn}\\
&\rho^{*}_{n}((T\circ \phi^{!}_{m})\circ \rho^{*}_{m}i_{mn}^{n})^{C_{n}}\ar@{{}{=}{}}[r]
&\rho^{*}_{n}(\Fix{n}T)\circ \rho^{*}_{n}\phi^{!}_{m;mn}\ar[ur]_{r_{n}}
}
\]
Here the maps labelled $r_{m}$ and $r_{n}$ are really the maps induced
by $r_{m}$ and $r_{n}$ inside the fixed points and compositions.  The
functor labelled $i_{mn}^{n}$ is the inclusion of
$\sJ_{E_{mn}}=\sJ_{G}^{C_{mn}}$ into $\sJ_{E_{m}}=\sJ_{G}^{C_{m}}$ and
the functor labelled $\phi^{!}_{m;mn}$ is the map
$\rho^{*}_{m}\sJ_{E_{mn}}\to \sJ_{E_{n}}$ obtained from the restriction
of $\phi^{!}_{m}$.  The functor $\rho^{*}_{n}\phi^{!}_{m;mn}$ is then
the induced functor 
\[
\rho^{*}_{n}(\rho^{*}_{m}\sJ_{E_{mn}})\to \rho^{*}_{n}\sJ_{E_{n}}.
\]
The diagonal equality on the left arises from the equation
\[
i_{mn}=i_{m}\circ i_{mn}^{n},
\]
the horizontal equality on the bottom
arises from the equation
\[
\phi^{!}_{m}\circ \rho^{*}_{m}i_{mn}^{n}=i_{n}\circ \phi^{!}_{m;mn},
\]
and the diagonal equality on the right arises from the equation
\[
\phi^{!}_{mn}=\phi^{!}_{n}\circ \rho^{*}_{n}\phi^{!}_{m;mn}.
\]
In this case, the spacewise definition above seems less complicated
than the equivalent diagram space definition here. 
\end{rem}

Although we cannot expect the category of cyclotomic spectra to admit
general limits and colimits (because of condition~(iv)), it is closed
under homotopy colimits and finite homotopy limits in orthogonal
spectra.

\begin{prop}
If $\Psi$ is a functor from a small category $\aD$ to cyclotomic spectra,
then the homotopy colimit of $\Psi$ in orthogonal $G$-spectra admits the
natural structure of a cyclotomic spectrum.  If $\aD$ has only
finitely many objects and finitely many sequences of composable
non-identity morphisms, then the homotopy limit of $\Psi$ in
orthogonal $G$-spectra likewise admits the natural structure of a
cyclotomic spectrum.
\end{prop}

\begin{proof}
Under the hypotheses, we can compute the homotopy colimit and homotopy
limit spacewise using the usual bar and cobar constructions, respectively.
The structure maps $r_{n,V}$ then commute with these constructions
and conditions (i)--(iii) in the definition of cyclotomic spectra are
clearly preserved.  The standard properties of homotopy groups and
geometric fixed points in the equivariant stable category imply that
condition~(iv) is also preserved.
\end{proof}

In this section, we write $T^{C_{n}}$ for the point-set fixed point
orthogonal spectrum of an orthogonal $G$-spectrum $T$: It is a
non-equivariant orthogonal spectrum with $T^{C_{n}}(V)=(T(V))^{C_{n}}$
where we regard a (non-equivariant) inner product space $V$ as a
$G$-inner product space with trivial $G$-action.  The inclusion of
$C_{m}$ in $C_{mn}$ then induces a map of (non-equivariant) orthogonal
spectra 
\[
F_{n}\colon T^{C_{mn}}\to T^{C_{m}}.
\]
This system of maps is compatible in the sense that 
\[
F_{p}\circ F_{n}=F_{np}\colon T^{C_{mnp}}\to T^{C_{m}}
\]
for all $m$, $n$, $p$.  Now let $T$ be a cyclotomic spectrum.
Restricting to the inner
product spaces $V$ with trivial $G$-action, the maps $r_{n,V}$
assemble to a map of (non-equivariant) orthogonal spectra
$T^{C_{n}}\to T$; more generally, looking at the maps 
\[
r_{n,V}{}^{C_{m}}\colon T(V)^{C_{mn}}=(\rho^{*}_{n}T(V)^{C_{n}})^{C_{m}}
\to T(V)^{C_{m}},
\]
we get maps of (non-equivariant) orthogonal spectra
\[
R_{n}\colon T^{C_{mn}}\to T^{C_{m}}.
\]
This system of maps is compatible in the sense that 
\[
R_{p}\circ R_{n}=R_{np}\colon T^{C_{mnp}}\to T^{C_{m}}
\]
for all $m$, $n$, $p$.  Moreover, by construction, the maps $R_{n}$
and $F_{p}$ commute,
\[
F_{n}\circ R_{p}=R_{p}\circ F_{n}\colon T^{C_{mnp}}\to T^{C_{m}}.
\]
For $p$ a fixed prime $F_{p}$ and $R_{p}$ provide the maps $F$ and $R$
for constructing $TR$ and $TC$.  To get the correct homotopy type, we
need the fixed point orthogonal spectra to have the correct homotopy
type, and we can only expect this to happen when $T$ satisfies
additional constraints, such as being an equivariant
$\Omega$-spectrum.  To arrange this, we use an $\Omega$-spectrum
replacement functor.

\begin{defn}\label{deffibap}
An $\Omega$-spectrum replacement functor in the category of cyclotomic
spectra consists of a functor $Q$ from cyclotomic spectra to itself and a
natural transformation of cyclotomic spectra $\theta \colon \Id\to Q$ such that
for any cyclotomic spectrum $T$, the map $\theta \colon T\to QT$ is a
weak equivalence of orthogonal $G$-spectra and $QT$ is a fibrant
orthogonal $G$-spectrum, i.e., for any orthogonal $G$-representations
$V,W$, the structure map $QT(V)\to \Omega^{W}QT(V\oplus W)$ is a weak
equivalence of $G$-spaces.
\end{defn}

As indicated in the previous section, we then construct $TR$ and $TC$
for a cyclotomic spectrum $T$ using the maps $F$ and $R$ for the
cyclotomic spectrum $QT$.  For this to work, we need to know that an
$\Omega$-spectrum replacement functor in the category of cyclotomic
spectra exists.  We describe two, $Q^{\aN}$ and $Q^{\aI}$, both of
which are lifts of fibrant replacement functors on the category of
orthogonal $G$-spectra.  Recall that $\aI$ is the category whose
objects are the finite sets $\bm=\{1,\dotsc,m\}$ (including
$\bz=\{\}$) and whose maps are the inclusions.  Let $\aN$ be the
subset of standard inclusions, so that there is a unique map $\bm\to
\bn$ whenever $m\leq n$.  Let $\aI^{\omega}$ be the subcategory of
$\aI^{\infty}=\prod_{i=0}^{\infty}\aI$ of sequences
$\vec\bm=(\bm_{0},\bm_{1},\dotsc )$ such that all but finitely many of
the $\bm_{i}$ are $\bz$, and let $\aN^{\omega}$ be the corresponding
subcategory of $\aN^{\infty}$.  For $\vec\bm$ in $\aI^{\omega}$, let
\[
\bC(\vec\bm)=\bC(0)^{m_{0}}\oplus \bC(1)^{m_{1}}\oplus
\bC(2)^{m_{2}}\oplus \dotsb <\aU.
\]
Then $Q^{\aN}T$ and $Q^{\aI}T$ are the orthogonal $G$-spectra with
$V$-th spaces
\begin{align*}
Q^{\aN}T(V)&=\hocolim_{\vec\bm\in \aN} 
 \Omega^{\bC(\vec \bm)}(T(V\oplus \bC(\vec \bm)))\\
Q^{\aI}T(V)&=\hocolim_{\vec\bm\in \aI} 
 \Omega^{\bC(\vec \bm)}(T(V\oplus \bC(\vec \bm)))
\end{align*}
with structure maps induced by the $V$ variable.  The inclusion of the
object 
$\vec\bz=(\bz,\bz,\dotsc )$ in the homotopy colimit then induces a map
$T(V)\to Q^{\aN}T(V)\to 
Q^{\aI}T(V)$ that induces the natural transformations $T\to Q^{\aN}T$
and $T\to Q^{\aI}T$. We then have the
following well-known fact.

\begin{prop}
The functors $Q^{\aN}$ and $Q^{\aI}$ are fibrant replacement functors
in orthogonal $G$-spectra.
\end{prop}

\begin{proof}
This is clear for $Q^{\aN}$ by calculating homotopy groups.  The
space-level fixed point functors commute with the homotopy colimits
and the obvious generalization of the argument of
\cite[2.2.9]{ShipleyD} to $\aI^{\omega}$ shows that the map
$(Q^{\aN}T(V))^{C_{n}}\to(Q^{\aI}T(V))^{C_{n}}$ is a weak equivalence
for all $V$ and $n$.
\end{proof}

Let $\gamma_{n}\colon \aI^{\omega}\to \aI^{\omega}$ be the functor
that sends $\vec\bm\in \aI$ to
\[
\gamma_{n}(\vec \bm)=(\bm_{0},\bm_{n},\bm_{2n},\dotsc);
\]
Then $\bC(\gamma_{n}(\vec \bm))=\rho^{*}_{n}(\bC(\vec\bm)^{C_{n}})$,
and we get maps $r_{n,V}$ for $Q^{\aN}T$ and $Q^{\aI}T$ as the map
induced by the functor $\gamma_{n}$ and the natural transformation
\begin{multline*}
\rho^{*}_{n}((\Omega^{\bC(\vec \bm)}(T(V\oplus \bC(\vec \bm))))^{C_{n}})
\to
\rho^{*}_{n}(\Omega^{\bC(\vec \bm)^{C_{n}}}(T(V\oplus \bC(\vec \bm))^{C_{n}}))\\
=
\Omega^{\bC(\gamma_{n}(\vec \bm))}\rho^{*}_{n}(T(V\oplus \bC(\vec \bm))^{C_{n}})
\xrightarrow{r_{n}}
\Omega^{\bC(\gamma_{n}(\vec \bm))}
T(\rho^{*}_{n}V^{C_{n}}\oplus \bC(\gamma_{n}(\vec \bm))).
\end{multline*}
A check of the diagrams then proves the following theorem.

\begin{thm}\label{thmQRF}
The natural transformations above make $Q^{\aN}$ and $Q^{\aI}$ into
$\Omega$-spectrum replacement functors in the category of cyclotomic
spectra. 
\end{thm}

We also make the (trivial) observation for use in
Section~\ref{apptrace} that as endofunctors on orthogonal $G$-spectra,
both $Q^{\aN}$ and $Q^{\aI}$ are enriched in based spaces (i.e., they
are continuous and preserve the base point on mapping spaces).

The advantage of $Q^{\aN}$ is that it is very closely related to the
Lewis-May spectrification functor for ``good'' prespectra used in
\cite{HM2} to construct the $TR$ and $TC$ of a cyclotomic spectrum;
the proof of the following proposition is easy, but is omitted as it
would require a long review of \cite[\S2.2]{HM2}.

\begin{prop}
For any cyclotomic spectrum $T$, the underlying Lewis-May prespectra
of $TR$ and $TC$ constructed using $Q^{\aN}T$ above is weakly
equivalent through a natural zigzag to the underlying prespectrum of
$TR$ and $TC$, respectively, constructed from the underlying
equivariant Lewis-May 
prespectrum of $T$ using the construction of \cite{HM2}.
\end{prop}

The advantage of $Q^{\aI}$ is that viewed as a fibrant replacement
functor on orthogonal $G$-spectra, it is lax monoidal.  A famous
observation of Lewis \cite{LewisNotConvenient} (in modern terms) is
that no fibrant replacement functor in orthogonal $G$-spectra can be
symmetric monoidal.  There are symmetric monoidal fibrant replacement
functors for the positive stable model structure \cite{Kro} (and this
would suffice to get homotopically correct fixed point functors), but
we do not know how to construct one that provides a positive
$\Omega$-spectrum replacement functor in the category of cyclotomic
spectra.

Finally, we need to check that for a small spectral category $\aC$,
$THH(\aC)$ is naturally a cyclotomic spectrum.  The remainder of the
section is devoted to the proof of the following theorem.

\begin{thm}\label{thmthhcyc}
$THH(-)$, with the natural structure maps $r_{n,V}$ described below, is a
functor from the category of small spectral categories to the category of
cyclotomic spectra. 
\end{thm}

\begin{proof}
Recall from the previous section that $THH(\aC)$ is defined as the
geometric realization of the cyclic object
\[
THH_{q}(\aC)(V)=\hocolim_{\vec\bn\in\aI^{q+1}} \aG(\aC;S^{V})_{\vec\bn}.
\]
where
\[
\aG(\aC;S^{V})_{\vec\bn}=
\Omega^{n_{0}+\dotsb+n_{q}}
(\bigvee |\aC(c_{q-1},c_{q})_{n_{q}} \sma \dotsb \sma
\aC(c_{0},c_{1})_{n_{1}}\sma \aC(c_{q},c_{0})_{n_{0}}|\sma S^{V}).
\]
Here $V$ is meant to be an orthogonal $G$-representation, and the
$G$-action on $THH(\aC)(V)$ is the diagonal action of the $G$-action
on the cyclic structure and the $G$-action on $S^{V}$.
Restricting to the subgroup $C_{n}$, the action then has a concrete
description in terms of the $n$-th edgewise subdivision
$\sd_{n}THH(\aC)(V)$.  This is the geometric
realization of the simplicial $C_{n}$-space
\[
\sd_{n}THH\subdot(\aC)(V)=
THH_{n(\ssdot+1)-1}(\aC)(V)=
\hocolim_{\vec\bp\in\aI^{n(\ssdot+1)}} \aG(\aC;S^{V})_{\vec\bp},
\]
where the $C_{n}$-action is induced both by rotating the coordinates
of $\aI^{n(\ssdot+1)}$ (and the corresponding loops and factors of
$\aC$ inside $\aG(\aC;S^{V})$) and by the action on $S^{V}$.  An
element of the homotopy colimit can only be a $C_{n}$-fixed point
when it comes from a $\vec\bp$ of the form
\[
n\vec\bm =
(\bm_{0},\dotsc,\bm_{q},\bm_{0},\dotsc,\bm_{q},\dotsc,\bm_{0},\dotsc,\bm_{q})
\]
for a sequence $\vec\bm=(\bm_{0},\dotsc,\bm_{q})$ repeated $n$ times.
For such a $\vec\bp$, 
\begin{multline*}
\aG(\aC;S^{V})_{\vec\bp}=\Omega^{n(m_{0}+\dotsb+m_{q})}
(\bigvee |\aC(c_{n(q+1)-2},c_{n(q+1)-1})_{m_{q}} \sma \dotsb \\
\sma
\aC(c_{(n(q+1)-q-1},c_{n(q+1)-q})_{m_{0}}\sma
\aC(c_{(n-1)(q+1)-2},c_{(n-1)(q+1)-1})_{m_{q}}\sma \dotsb \\
\sma
\aC(c_{0},c_{1})_{m_{1}}\sma \aC(c_{n(q+1)-1},c_{0})_{m_{0}}|\sma S^{V})
\end{multline*}
and the $C_{n}$-action on the homotopy colimit restricts to a
$C_{n}$-action on $\aC(\aC;S^{V})_{n\vec \bm}$.
Viewing $\Omega^{nm}$ as based maps out of $S^{nm}$, the $C_{n}$ fixed
points are the $C_{n}$-equivariant maps out of $S^{nm}$; for such a
map, restricting to fixed points gives a based map from
$S^{m}=(S^{nm})^{C_{n}}$ to 
\begin{multline*}
(\bigvee
|\aC(c_{n(q+1)-2},c_{n(q+1)-1})_{m_{q}} \sma \dotsb \sma
\aC(c_{n(q+1)-1},c_{0})_{m_{0}}|\sma S^{V})^{C_{n}}\\
=\bigvee |\aC(c_{q-1},c_{q})_{m_{q}}\sma \dotsb \sma
\aC(c_{0},c_{q})_{m_{0}}|\sma S^{V^{C_{n}}}.
\end{multline*}
Thus, restricting to fixed points induces a map
\[
\aG(\aC;S^{V})_{n\vec\bm}^{C_{n}}\to 
\aG(\aC;S^{V^{C_{n}}})_{\vec\bm}.
\]
We define the map $r_{n,V}$ to be the induced simplicial map
\begin{multline}\label{eqdefr}
(\hocolim_{\bp\in \aI^{n(\ssdot+1)}}\aG(\aC;S^{V})_{\bp})^{C_{n}}
\iso
\hocolim_{\bm\in \aI^{\ssdot+1}}\aG(\aC;S^{V})_{n\bm}^{C_{n}}\\
\to
\hocolim_{\bm\in \aI^{\ssdot+1}}\aG(\aC;S^{V^{C_{n}}})_{\bm}.
\end{multline}
Following through the resulting $G$-action, we get as required a map
\[
r_{n,V}\colon \rho^{*}_{n}(THH(\aC,V)^{C_{n}})\to
THH(\aC,\rho^{*}_{n}(V^{C_{n}})),
\]
natural in maps of the small spectral category $\aC$.
Conditions~(i), (ii), and~(iii) of the definition of cyclotomic
spectra are straightforward checks of the diagrams.  It remains to check
condition~(iv). 

For any $W$ in $\aU^{C_{n}}$, let $\bar W$ be the union of the
$G$-subspaces $V$ of $\aU$ with 
$V^{C_{n}}=W$.  Then $\bar W$ is an infinite dimensional subspace of
$\aU$ but we can generalize the notation $THH(\aC)(V)$ to obtain the $G$-space
$THH(\aC)(\bar W)$ and the map $r_{n,V}$ generalizes to a map
\[
r_{n,\bar W}\colon \rho^{*}_{n}(THH(\aC,\bar W)^{C_{n}})\to
THH(\aC,\rho^{*}_{n}(\bar W^{C_{n}}))
=THH(\aC,\rho^{*}_{n}W).
\]  
We then have
\[
THH(\aC)(\bar W)=\lcolim_{V<\aU, V^{C_{n}}=W} THH(\aC)(V),
\]
and we get a canonical isomorphism 
\[
\lcolim_{V< \aU} \pi_{q}\Omega^{V^{C_{n}}}(THH(\aC)(V)^{C_{n}})
\iso
\lcolim_{W<\aU^{C_{n}}} \pi_{q}\Omega^{W}((THH(\aC)(\bar W))^{C_{n}})
\]
for $q\geq 0$ and likewise for the colimit for negative $q$. Thus, it
suffices to show that the map $r_{n,\bar W}$ is a non-equivariant weak
equivalence. 

The space $S^{\bar W}$ is a model for the space $\Sigma^{W}\tilde
E\aF[C_{n}]$, meaning that for $H$ a closed subgroup of $G$, the fixed
point space 
$(S^{\bar W})^{H}$ is $(S^{W})^{H}$ if $H$ contains $C_{n}$
and is contractible otherwise.  It follows that for any based CW
$C_{n}$-spaces $Y$, $Z$, the map
\[
F(Y,Z\sma S^{\bar W})^{C_{n}}\to F(Y^{C_{n}}, Z^{C_{n}}\sma S^{W})
\]
is a weak-equivalence, where $F$ denotes the
$C_{n}$-space of 
based maps (cf. \cite[II.9.3]{LMS}).  In particular, since the map 
$r_{n,\bar W}$ is up to isomorphism the geometric realization of a
homotopy colimit of maps of this form~\eqref{eqdefr}, it is 
a weak equivalence, and 
this completes the argument.
\end{proof}

\section[Invariance of $THH$]{Spectral categories, homotopy categories,
and invariance of $THH$}\label{secdkequiv}

In this section, we continue the discussion of the basic properties of
spectral categories and study the natural conditions under which
spectral functors induce equivalences on $THH$.  We review the concept
of ``Dwyer-Kan equivalence'' (Definition~\ref{defdkequiv}) of spectral
categories, which provides a more sophisticated notion of weak
equivalence of spectral categories; Theorem~\ref{thmdkequiv}
below indicates that Dwyer-Kan equivalences induce equivalences of
$THH$.  The mapping spectra of a spectral category $\aC$ give rise to
an associated ``homotopy category'' that is an invariant of the
Dwyer-Kan equivalences.  Under rather general conditions
(q.v. Definition~\ref{deftriang}), the homotopy category has a
triangulated structure and this allows us to formulate useful
``cofinality'' and ``thick subcategory'' criteria for spectral
functors to induce equivalences of $THH$ in Theorems~\ref{thmfactor}
and~\ref{thmthick}.  Proofs of Theorems~\ref{thmdkequiv},
\ref{thmfactor}, and~\ref{thmthick} require the technical tools
developed in the next section and are given there.

\begin{defn}\label{defdkequiv}
Let $F\colon \aC\to \aD$ be a spectral functor.  We say that $F$ is a
\term{Dwyer-Kan embedding} or \term{DK-embedding} when for every
$a,b\in \ob\aC$, the map $\aC(a,b)\to \aD(Fa,Fb)$ is a weak
equivalence.

We say that $F$ is a \term{Dwyer-Kan equivalence} or
\term{DK-equivalence} when $F$ is a DK-embedding and for every $d\in
\ob\aD$, there exists a $c\in\ob\aC$ such that $\aD(-,d)$ and
$\aD(-,Fc)$ represent naturally isomorphic enriched functors from
$\aD^{\op}$ to the stable category.
\end{defn}

We can rephrase this definition in terms of ``homotopy categories'':
Associated to a spectral category $\aC$, we have the following notion
of homotopy category.

\begin{defn}\label{defhtycat}
For a spectral category $\aC$, the \term{homotopy category}
$\pi_{0}\aC$ is the Ab-category with the same objects, with morphism
abelian groups $\pi_{0}\aC(a,b)$, and with units and composition
induced by the unit and composition maps of $\aC$.  
The \term{graded homotopy category} is the $\mathrm{Ab}_{*}$-category
with objects $\ob\aC$ and morphisms $\pi_{*}\aC(a,b)$.
\end{defn}

We remind the reader that by convention, $\pi_{0}\aC(a,b)$ and
$\pi_{*}\aC(a,b)$ denote the homotopy groups of $\aC(a,b)$ viewed as
an object of the stable category.

Without any further hypotheses on the spectral categories in question,
the following proposition is a straightforward consequence of the
definitions and the Yoneda lemma for enriched
functors \cite[2.4]{kelly}.

\begin{prop}\label{propdkequiv}
A spectral functor $\aC\to \aD$ is a Dwyer-Kan equivalence if and only
if it induces an equivalence of graded homotopy categories
$\pi_{*}\aC\to\pi_{*}\aD$.
\end{prop}

As we will see in Theorems~\ref{thmtriangenv} and~\ref{thmtriang} below,
the homotopy category in practice often has a triangulated structure
compatible with the mapping spectra.  We formalize this in the
following definition.

\begin{defn}\label{deftriang}
A spectral category $\aC$ is \term{pretriangulated} means:
\begin{enumerate}
\item There is an object $0$ in $\aC$ such that the right $\aC$-module
$\aC(-,0)$ is homotopically trivial (weakly 
equivalent to the constant functor with value the one-point symmetric
spectrum $*$).
\item Whenever a right $\aC$-module $\fM$ has the property that $\Sigma\fM$
is weakly equivalent to a representable $\aC$-module $\aC(-,c)$ (for
some object $c$ in $\aC$), then
$\fM$ is weakly equivalent to a representable $\aC$-module
$\aC(-,d)$ for some object $d$ in $\aC$.
\item Whenever the right $\aC$-modules $\fM$ and $\fN$ are weakly equivalent to
representable $\aC$-modules $\aC(-,a)$ and $\aC(-,b)$ respectively, then the
homotopy cofiber of any map of right $\aC$-modules $\fM\to \fN$ is weakly
equivalent to a representable $\aC$-module.
\end{enumerate}
\end{defn}

The first condition guarantees the existence of a zero object in the
homotopy category $\pi_{0}\aC$; the usual argument shows that the left module
$\aC(0,-)$ is also homotopically trivial.  The second condition gives
a desuspension functor on $\pi_{0}\aC$ and the third condition
in particular produces a suspension functor on $\pi_{0}\aC$:
We choose $\Sigma^{-1}a$ and $\Sigma a$ representing
$\Sigma^{-1}\aC(-,a)$ and $\Sigma \aC(-,a)$, respectively, in the
derived category of right $\aC$-modules.  Then $\Sigma^{-1}a$ and
$\Sigma a$ in particular represent the functors $\pi_{1}\aC(-,a)$ and
$\pi_{-1}\aC(-,a)$, respectively, from $\pi_{0}\aC$ to sets, and so
are unique up to unique isomorphism in $\pi_{0}\aC$.

In the third condition, note that maps of right modules from
$\aC(-,a)$ to $\aC(-,b)$ are in one to one correspondence with the
vertices of the zeroth simplicial set of $\aC(a,b)$; using weakly
equivalent $\fM$ and $\fN$, the maps then represent arbitrary elements
of $\pi_{0}\aC(a,b)$.  In the case when all of the mapping spectra of
$\aC$ are fibrant symmetric spectra (e.g., after replacing $\aC$ by
$\aC^{\Omega}$), condition~(iii) can be simplified to considering just
the homotopy cofibers of maps $\aC(-,a)\to\aC(-,b)$.  We explain this
interpretation of condition~(iii) in more detail at the end of the
section in the proof of the following theorem.

\begin{thm}\label{thmtriangenv}
Any small spectral category $\aC$ DK-embeds in a small pretriangulated
spectral category $\tilde \aC$.
\end{thm}

The category $\tilde \aC$ is closely related to the category of
right $\aC$-modules, essentially the closure of (the Yoneda
embedding of) $\aC$ under desuspensions and cofibration
sequences.  The third condition in the definition of pretriangulated
spectral category then indicates the sequences in $\aC$ that are
equivalent to cofibration sequences in $\tilde \aC$.  We can therefore use
the third condition to define the analogue of Puppe (cofibration)
sequences:  We say that a sequence  
\[
a\to b\to c \to \Sigma a
\]
in $\pi_{0}\aC$ is a \term{four term Puppe sequence} if there exists
right $\aC$-modules $\aM$ and $\aN$ and a map of right $\aC$-modules
$f\colon \aM\to \aN$ such that the four term Puppe sequence of $f$
\[
\aM\to \aN\to Cf\to \Sigma \aM
\]
in the category of right $\aC$-modules
is isomorphic in the derived category of right $\aC$-modules to the
sequence  
\[
\aC(-,a)\to \aC(-,b)\to \aC(-,c) \to \aC(-,\Sigma a)
\]
such that the isomorphism $\Sigma \aM\to \aC(-,\Sigma a)\iso \Sigma
\aC(-,a)$ is the suspension of the isomorphism $\aM\to \aC(-,a)$.
We prove the following theorem at the end of the section.

\begin{thm}\label{thmtriang}
If the spectral category $\aC$ is pretriangulated, then its homotopy
category is triangulated with distinguished triangles the four term
Puppe sequences.
A spectral functor between pretriangulated spectral
categories induces a triangulated functor on homotopy categories.
\end{thm}

\begin{cor}\label{cortriang}
A spectral functor $\aC\to \aD$ between pretriangulated spectral
categories is a Dwyer-Kan equivalence if and only
if it induces an equivalence of homotopy categories
$\pi_{0}\aC\to\pi_{0}\aD$. 
\end{cor}

In the context of DG-categories, various analogous conditions have
been given to ensure that the homotopy category of the DG-category is
triangulated \cite{KellerCyclic,BLL,DrinfeldDG}.  Following \cite{BLL}
and \cite{DrinfeldDG}, we refer to such DG-categories as
pretriangulated.  We have the following consistency result.  It
follows from the fact that the various functors relating chain
complexes to $H\bZ$-modules to symmetric spectra all preserve homotopy
cofibration sequences.

\begin{prop}\label{propDGtriangSpec}
If $\aD$ is a pretriangulated DG-category, then its associated
spectral category is pretriangulated.
\end{prop}

As indicated by Proposition~\ref{propdkequiv} and
Corollary~\ref{cortriang}, we take the perspective that the mapping
spectra encode the homotopy theory of the spectral category.  From
this viewpoint, DK-equivalences clearly represent the correct general
notion of weak equivalence of spectral categories.  An alternative
perspective would not require the mapping spectra of a spectral
category $\aC$ to encode all of the homotopy theory, but rather also
include an additional notion of weak equivalence of objects of
$\aC$. For example, this is appropriate in the context of enriched
model categories.  For model categories enriched over symmetric
spectra, the homotopy theory is a localization of the intrinsic
homotopy theory of the associated spectral category.  The full
spectral subcategory of the cofibrant-fibrant objects is the spectral
category whose mapping spectra encode the homotopy theory of the
enriched model category.  This subcategory tends not to be preserved
under most interesting functors.  Under properness hypotheses, a
``cofiber'' version of $THH$ works somewhat better; see
Section~\ref{appcofiber} for more details.

We prove the following invariance theorem for DK-equivalences in the
next section.

\begin{thm}\label{thmdkequiv}
A DK-equivalence $\aC\to \aD$ induces a weak equivalence $THH(\aC)\to
THH(\aD)$. 
\end{thm}

We also prove the following more general theorem for bimodule
coefficients.  In the statement, the $(\aC,\aC)$-bimodule $F^{*}\fN$ is the
bimodule obtained by restriction of scalars; it is the spectral
functor from $\aC^{\op}\sma \aC$ to symmetric spectra defined by first 
applying $F$ to each variable and then applying $\fN$.

\begin{thm}\label{thmbimodequiv}
Let $F\colon \aC\to \aD$ be a DK-equivalence of small spectral categories, $\fM$ a
$(\aC,\aC)$-bimodule and $\fN$ a $(\aD,\aD)$-bimodule.  A weak
equivalence $\fM\to F^{*}\fN$ induces a weak equivalence
$THH(\aC;\fM)\to THH(\aD;\fN)$.
\end{thm}

We now move on from weak equivalences to Morita equivalences.  For
objects $a$ and $c$ of $\aD$, say that $c$ is a \term{homotopy factor}
of $a$ if it is a factor in the graded homotopy category $\pi_{*}\aD$,
i.e., if there exists an object $b$ in $\aD$ and a natural isomorphism
$\pi_{*}\aD(-,c)\iso \pi_{*}\aD(-,a) \times \pi_{*}\aD(-,b)$ of
contravariant functors from $\pi_{*}\aD$ to the category of graded
abelian groups.  We say that a spectral functor $F\colon \aC\to \aD$
is \term{homotopy cofinal} if it induces weak equivalences on mapping
spaces and each object of $\aD$ is a homotopy factor of the image of
some object in $\aC$.  We prove the following theorem in the next
section.

\begin{thm}\label{thmfactor}
A homotopy cofinal spectral functor $\aC\to \aD$ of small spectral
categories induces a weak
equivalence $THH(\aC)\to THH(\aD)$.
\end{thm}

The previous theorem admits a more sophisticated variant.  Given a
set $C$ of objects in a pretriangulated spectral category
$\aD$, the \term{thick closure} of $C$ is the set of objects in
the thick subcategory of $\pi_{0}\aD$ generated by $C$.  In terms of
the spectral category $\aD$, the thick closure of $C$ is the smallest
set $\bar C$ of objects of $\aD$ containing $C$ and satisfying:
\begin{enumerate}
\item If $a$ is a homotopy factor of an object of $\bar C$, then $a$ is
in $\bar C$.
\item If the right $\aD$-module $\Sigma \aD(-,a)$ is weakly equivalent
to $\aD(-,c)$ for some $c$ in $\bar C$, then $a$ is in $\bar C$. 
\item If the right $\aD$-module $\aD(-,a)$ is weakly equivalent to the
cofiber of a map of right $\aD$-modules $\fM\to \fM'$ with $\fM$, $\fM'$
weakly equivalent to $\aD(-,c)$, $\aD(-,c')$ for $c,c'$
in $\bar C$, then $a$ is in $\bar C$.
\end{enumerate}
A set is \term{thick} if it is its own thick closure.  Since
any small spectral category $\aC$ embeds as a full spectral subcategory of a
small pretriangulated spectral category $\aD$, the following theorem, proved
in the next section, in particular allows us to always reduce
questions in $THH$ to the case of small pretriangulated spectral
categories.

\begin{thm}\label{thmthick}
Let $\aD$ be a pretriangulated spectral category.  Let $C$ be a set of
objects of $\aD$, $\bar C$ its thick closure, and $C'$ a set
containing $C$ and contained in $\bar C$. Let $\aC$ and $\aC'$ be the
full spectral subcategories of $\aD$ on the objects in $C$ and $C'$
respectively.  Then the inclusion $\aC\to \aC'$ induces a weak 
equivalence $THH(\aC) \to THH(\aC')$.
\end{thm}

We close the section with the proofs of Theorems~\ref{thmtriangenv}
and~\ref{thmtriang}.  The argument involves the well-known properties
of categories of enriched functors into a Quillen closed model
category.  For any small spectral category $\aC$, the category
$\RMod{\aC}$ of right $\aC$-modules has a \term{standard model
structure} (or projective model structure) that is proper and
compactly generated, where the generating cofibrations and generating
acyclic cofibrations are the maps $\aC(-,c)\sma f$ for $c$ in $\aC$
and $f$ varying through the generating cofibrations and generating
acyclic cofibrations (respectively) of the stable model structure on
symmetric spectra described in \cite[3.3.2,3.4.9]{HSS} (see also
\cite[3.4.2.1,3.4.16]{HSS}).  Consequently, the weak equivalences and
fibrations are the maps that are objectwise weak equivalences and
objectwise fibrations (respectively) in the stable model structure on
symmetric spectra.  The representable right $\aC$-modules $\aC(-,c)$
are cofibrant and \term{compact}, meaning that maps out of $\aC(-,c)$ preserve
sequential colimits.  In fact, the set of maps, simplicial set of
maps, and symmetric spectrum of maps out of $\aC(-,c)$ preserve
arbitrary colimits, by the enriched Yoneda lemma \cite[2.4]{kelly}.

\begin{proof}[Proof of Theorem~\ref{thmtriang}]
In the case when $\aC$ is small, we use the model theory above as
follows.  Using Quillen's theory of cofibration sequences, we obtain a
triangulated structure on the Quillen homotopy category
$\Ho\RMod{\aC}$ of $\RMod{\aC}$.  The homotopy category $\pi_{0}\aC$
embeds as a full subcategory of $\Ho\RMod{\aC}$, and the conditions in
the definition of pretriangulated spectral category imply that
$\pi_{0}\aC$ is closed under desuspensions, suspensions, and triangles
in $\Ho\RMod{\aC}$.  In the case when $\aC$ is not small, $\RMod{\aC}$
does not typically have small Hom-sets.  To handle this, one restricts
to the cell modules and uses small object arguments to construct
replacements whose values on a given small set of objects are fibrant
symmetric spectra.

Now given a spectral functor $F\colon \aC\to \aD$ between small
pretriangulated spectral categories, left Kan extension produces a
functor $\Lan_{F}\colon \RMod{\aC}\to\RMod{\aD}$ left adjoint to the
functor $F^{*}\colon \RMod{\aD}\to \RMod{\aC}$.  Since $F^{*}$
preserves fibrations and weak equivalences in the model structure
above, $\Lan_{F}$ and $F^{*}$ form a Quillen adjoint pair.  In
particular $\Lan_{F}$ preserves Quillen cofibration sequences and Quillen
suspensions.  It follows that the left derived functor of $\Lan_{F}$
on Quillen homotopy categories is triangulated; on the representable
functors, the left derived functor $\Lan_{F}$ is just $\pi_{0}F\colon
\pi_{0}\aC\to\pi_{0}\aD$.  Again, the case when $\aC$ or $\aD$ is not
small may be handled by a straightforward direct argument in terms of
the cell modules.
\end{proof}

\begin{proof}[Proof of Theorem~\ref{thmtriangenv}]
We again use the model theory on categories of modules.  By
Proposition~\ref{propfibrep}, we can assume without loss of generality
that all of the mapping spectra $\aC(x,y)$ are fibrant in the stable
model structure on symmetric spectra, and so the representable right
$\aC$-modules $\aC(-,c)$ are both cofibrant and fibrant in the model
structure on $\RMod{\aC}$.  In order to remain in the setting of small
categories, we restrict to a small full subcategory of $\RMod{\aC}$ as
follows.  For any set $U$, write $U\RMod{\aC}$ for the full
subcategory of $\RMod{\aC}$ consisting of the functors that take
values in symmetric spectra whose underlying sets are in $U$.  Then
$U\RMod{\aC}$ is a small spectral category, and if we choose $U$ to be
the power set of a sufficiently large cardinal, then $U\RMod{\aC}$
will be closed (up to universal isomorphisms) under the usual
(cardinal limited) constructions of homotopy theory in $\RMod{\aC}$,
including the small objects argument constructing factorizations.  In
particular, $U\RMod{\aC}$ is a Quillen model category (assuming just
finite limits and colimits) with cofibrations, fibrations, and weak
equivalences the maps that are such in $\RMod{\aC}$.  We then get a
(closed model) category $\RMod{U\RMod{\aC}}$.

Let $\tilde\aC$ be the full spectral subcategory of $U\RMod{\aC}$
consisting of the cofibrant-fibrant objects.  The enriched Yoneda
lemma embeds $\aC$ as a full spectral subcategory of $\tilde\aC$.
Properties~(i) and (ii) for $\tilde\aC$ in the definition of
pretriangulated spectral category are clear.  For property~(iii),
consider a map of right $\tilde\aC$-modules $\fM\to \fN$.  Since the
model structure on $\RMod{\tilde\aC}$ is left proper, after replacing
$\fM$ and $\fN$ with fibrant replacements, we obtain an equivalent
homotopy cofiber, and so we can assume without loss of generality that
$\fM$ and $\fN$ are fibrant.  We assume that $\fM$ is weakly
equivalent to $\tilde \aC(-,a)$ and $\fN$ is weakly equivalent to
$\tilde \aC(-,b)$ for objects $a,b$ in $\tilde \aC$; since $\tilde
\aC(-,a)$ and $\tilde \aC(-,b)$ are cofibrant and $\fM$ and $\fN$ are
fibrant, we can choose weak equivalences $\tilde \aC(-,a)\to \fM$ and
$\tilde \aC(-,b)\to \fN$.  Furthermore, as $\tilde\aC(a,b)$ and
$\fN(a)$ are both fibrant, we can lift the composite map
$\tilde\aC(-,a)\to \fN$ to a homotopic map $\tilde\aC(-,a)\to
\tilde\aC(-,b)$.  We get a weak equivalence on the homotopy cofibers.
The map $\tilde\aC(-,a)\to \tilde\aC(-,b)$ comes from a map $a\to b$
by the Yoneda lemma.  A fibrant replacement of the homotopy cofiber
in $U\RMod{\aC}$ is in $\tilde\aC$ and represents the homotopy cofiber
of $\fM\to \fN$ in $\RMod{\tilde\aC}$.  This completes the proof of
Theorem~\ref{thmtriangenv}.
\end{proof}

\section{The Dennis-Waldhausen Morita Argument}\label{secdwm}

In this section, we consider the invariance properties of $THH$ from
the perspective of generalized Morita theory.  Dennis and Waldhausen
gave a very concrete argument for the Morita invariance of the
Hochschild homology of rings using an explicit bisimplicial
construction \cite[p.~391]{WaldhausenA2} .  We give a broad
generalization of this argument to the setting of spectral categories
that provides the technical foundations for the proofs of the theorems
of the previous section as well as the arguments in the remainder of
the paper.

The argument uses the two-sided bar construction.

\begin{defn}
Let $\aC$ be a small spectral category, $\fM$ a right $\aC$-module, and $\fN$
a left $\aC$-module.  The \term{two-sided bar construction}
$\TB(\fM;\aC;\fN)$ is the diagonal (geometric realization) of the simplicial
symmetric spectrum $\TB\subdot(\fM;\aC;\fN)$, where
\[
\TB_{q}(\fM;\aC;\fN)
=\bigvee \fM(c_{q})\sma \aC(c_{q-1},c_{q}) \sma \dotsb \sma
\aC(c_{0},c_{1}) \sma \fN(c_{0}),
\]
where the sum is over the $(q+1)$-tuples $(c_{0},\dotsc,c_{q})$ of
objects of $\aC$.  We make this a simplicial object with the usual
two-sided bar construction face and degeneracy maps: the zeroth
face map is induced by the action of $\aC$ on $\aN$, the last face map
is induced by the action of $\aC$ on $\aM$, and the remaining face maps
are induced by the composition in $\aC$. The degeneracy maps are induced
by the unit maps $S\to \aC(c_{i},c_{i})$.
\end{defn}

The following is the main technical proposition of this section.
In it and elsewhere when necessary for clarity, we write 
\[
\TB(\fM(x); x,y \in \aC; \fN(y))
\qquad \text{and}\qquad 
\THM(x,y\in \aC;\fP(x,y))
\]
for $\TB(\fM;\aC;\fN)$ and $\THM(\aC;\fP)$, especially when $\fM$,
$\fN$, and/or $\fP$ depend on other variables.   

\begin{prop}[Dennis-Waldhausen Morita Argument]\label{propcoremorita}
Let $\aC$ and $\aD$ be small spectral categories. Let $\fP$ be
a $(\aD,\aC)$-bimodule and $\fQ$ a $(\aC,\aD)$-bimodule.  Then there
is a natural isomorphism of symmetric spectra 
\[
\THM(\aC,\TB(\fP,\aD,\fQ))\iso \THM(\aD;\TB(\fQ,\aC,\fP)),
\]
that is, 
\begin{multline*}
\THM(x,y\in \aC;\TB(\fP(w,y);w,z\in \aD;\fQ(x,z)))\\
\iso
\THM(w,z\in \aD;\TB(\fQ(x,z);x,y\in \aC;\fP(w,y))).
\end{multline*}
\end{prop}

\begin{proof}
We can identify both symmetric spectra
\[
\THM(\aC;B(\fP;\aD;\fQ))\qquad \text{and}\qquad
\THM(\aD;B(\fQ;\aC;\fP))
\]
as the diagonal of the bisimplicial spectrum
with $(q,r)$-simplices as pictured.
\[
\xymatrix@=1pc{%
&\aC(c_{q-1},x)\sma \dotsb \sma \aC(y,c_{1})\\
\fQ(x,z)
\ar@{{}{}{}}[dr]_(.3){\relax\textstyle\sma}
\ar@{{}{}{}}[ur]^(.3){\relax\textstyle\sma}
&&\fP(w,y)
\ar@{{}{}{}}[dl]^(.3){\relax\textstyle\sma}
\ar@{{}{}{}}[ul]_(.3){\relax\textstyle\sma}
\\
&\aD(z,d_{1})\sma \dotsb \sma \aD(d_{r-1},w)
}
\]
These two constructions are therefore canonically isomorphic in the
point-set category of symmetric spectra.
\end{proof}

The following lemma complements Proposition~\ref{propcoremorita}.  Its
proof is the usual simplicial contraction (see for example
\cite[9.8]{GILS}).

\begin{lem}[Two-Sided Bar Lemma]\label{lemtwobar}
Let $\aC$ be a small spectral category, let $\aM$ be a right
$\aC$-module, and let $\aN$ be a left $\aC$-module. For any object $c$
in $\aC$, the composition maps
\[
\TB\subdot(\fM;\aC;\aC(c,-))\to \fM(c)
\qquad\text{and}\qquad 
\TB\subdot(\aC(-,c);\aC;\fN)\to \fN(c)
\]
are simplicial homotopy equivalences.
\end{lem}

We use the previous proposition and lemma to prove the following
theorem, which gives a criterion for converting objectwise equivalence
conditions into equivalences on $THH$.

\begin{thm}\label{thmcompcrit}
Let $\aC$ and $\aD$ be small spectral categories and let
$F\colon \aC \to \aD$ be a spectral functor. Let $\fM$ be a
$(\aC,\aC)$-bimodule, $\fN$ a $(\aD,\aD)$-bimodule and $\fM\to
F^{*}\fN$ a weak equivalence.  Assume that $\aC$ and $\aD$ are
pointwise cofibrant.  Suppose that the map of symmetric spectra
\begin{equation}
\tag{*}
\TB(\aD(F-,z);\aC;\fN(w,F-))\to 
\TB(\aD(-,z);\aD;\fN(w,-))
\end{equation}
is a weak equivalence for each fixed $w$,$z$ in $\aD$.  
Then the map
\[
THH(\aC;\fM)\to THH(\aD;\fN)
\]
is a weak equivalence.
\end{thm}

\begin{proof}
Consider the commutative diagram
\begin{equation}\label{eqpfdkequiv}
\begin{gathered}
\def\shift{.65}
\xymatrix{%
THH(\aC;\TB(\fM(-,-);\aC;\aC(-,-)))\ar[r]^(\shift){\simeq}\ar[d]&THH(\aC;\fM)\ar[d]^{\simeq}\\
THH(\aC;\TB(\fN(-,F-);\aD;\aD(F-,-)))\ar[r]^(\shift){\simeq}\ar[d]&THH(\aC;F^{*}\fN)\ar[d]\\
THH(\aD;\TB(\fN(-,-);\aD;\aD(-,-)))\ar[r]_(\shift){\simeq}&THH(\aD;\fN).
}
\end{gathered}
\end{equation}
The arrows marked ``$\simeq$'' are weak equivalences by the Two-Sided
Bar Lemma and Proposition~\ref{propnibus} above.  The composite map on the
right is the map we are interested in, and so our goal is to prove
that the two maps on the left are weak 
equivalences.  The first map is induced by a map of $\aC$-bimodules 
\[
\TB(\fM(-,-);\aC;\aC(-,-))\to \TB(\fN(-,F-);\aD;\aD(F-,-)),
\]
which is easily seen to be a weak equivalence by the Two-Sided Bar
Lemma and the hypothesis that $\fM\to F^{*}\fN$ is a weak equivalence.
Since we have assumed that $\aC$ and $\aD$ are pointwise cofibrant, by
Proposition~\ref{propTHHvsTHM}, to show that the second map is a weak
equivalence, it suffices to show that the map
\[
\THM(\aC;\TB(\fN(-,F-);\aD;\aD(F-,-)))\to
\THM(\aD;\TB(\fN(-,-);\aD;\aD(-,-)))
\]
is a weak equivalence and here we apply  
Proposition~\ref{propcoremorita}:  We obtain a commutative diagram with
the horizontal maps isomorphisms
\[
\xymatrix@C=1pc{%
\THM(\aC;\TB(\fN(-,F-);\aD;\aD(F-,-)))\ar[d]\ar[r]^{\iso}
&\THM(\aD;\TB(\aD(F-,-);\aC;\fN(-,F-)))\ar[d]\\
\THM(\aD;\TB(\fN(-,-);\aD;\aD(-,-)))\ar[r]_{\iso}
&\THM(\aD;\TB(\aD(-,-);\aD;\fN(-,-)))
}
\]
by applying Proposition~\ref{propcoremorita} with
$\fP=\fN(-,F-)$ and $\fQ=\aD(F-,-)$ on the top
and $\fP=\fN(-,-)$ and $\fQ=\aD(-,-)$ (for $\aC=\aD$) on the bottom. 
By the hypothesis on the map~(*), the map of $(\aD,\aD)$-bimodules
\[
\TB(\aD(F-,-);\aC;\fN(-,F-))\to 
\TB(\aD(-,-);\aD;\fN(-,-))
\]
is a weak equivalence, and it follows that
the vertical maps above are weak equivalences.
\end{proof}

We now apply this criterion in the proof of Theorem~\ref{thmdkequiv},
Theorem~\ref{thmbimodequiv}, Theorem~\ref{thmfactor}, and
Theorem~\ref{thmthick}.

\begin{proof}[Proof of Theorem~\ref{thmdkequiv}]
Using Proposition~\ref{propcofrep}, we can assume without loss of
generality that $\aC$ and $\aD$ are pointwise cofibrant, and then 
apply Theorem~\ref{thmcompcrit} with $\fM=\aC$ and $\fN=\aD$.
We need to show that the map
\[
\TB(\aD(F-,z);\aC;\aD(w,F-))\to 
\TB(\aD(-,z);\aD;\aD(w,-))
\]
is a weak equivalence for each $w,z$ in $\aD$, and using the Two-Sided
Bar Lemma, it suffices to show that the composite map
\begin{equation}\label{eqtempdkequiv}
\TB(\aD(F-,z);\aC;\aD(w,F-))\to 
\aD(w,z)
\end{equation}
is a weak equivalence.  Viewing
$\TB\subdot(\aD(F-,z);\aC;\aD(w,F-))$ as a simplicial object in the
stable category, up to simplicial isomorphism, it only depends on
$\aD(F-,z)$ as a functor from $\aC^{\op}$ to the stable category.  By
hypothesis, there exists an object $c$ in $\aC$ such that $\aD(-,z)$
and $\aD(-,Fc)$ are isomorphic as functors from $\aD^{\op}$ to the stable
category.  It follows that $\aD(F-,z)$ and $\aD(F-,Fc)\simeq \aC(-,c)$
are isomorphic as functors from $\aC^{\op}$ to the stable category.
Since this is just a comparison of simplicial objects in
the stable category, we do not get a direct comparison on geometric
realizations (but see also the proof of Theorem~\ref{thmbimodequiv}
below).  Nonetheless, the 
homotopy groups of $\TB\subdot$ are the $E_{1}$-term of a spectral
sequence that computes the homotopy groups of $\TB$. 
The $E_{1}$ differential comes from the simplicial face maps, and
applying the Two-Sided Bar Lemma, we see that this spectral sequence 
degenerates at $E_{2}$ and that \eqref{eqtempdkequiv} is a weak equivalence.
\end{proof}

For the proof of Theorem~\ref{thmbimodequiv}, we note that the
map~(*) in the statement of Theorem~\ref{thmcompcrit} still makes
sense when we replace $\aD(-,z)$ with an arbitrary right $\aD$-module $\phi$:
\begin{equation}\label{eqextmod}
\TB(F^{*}\phi;\aC;\fN(w,F-))\to 
\TB(\phi;\aD;\fN(w,-)).
\end{equation}
We have used $\phi$ to denote the right $\aD$-module to avoid possible
confusion between the different roles played by the right module $\phi$ and
the bimodule $\fN$.  If $\fN(w,-)$ is objectwise cofibrant, then
$\TB(F^{*}\phi;\aC;\fN(w,F-))$ and $\TB(\phi;\aD;\fN(w,-))$ preserve
weak equivalences in $\phi$.

\begin{proof}[Proof of Theorem~\ref{thmbimodequiv}]
We can assume without loss of generality that $\aC$ and $\aD$ are
pointwise cofibrant and that $\fM$ and $\fN$ are objectwise
cofibrant.  Applying Theorem~\ref{thmcompcrit}, it suffices to show
that the map~(*)
\[
\TB(\aD(F-,z);\aC;\fN(w,F-))\to 
\TB(\aD(-,z);\aD;\fN(w,-))
\]
is a weak equivalence for each $w,z$ in $\aD$.  Fixing $w,z$, it
suffices to show that the
map~\eqref{eqextmod} is a weak equivalence for $\phi =\aD(-,z)$.
Let $\psi$ be a right $\aD$-module fibrant replacement of
$\aD(-,z)$.  By hypothesis, viewing $\psi$ as an enriched functor from
$\aD$ to the stable category, we have a natural isomorphism $\tilde
f\colon \aD(-,Fz')\to \psi$ for some $z'$ in $\aC$; by the Yoneda
lemma for enriched functors, this corresponds to an element $\tilde
f\in \pi_{0}(\psi(Fz'))$.  Since $\psi(Fz')$ is fibrant, we can choose a
vertex $f$ in $\psi(Fz')_{0}$ representing $\tilde f$.  Again by the
Yoneda lemma, $f$ represents a map of right $\aD$-modules
$\aD(-,Fz')\to \psi$ that induces the natural isomorphism $\tilde f$ of
enriched functors to the stable category.  In particular, $f$ is a
weak equivalence.  The map~\eqref{eqextmod} is a weak equivalence for
$\phi =\aD(-,Fz')$, and so is a weak equivalence for $\phi =\psi$ and for
$\phi =\aD(-,z)$.
\end{proof}

Theorem~\ref{thmfactor} can be proved using essentially the same
argument as the proof of Theorem~\ref{thmdkequiv} above, using the
fact that a direct sum of maps in the stable category is an
isomorphism if and only if it is an isomorphism on each factor.  On
the other hand, given Theorem~\ref{thmtriangenv},
Theorem~\ref{thmfactor} follows from Theorem~\ref{thmdkequiv} and
Theorem~\ref{thmthick}, which we now prove.

\begin{proof}[Proof of \ref{thmthick}]
Without loss of generality, we can assume that $\aD$ is pointwise
cofibrant, and then $\aC$ and $\aC'$ are also pointwise cofibrant.
Applying Theorem~\ref{thmcompcrit}, it suffices to show that the map~(*)
\[
\TB(\aC'(-,z);\aC;\aC'(w,-))\to 
\TB(\aC'(-,z);\aC';\aC'(w,-))
\]
is a weak equivalence for every $w,z$ in $\aC'$, or more generally,
that the map~\eqref{eqextmod} 
\begin{equation}\label{eqtmpthmthick}
\TB(\phi;\aC;\aC'(w,-))\to 
\TB(\phi;\aC';\aC'(w,-))
\end{equation}
is a weak equivalence for any right $\aC'$-module $\phi$ that is weakly
equivalent to $\aC'(-,z)$ for $z$ in $\aC'$. 
By the Two-Sided Bar Lemma, we know that~\eqref{eqtmpthmthick} is a
weak equivalence when $\phi$ is $\aC(-,x)$ for $x$ in $\aC$.
Using the fact that both
sides preserve homotopy cofibration sequences in the $\phi$ variable, it
follows that \eqref{eqtmpthmthick} is a weak equivalence for
$\aD(-,z)$ for any $z$ in the thick subcategory of
$\pi_{0}\aD$ generated by $\pi_{0}\aC$.  This completes the
proof of Theorem~\ref{thmthick}.
\end{proof}

\section{The general localization theorem}\label{secgenloc}

In this section, we discuss a general theorem that produces
localization cofibration sequences in $THH$.  The basic strategy takes
advantage of the fact that $THH$ preserves (co)fibration sequences in the
bimodule variable: We apply the Dennis-Waldhausen Morita argument to
identify $THH$ of a small spectral category with $THH$ of another
small spectral
category with coefficients in a bimodule; see Lemma~\ref{lemldoesa}
below.  Roughly, we then obtain our 
localization cofibration sequences by reinterpreting sequences of $THH$ of
small spectral categories 
\[
THH(\aA)\to THH(\aB)\to THH(\aC)
\]
up to weak equivalence
as the $THH$ of a single small spectral category with
coefficients in a sequence of judiciously chosen bimodules
\[
THH(\aB;\fLL) \to THH(\aB; \aB) \to THH(\aB; \fQQ)
\]
where $\fQQ$ is the cofiber of a map of $(\aB,\aB)$-bimodules $\fLL \to
\aB$.  Although we can make more general
statements, the situation we 
are most interested in is when the sequence of small spectral categories
models a triangulated quotient (in the sense of Verdier).  We prove
the following theorem.

\begin{thm}\label{thmgenone}
Let $F\colon \aB\to \aC$ be a spectral functor between small
pretriangulated spectral categories, and let $\aA$ be the full spectral
subcategory of $\aB$ consisting of the objects $a$ such that $F(a)$ is
isomorphic to zero in the homotopy category $\pi_{0}\aC$. If the
induced map from the triangulated quotient $\pi_{0}\aB/\pi_{0}\aA$ to
$\pi_{0}\aC$ is cofinal, then
$THH(\aC)$ is weakly equivalent the homotopy cofiber of $THH(\aA)\to
THH(\aB)$.  The analogous relationship holds also for $TR(\aC)$ and
$TC(\aC)$. 
\end{thm}

In general, we call $(\aB,\aA)$ a \term{localization pair} when $\aB$
is a pretriangulated spectral category and $\aA$ is a full spectral
subcategory such that $\pi_{0}\aA$ is thick in $\pi_{0}\aB$; we say
that the localization pair is small when the spectral category $\aB$
is small.
This definition of localization pair differs slightly from that of
Keller \cite[2.4]{KellerCyclic} in that we do not require a well-behaved
ambient category (our additional requirement that $\aA$ be thick is
for convenience rather than necessity by Theorem~\ref{thmthick}).

In Theorem~\ref{thmgenone}, letting $\aZ$ be the full subcategory of
objects of $\aC$ in the thick 
closure of the image of $\aA$, then $(\aC,\aZ)$ is a localization pair
and $(\aB,\aA)\to (\aC,\aZ)$ is a map of localization pairs: It is a
spectral functor $\aB\to\aC$ that takes $\aA$ into $\aZ$.  Note that
for any objects $x,y$ in $\aZ$, the spectrum $\aZ(x,y)$ is
homotopically trivial, and so $THH(\aZ)$ is homotopically trivial.
The inclusion of $THH(\aC)$ in the homotopy cofiber of
$THH(\aZ)\to THH(\aC)$ is therefore a weak equivalence.
We have the map of homotopy cofibers 
\[
C(THH(\aA)\rightarrow THH(\aB)) \quad \to \quad
C(THH(\aZ)\rightarrow THH(\aC)),
\]
and we prove Theorem~\ref{thmgenone} by showing that this map is a
weak equivalence: We get the statement for $TR(\aC)$ and $TC(\aC)$
because the weak equivalence on $THH$ above implies that the map 
of cyclotomic spectra
\[
C(THH(\aA)\rightarrow THH(\aB)) \quad \to \quad
C(THH(\aZ)\rightarrow THH(\aC))
\]
is a weak equivalence.

Theorem~\ref{thmgenone} then naturally appears as a special case
of the following theorem, which essentially says that the cofiber of
$THH$ is an invariant of the localization pair.

\begin{thm}\label{thmgentwo}
Let $F\colon (\aB_{1},\aA_{1})\to (\aB_{2},\aA_{2})$ be a map of
small localization pairs.  If the induced
map of triangulated quotients 
\[
\pi_{0}\aB_{1}/\pi_{0}\aA_{1} \to \pi_{0}\aB_{2}/\pi_{0}\aA_{2}
\]
is cofinal, then the induced map of cyclotomic spectra
\[
C(THH(\aA_{1})\rightarrow THH(\aB_{1})) \quad \to \quad
C(THH(\aA_{2})\rightarrow THH(\aB_{2}))
\]
is a weak equivalence.
\end{thm}

For ease of reference, we also state the following immediate
corollary.

\begin{cor}\label{corgentwo}
Under the hypotheses and notation of the previous theorem, the maps 
\begin{gather*}
C(TR(\aA_{1})\rightarrow TR(\aB_{1})) \quad \to \quad
C(TR(\aA_{2})\rightarrow TR(\aB_{2}))\\
C(TC(\aA_{1})\rightarrow TC(\aB_{1})) \quad \to \quad
C(TC(\aA_{2})\rightarrow TC(\aB_{2}))
\end{gather*}
are weak equivalences.
\end{cor}

Following Keller \cite{KellerCyclic}, we define $THH$ of a
localization pair as the cofiber on $THH$.  The previous theorem
provides a perspective and justification for the following definition.

\begin{defn}
Let $(\aB,\aA)$ be a small localization pair.  We write $\CTHH(\aB /
\aA)$ for the cyclotomic spectrum obtained as the cofiber of the map 
$THH(\aA)\to THH(\aB)$.  
\end{defn}

We now move on to the proof of
Theorem~\ref{thmgentwo}.  Applying Proposition~\ref{propcofrep},
we assume without loss of generality that $\aB$ and $\aC$ are
pointwise cofibrant and it follows that $\aA$ is pointwise
cofibrant.  We write $\CTHM(\aB/\aA)$ for the cofiber of the map
$\THM(\aA)\to \THM(\aB)$.  Proposition~\ref{propTHHvsTHM} now reduces
the proof of Theorem~\ref{thmgentwo} to the following lemma, whose
proof occupies the remainder of the section.

\begin{lem}\label{lemlocmain}
Let $F\colon (\aB_{1},\aA_{1})\to (\aB_{2},\aA_{2})$ be a map of
small localization pairs with $\aB_{1}$ and $\aB_{2}$ pointwise
cofibrant.  If the induced map of triangulated quotients 
\[
\pi_{0}\aB_{1}/\pi_{0}\aA_{1} \to \pi_{0}\aB_{2}/\pi_{0}\aA_{2}
\]
is cofinal, then the induced map
$\CTHM(\aB_{1}/\aA_{1})\to \CTHM(\aB_{2}/\aA_{2})$
is a weak equivalence.
\end{lem}

As a first step, we use the Dennis-Waldhausen
Morita argument, Proposition~\ref{propcoremorita}, to rewrite
$\THM(\aA)$ as $\THM(\aB;\fLL)$ for an appropriate $(\aB,\aB)$-bimodule $\fLL$.

\begin{lem}\label{lemldoesa}
Let $(\aB,\aA)$ be a small localization pair with $\aB$ pointwise
cofibrant, and let 
$\fLL$ be the $(\aB,\aB)$-bimodule defined by
\[
\fLL(x,y)=\TB(\aB(-,y);\aA;\aB(x,-)).
\]
Then $\THM(\aA)$ is naturally (in maps of pointwise cofibrant
localization pairs) weakly equivalent to $\THM(\aB;\fLL)$.
\end{lem}

\begin{proof}
We apply Proposition~\ref{propcoremorita} with $\aC=\aA$,
$\aD=\aB$, $\aP=\aB$, and $\aQ=\aB$ to obtain a natural isomorphism
\[
\THM(\aB;\fLL)=
\THM(\aB;\TB(\aB;\aA;\aB))
\iso 
\THM(\aA;\TB(\aB;\aB;\aB)).
\]
The natural map
\[
THH(\aA;\TB(\aB;\aB;\aB))\to THH(\aA;\aB)= THH(\aA;\aA)=THH(\aA).
\]
is a weak equivalence since the composition map of
$(\aA,\aA)$-bimodules $\TB(\aB;\aB;\aB)\to \aB$ is a weak equivalence
by the Two-Sided Bar Lemma~\ref{lemtwobar}.
\end{proof}

For a small localization pair $(\aB,\aA)$, write $\fQQ$ for the
$(\aB,\aB)$-bimodule obtained as the cofiber of the composition map
$\fLL\to\aB$.  Then by the previous lemma, we have a natural weak equivalence 
\[
\THM(\aB;\fQQ)\simeq \CTHH(\aB / \aA).
\]
Naturality here refers to the fact that a map of small localization
pairs $F$ induces a map of $(\aB_{1},\aB_{1})$-bimodules $\fQQ[1]\to
F^{*}\fQQ[2]$ and therefore a map 
\begin{equation}\label{eqgentwo}
\THM(\aB_{1};\fQQ[1])\to \THM(\aB_{2};\fQQ[2]).
\end{equation}
Looking at the proof of Lemma~\ref{lemldoesa}, we see that this map is
compatible under the weak equivalences above with the map on homotopy
cofibers in the statement of Lemma~\ref{lemlocmain}.  Thus, to prove
Lemma~\ref{lemlocmain}, we just need to 
show that the map~\eqref{eqgentwo} is a weak equivalence.

For a small localization pair $(\aB,\aA)$ and fixed object $b$ in
$\aB$, the right $\aB$-module $\fLL(-,b)$ is the enriched homotopy
left Kan extension along $\aA\to \aB$ of the enriched functor
$\aB(-,b)$ from $\aA$ to symmetric spectra.  (By which we mean the
derived functor of the enriched left Kan extension \cite[\S 4]{kelly}.)
Philosophically, the 
homotopy cofiber of the map $\fLL(-,b)\to \aB(-,b)$ should then
represent the right $\aC$-module of maps into the image of $b$ in any
spectral category $\aC$ representing the triangulated quotient,
cf.~\cite[(1.3)]{DrinfeldDG}.  From this perspective and viewed
through the principles of the Dennis-Waldhausen Morita argument,
$THH(\aB;\fQQ)$ should be equivalent to $THH(\aC)$.  This is the
idea behind the following lemma proved at the end of the section.

\begin{lem}\label{lemNoisNt}
Under the hypotheses of Lemma~\ref{lemlocmain}, the map of
$(\aB_{1},\aB_{1})$-bimodules $\fQQ[1]\to F^{*}\fQQ[2]$ is 
a weak equivalence.
\end{lem}

A fundamental property of $\fQQ$ is that $\fQQ(a,-)$ and
$\fQQ(-,a)$ are homotopically trivial for any object $a$ in $\aA$: The
Two-Sided Bar Lemma~\ref{lemtwobar} implies 
that the composition maps $\fLL(a,-)\to \aB(a,-)$ and
$\fLL(-,a)\to \aB(-,a)$ are weak equivalences.  This leads to
the following technical observation needed below to analyze the
map~\eqref{eqgentwo}.

\begin{lem}\label{lemNNisBN}
For a small localization pair $(\aB,\aA)$, the maps of bimodules 
\[
\TB(\aB;\aB;\fQQ)\to \TB(\fQQ;\aB;\fQQ)
\quad \text{and}\quad 		
\TB(\fQQ;\aB;\aB)\to \TB(\fQQ;\aB;\fQQ)
\]
induced by $\aB \to \fQQ$ are weak equivalences.
\end{lem}

\begin{proof}
We prove the first equivalence; the argument for the second is similar.
Expanding $\fQQ$ in terms of its
definition, we see that 
$\TB(\fQQ;\aB;\fQQ)$ is the cofiber of the bimodule map
$\TB(\fLL;\aB;\fQQ)\to \TB(\aB;\aB;\fQQ)$, and so by the Two-Sided Bar
Lemma~\ref{lemtwobar}, it suffices to see that $\TB(\fLL;\aB;\fQQ)$ is
homotopically trivial.  Since $\fLL=\TB(\aB;\aA;\aB)$, expanding out 
$\TB(\fLL;\aB;\fQQ)$ as a bisimplicial object, we get the isomorphism  
\[
\TB(\fLL;\aB;\fQQ)=
\TB(\TB(\aB,\aA,\aB);\aB;\fQQ)\iso
\TB(\aB;\aA;\TB(\aB;\aB;\fQQ)).
\]
Applying the Two-Sided Bar
Lemma~\ref{lemtwobar} again, we see that
$\TB(\aB;\aA;\TB(\aB;\aB;\fQQ))$ is weakly equivalent to 
$\TB(\aB;\aA;\fQQ)$.  This is homotopically trivial because the restriction of
$\fQQ(x,-)$ to $\aA$ is homotopically trivial for any $x$. 
\end{proof}

We can extend $\fQQ$ to be a $(\aB,\RMod{\aB})$-bimodule,
where $\RMod{\aB}$ denotes the category of right
$\aB$-modules.  For $x$ an object of $\aB$ and $\phi$ a right
$\aB$-module, let $\fQQ(x,\phi)$ be the cofiber of the
composition map
\[
\TB(\phi(-);\aA;\aB(x,-))\to \phi(x).
\]
Clearly, $\fQQ(x,\phi)$ is isomorphic to $\fQQ(x,y)$ when $\phi
=\aB(-,y)$, and $\fQQ(x,-)$ sends cofibration sequences of right
$\aB$-modules to cofibration sequences of symmetric spectra and sends weak
equivalences of right $\aB$-modules to weak equivalences of symmetric
spectra.  The usual category of 
fractions description of the triangulated quotient
$\pi_{0}\aB/\pi_{0}\aA$ and the fact that $\fQQ$ is homotopically trivial when
either variable is in $\aA$ then implies that $\fQQ(-,y)$ and
$\fQQ(-,y')$ are weakly equivalent right $\aB$-modules when $y$
and $y'$ are isomorphic in $\pi_{0}\aB/\pi_{0}\aA$.
Moreover, when $z$ is isomorphic to $w\vee y$ in
$\pi_{0}\aB/\pi_{0}\aA$, $\fQQ(-,z)$ is weakly equivalent
as a right $\aB$-module to $\fQQ(-,w)\vee \fQQ(-,y)$.
Using these observations and the lemmas above, we
can now prove Lemma~\ref{lemlocmain}.

\begin{proof}[Proof of Lemma~\ref{lemlocmain}]
We need to show that the map~\eqref{eqgentwo} is a weak equivalence.
Consider the following commutative diagram
\[ 
\xymatrix@C-1pc{%
\THM(\aB_{1};F^{*}\TB(\fQQ[2];\aB_{2};\fQQ[2]))\ar[d]
&\THM(\aB_{1};F^{*}\TB(\aB_{2};\aB_{2};\fQQ[2]))
 \ar[d]\ar[l]_-{\simeq}\ar[r]^-{\simeq}
&\THM(\aB_{1};F^{*}\fQQ[2])\ar[d]\\
\THM(\aB_{2};\TB(\fQQ[2];\aB_{2};\fQQ[2]))
&\THM(\aB_{2};\TB(\aB_{2};\aB_{2};\fQQ[2]))\ar[r]_-{\simeq}\ar[l]^-{\simeq}
&\THM(\aB_{2};\fQQ[2]).
}
\]
The lefthand horizontal maps are weak equivalences by
Lemma~\ref{lemNNisBN} and the righthand horizontal maps are weak
equivalences by the Two-Sided Bar
Lemma~\ref{lemtwobar}.  The map~\eqref{eqgentwo} is the composite of
the righthand vertical map and the induced map on $\THM$ of the map of
bimodules 
$\fQQ[1]\to F^{*}\fQQ[2]$, which is a weak equivalence by
Lemma~\ref{lemNoisNt}.
Thus, to see that \eqref{eqgentwo} is a weak equivalence, it suffices
to show that one of  the
vertical maps is a weak equivalence.

Focusing on the lefthand vertical map and
applying the Dennis-Waldhausen
Morita argument~\ref{propcoremorita}, it suffices to show that the map
\begin{equation}\label{eqpfgentwo}
\TB(\fQQ[2](F-,y);\aB_{1};\fQQ[2](x,F-))\to 
\TB(\fQQ[2](-,y);\aB_{2};\fQQ[2](x,-))
\end{equation}
is a weak equivalence for every pair of objects $x,y$ in $\aB_{2}$.
It is clear from Lemmas~\ref{lemNoisNt} and~\ref{lemNNisBN} that
\eqref{eqpfgentwo} is an 
equivalence when either $x$ or $y$ is in the image of $\aB_{1}$.  By
the remarks above, if an object
$y$ in $\aB_{2}$ is isomorphic in
$\pi_{0}\aB_{2}/\pi_{0}\aA_{2}$ to $Fy'$ for some object $y'$ in 
$\aB_{1}$, then $\fQQ[2](-,y)$ is weakly equivalent as a right
$\aB_{2}$-module to $\fQQ[2](-,Fy')$ and the map \eqref{eqpfgentwo} is a weak
equivalence for all $x$.  Since $\pi_{0}\aB_{1}/\pi_{0}\aA_{1}$ is
cofinal in $\pi_{0}\aB_{2}/\pi_{0}\aA_{2}$, for any $y$ in $\aB_{2}$,
there exists $w$ in $\aB_{2}$ such that the sum $w\vee y$ in
$\pi_{0}\aB_{2}/\pi_{0}\aA_{2}$ is isomorphic to $Fz$ for some $z$ in
$\aB_{1}$; then as noted above, the right $\aB_{2}$-module
$\fQQ[2](-,w)\vee \fQQ[2](-,y)$ is 
weakly equivalent to $\fQQ[2](-,Fz)$.  We get compatible weak equivalences
\begin{gather*}
\TB(\fQQ[2](F-,w);\aB_{1};\fQQ[2](x,F-))\ \vee\ 
\TB(\fQQ[2](F-,y);\aB_{1};\fQQ[2](x,F-))\\
\simeq\\
\TB(\fQQ[2](F-,Fz);\aB_{1};\fQQ[2](x,F-))
\end{gather*}
and
\begin{gather*}
\TB(\fQQ[2](-,w);\aB_{2};\fQQ[2](x,-))\ \vee\ 
\TB(\fQQ[2](-,y);\aB_{2};\fQQ[2](x,-))\\
\simeq\\
\TB(\fQQ[2](-,Fz);\aB_{2};\fQQ[2](x,-)),
\end{gather*}
and we see that \eqref{eqpfgentwo}
is a weak equivalence for $x$
and $y$.
\end{proof}

It still remains to prove Lemma~\ref{lemNoisNt}.  The proof makes use
of Bousfield localization \cite{BousfieldLocalization}, \cite[\S
3.3]{Hirschhorn} in the category of right $\aB$-modules for a 
small pretriangulated spectral category $\aB$.  As discussed in
Section~\ref{secdkequiv}, the category $\RMod{\aB}$ of right
$\aB$-modules has a standard compactly generated stable model
structure where the weak equivalences and fibrations are the maps that
are objectwise weak equivalences and fibrations in the stable model
structure on symmetric spectra.  The generating cofibrations and
generating acyclic cofibrations are the maps $\aB(-,b)\sma f$ for $b$
in $\aB$ and $f$ varying through the generating cofibrations and
generating acyclic cofibrations, respectively, of the stable model
structure on symmetric spectra.  The representable right $\aB$-modules
$\aB(-,b)$ are both cofibrant and compact.

Now let $(\aB,\aA)$ be a localization pair.  We say that a right
$\aB$-module $\psi$ is \term{$\aA$-local} if it is fibrant and
$\psi(a)$ is homotopically trivial for every object $a$ of $\aA$.  In
this context, we say that a map of right $\aB$-modules
$f\colon \phi \to \phi'$ is an \term{$\aA$-local equivalence} if it
induces a bijection of morphism sets in
the Quillen homotopy category, $[\phi',\psi]\to [\phi,\psi]$, for every
$\aA$-local right $\aB$-module $\psi$.  The $\aA$-local model
structure on right $\aB$-modules has the same cofibrations as the
standard stable model structure but has weak equivalences the
$\aA$-local equivalences \cite[4.1.2]{Hirschhorn}.  This is a
compactly generated model structure with the acyclic cofibrations
generated by the acyclic cofibrations in the standard stable model
structure together with the maps of the form $\aB(-,a)\sma f$ for $a$
in $\aA$ and $f$ a generating cofibration in the stable model
structure on symmetric spectra.  The fibrant objects in the
$\aA$-local model structure are the $\aA$-local right $\aB$-modules.

More specifically, every acyclic cofibration in the $\aA$-local model
structure is a retract of a sequential colimit of pushouts along
arbitrary coproducts of the generating acyclic cofibrations indicated
above.  The cofiber of such a map in the colimit is weakly equivalent
(in the standard stable model structure) to a wedge of objects of
$\aA$.  Since modules represented by objects of $\aB$ are compact, a
standard argument 
\cite[2.3.17]{HPS}, \cite[2.1]{NeemanBousfield} shows that if a
representable right $\aB$-module $\aB(-,b)$
is $\aA$-acyclic ($\aA$-locally equivalent to $*$), then
in the Quillen homotopy category of the standard stable model
structure, $\aB(-,b)$ is in the thick subcategory generated by the
representables from $\aA$, and so $b$ is in $\aA$.  This implies the
following proposition. 

\begin{prop}\label{proplocmod}
Let $(\aB,\aA)$ be a localization pair.
The spectral Yoneda functor that includes $\aB$ in $\RMod{\aB}$ as the 
representable functors induces a triangulated embedding of
$\pi_{0}\aB/\pi_{0}\aA$ in the Quillen homotopy category of the
$\aA$-local model structure on $\RMod{\aB}$.
\end{prop}

Using the theory reviewed above and the previous proposition, we can now
prove Lemma~\ref{lemNoisNt}.

\begin{proof}[Proof of Lemma~\ref{lemNoisNt}]
Fixing objects $x,y$ in $\aB_{1}$, it suffices to show that the map
$\fQQ[1](x,y)\to \fQQ[2](Fx,Fy)$ is a weak equivalence.  

We take advantage of the functoriality of $\fQQ[i]$ generalized to
modules and the previous proposition.  Choose an $\aA_{1}$-local
$\aB_{1}$-module $\psi_{1}$ and an $\aA_{1}$-local acyclic cofibration
$q\colon \aB_{1}(-,y)\to \psi_{1}$.  It is clear from the
characterization of the generating $\aA_{1}$-local acyclic
cofibrations that $q$ induces a weak equivalence
$\fQQ[1](x,y)\to \fQQ[1](x,\psi_{1})$. Moreover, since $\psi_{1}(a)$
is homotopically trivial for every object $a$ in $\aA_{1}$, we have that the map
$\psi_{1}(x)\to \fQQ[1](x,\psi_{1})$ is a weak equivalence.

The functor $F^{*}$ from right $\aB_{2}$-modules to right
$\aB_{1}$-modules has a left adjoint $\Lan_{F}$ defined by left Kan
extension.  Since $\Lan_{F}$ takes $\aB_{1}(-,b)$ to $\aB_{2}(-,Fb)$
for any object $b$ in $\aB_{1}$, $\Lan_{F}$ takes the generating
cofibrations and generating acyclic cofibrations of the
$\aA_{1}$-local model structure into the generating cofibrations and
generating acyclic cofibrations of the $\aA_{2}$-local model
structure, i.e., $\Lan_{F},F^{*}$ is a Quillen adjunction on the local
model structures.  In particular, $\Lan_{F}$ takes $q$ to an
$\aA_{2}$-local acyclic cofibration $\aB_{2}(-,Fy)\to
\Lan_{F}\psi_{1}$.  Choose an $\aA_{2}$-local object $\psi_{2}$ and an
$\aA_{2}$-local acyclic cofibration $\Lan_{F}\psi_{1}\to \psi_{2}$.
Now we have weak equivalences
\[
\fQQ[2](Fx,Fy)\to
\fQQ[2](Fx,\Lan_{F}\psi_{1})\to \fQQ[2](Fx,\psi_{2}).
\]
Moreover, since $\psi_{2}(a)$ is homotopically trivial for every object $a$ in
$\aA_{2}$, we have that the map $\psi_{2}(Fx)\to \fQQ[2](Fx,\psi_{2})$
is a weak equivalence.

Applying Proposition~\ref{proplocmod} and the hypothesis that
$\aB_{1}\to \aB_{2}$ induces an embedding of
$\pi_{0}\aB_{1}/\pi_{0}\aA_{1}$ into $\pi_{0}\aB_{2}/\pi_{0}\aA_{2}$,
we see that $\psi_{1}(b)\to \psi_{2}(Fb)$ is a weak equivalence for
every object $b$ in $\aB_{1}$.  Thus, we have shown that in the
commutative diagram
\[
\xymatrix{%
\psi_{1}(x)\ar[r]^-{\simeq}\ar[d]_-{\simeq}
&\fQQ[1](x,\psi_{1})\ar[d]
&\fQQ[1](x,y)\ar[l]_-{\simeq}\ar[d]\\
\psi_{2}(Fx)\ar[r]_-{\simeq}
&\fQQ[2](Fx,\psi_{2})
&\fQQ[2](Fx,Fy),\ar[l]^-{\simeq}
}
\]
the arrows marked ``$\simeq$'' are weak equivalences.  It follows that
the map $\fQQ[1](x,y)\to \fQQ[2](Fx,Fy)$ is a weak equivalence.
\end{proof}

\section{Applications of the general localization theorem}\label{secappl}

We now turn to the applications of the general theory of the preceding
sections to $THH$ and $TC$ of schemes.  We begin with a discussion of
the spectral enrichment of the derived category of a scheme.  Recent
work shows that any stable category can be regarded as enriched in
symmetric spectra \cite{Dugger,DuggerShipleyEnriched,SSStable} and one
approach would be to apply this theory in the setting of categories of
unbounded complexes to construct a spectral derived category from
first principles.  On the other hand, such an approach would
demand a comparison with the DG-category structures that arise in
nature on categories of complexes.  For this reason, we take
the simpler approach of lifting DG-categories to associated
spectral categories.

For a scheme $X$, let $\aK^{DG}(X)$ denote the pretriangulated
DG-category of unbounded (cohomologically graded) complexes of sheaves
of $\sO_{X}$-modules; 
its homotopy category $\pi_{0}\aK^{DG}(X)$ is the triangulated
category typically denoted $\aK(X)$ of unbounded complexes and chain
homotopy classes of maps.  The derived category $\aD(X)$ is the
localization of $\aK(X)$ obtained by inverting the quasi-isomorphisms.
The derived category of perfect complexes $\aD_{\parf}(X)$ is the full
triangulated subcategory of unbounded complexes locally
quasi-isomorphic to strictly bounded complexes of vector bundles.  By
choosing a large enough cardinal $\aleph$ and restricting to the
perfect complexes whose underlying sets are in $\aleph$, we can find a
small full pretriangulated subcategory $\aK^{DG}_{\parf}(X)$ of
$\aK^{DG}(X)$ consisting of perfect complexes and having the property
that the triangulated quotient of the homotopy category by the full
subcategory of acyclics is equivalent to $\aD_{\parf}(X)$ via the
canonical map.  Moreover, when $X$ is quasi-compact and
quasi-separated, the full subcategory $\aK^{DG}_{\parf}(X)_{\flat}$
consisting of those complexes in $\aK^{DG}_{\parf}(X)$ that are
strictly bounded above and degreewise flat $\sO_{X}$-modules also has
the property that the triangulated quotient of the homotopy category
by the full subcategory of acyclics is equivalent to $\aD_{\parf}(X)$
via the canonical map.

Keller \cite{KellerCyclic} and Drinfeld \cite{DrinfeldDG} described
``quotient'' DG-categories whose homotopy categories model the
quotients of triangulated categories.  We obtain small DG-categories
$\DGparf$ and (for $X$ quasi-compact and quasi-separated) $\DGparflat$
whose homotopy categories are equivalent to the derived category
$\aD_{\parf}$.  We obtain small spectral categories $\DSparf$ and
$\DSparflat$ associated to $\DGparf$ and $\DGparflat$.  We now prove
Theorems~\ref{inttt}, \ref{intgh}, and~\ref{intbl} from the
introduction; in all cases, the results for $TR$ and $TC$ follow from
the corresponding results for $THH$ by Proposition~\ref{proplazy}.


\begin{proof}[Proof of Theorem~\ref{inttt}]
The theorem follows from Theorem~\ref{thmgenone}: For the
first statement, the Thomason-Trobaugh localization sequence, we apply 
Theorem~\ref{thmgenone} with $\aB=\DSparf$, $\aA$ the full spectral
subcategory of $\DSparf$ consisting of those complexes that are
supported on $X-U$, and $\aC=\DSparf[U]$, using a lift $\aB\to \aC$ of
the DG-functor $j^{*}\colon \DGparf\to \DGparf[U]$.   
The Mayer-Vietoris statement follows from the 
localization statement and Corollary~\ref{cortriang} since the
inclusion in $X$ of any open set $V$ containing $X-U$ induces an
equivalence on the derived categories of perfect complexes supported
on $X-U=V-U\cap V$.
\end{proof}


\begin{proof}[Proof of Theorem~\ref{intgh}]
We choose an affine open cover 
$\{U_{1},\dotsc,U_{r}\}$ of 
$X$.  For each $i_{1},\dotsc,i_{n}$ let
$U_{i_{1},\dotsc,i_{n}}=U_{i_{1}}\cap \dotsb \cap U_{i_{n}}$ and let
$A_{i_{1},\dotsc,i_{n}}=\sO_{U_{i_{1},\dotsc,i_{n}}}$.  Since $X$ is
semi-separated, without loss of generality, 
$U_{i_{1},\dotsc,i_{n}}=\Spec A_{i_{1},\dotsc,i_{n}}$.  We now
construct a \v{C}ech complex on $THH$ associated to this cover as 
follows.

Let $\aA$ denote the full subcategory of $\aK^{DG}_{\parf}(X)$
consisting of the acyclic complexes, and for each (non-empty)
$U_{i_{1},\dotsc,i_{n}}$, let $\aA_{i_{1},\dotsc,i_{n}}$ denote the
full subcategory of $\aK^{DG}_{\parf}(X)$ of objects acyclic on
$U_{i_{1}, \dotsc ,i_{n}}$.  We have defined $\DGparf$ by Drinfeld's
quotient category construction $\aK^{DG}_{\parf}(X)/\aA$.  For all
$i_{1},\dotsc,i_{n}$ and $U=U_{i_{1}}\cap\dotsb\cap U_{i_{n}}$, the
DG-functor $j^{*}\colon \aK^{DG}_{\parf}(X)\to \aK^{DG}_{\parf}(U)$
(associated to $j\colon U\subset X$) induces a DG-functor
$\aK^{DG}_{\parf}(X)/\aA_{i_{1},\dotsc,i_{n}}\to \DGparf[U]$ that is a
DG-equivalence onto its image.  Moreover, this 
functor is cofinal in that every perfect complex on $U$ is a direct
summand of $j^{*}$ of a perfect complex on $X$.
We apply functorial factorization to construct the associated
(fibrant) spectral categories $\aD_{i_{1},\dotsc,i_{n}}$.  This
constructs a strictly commuting diagram of spectral functors
associated to intersections of the open sets in the cover.  Moreover,
comparing this construction with the argument above for
Theorem~\ref{inttt}, we see that the map
\[
THH(\DSparf)\to \holim_{\aS_{r}} THH(\aD_{i_{1},\dotsc,i_{n}})
\]
is a weak equivalence, where $\aS$ is the partially-ordered set of
non-empty subsets of $1,\dotsc,r$. 

Each of the categories $\aD_{i_{1},\dotsc,i_{n}}$ has an object called
$\sO_{X}$ whose endomorphism spectrum is an Eilenberg-Mac Lane ring
spectrum for $A_{i_{1},\dotsc,i_{n}}$.  Write $HA_{i_{1},\dotsc,i_{n}}$
for 
$\aD_{i_{1},\dotsc,i_{n}}(\sO_{X},\sO_{X})$.  Since the objects
called $\sO_{X}$ in $\aD_{i_{1},\dotsc,i_{n}}$ are compatible under
inclusion of intersections, we obtain a map
\[
\holim_{\aS_{r}} THH(HA_{i_{1},\dotsc,i_{n}}) \to \holim_{\aS_{r}}
THH(\aD_{i_{1},\dotsc,i_{n}}).
\]
The lefthand spectrum is easily seen to be equivalent to the \v{C}ech
cohomology associated to the cover $\{ U_{1},\dotsc,U_{r} \}$ of the
Zariski presheaf of symmetric spectra 
$THH(\sO_{(-)})$.
Geisser and Hesselholt \cite[3.2.1]{GeisserHesselholt} showed that the
homotopy groups of $THH$ of a commutative ring form a quasi-coherent
sheaf, and so the lefthand homotopy limit computes the hypercohomology
spectrum of $THH(\sO_{(-)})$, i.e., $THH(X)$ as defined by
\cite[3.2.3]{GeisserHesselholt}.

Thus, Theorem~\ref{intgh} for quasi-compact semi-separated schemes
reduces to showing that each map 
$THH(HA_{i_{1},\dotsc,i_{n}})\to THH(\aD_{i_{1},\dotsc,i_{n}})$ is a
weak equivalence.  This follows from Theorem~\ref{thmthick}.
\end{proof}


\begin{proof}[Proof of Theorem~\ref{intbl}]  
We have functors  
\begin{gather*}
Li^{*}\colon \DGparflat[X]\to \DGparflat[X']\\
Lj^{*}\colon \DGparflat[Y]\to \DGparflat[Y'].
\end{gather*}
Each of these is a
DG-equivalence to its image.  Let $\aB_{1}=\DSparflat[X']$ and 
$\aB_{2}=\DSparflat[Y']$, and let
$\aA_{1}$ denote the full spectral subcategory of $\aB_{1}$ consisting
of objects equivalent to those in the image of $Li^{*}$ and $\aA_{2}$
denote the full spectral subcategory of $\aB_{2}$ of objects in the
image of $Lj^{*}$.  The map $Lp^{*}$ lifts to a map $\DSparflat[X]\to
\DSparflat[X']$ that lands in $\aA_{1}$ and is a DK-equivalence to
$\aA_{1}$.  Likewise, $Lq^{*}$ induces a DK-equivalence of
$\DSparf[Y']$ with $\aA_{2}$.  In this way, we obtain a strictly
commuting spectral model for the DG-functors $Lp^{*}$, $Lq^{*}$,
$Li^{*}$, and $Lj^{*}$ as a map of localization pairs
$(\aB_{1},\aA_{1})\to (\aB_{2},\aA_{2})$.  By \cite[\S
2.7]{ThomasonBlowUp} or \cite[1.5]{CHSW}, this map induces an
equivalence on quotient triangulated categories, and therefore a weak
equivalence $\CTHH(\aB_{1}/\aA_{1})\to \CTHH(\aB_{2}/\aA_{2})$ by
Theorem~\ref{thmgentwo}.  Corollary~\ref{corgentwo} gives the
corresponding result on $TC$.  Theorem~\ref{intbl} now follows.
\end{proof}


\begin{proof}[Proof of Theorem~\ref{intproj}]
Let $\pi \colon \bP \sE_X \to X$ be the projective bundle of
an algebraic vector bundle $\sE$ of rank $r$.
Thomason \cite[\S 2.7]{ThomasonBlowUp} constructed a triangulated
filtration  
\[
0\htp \aA_{r}\subset \aA_{r-1} \subset \dotsb \subset \aA_{0} =
\aD_{\parf}(\bP\sE_{X})
\]
of the derived category as follows.  Let $\aA_k$ denote the
full subcategory of $\aD_{\parf}(\bP \sE_X)$ consisting of complexes
$Z$ such that $R\pi_*(Z \otimes \sO_{\bP \sE_X}(i)) = 0$ for $0 \leq i
< k$.  By definition, $\aA_0 = \DGparf(\bP \sE_X)$, and since $\sE$ is
rank $r$, $\aA_r$ consists of acyclic complexes
\cite[2.5]{ThomasonBlowUp}.  Furthermore, 
$\aA_k$ admits the alternate description as the thick subcategory of
$\aD_{\parf}(\bP \sE_X)$ 
generated by $L\pi^* (-) \otimes \sO_{\bP \sE_X}(-j)$ for $k \leq j <
r$ (see also \cite[1.2]{CHSW}).  Let
$\aA'_{r-1},\dotsc,\aA'_{0}$ denote the corresponding filtration on
$\DSparf[\bP\sE_{X}]$.

The functor $L\pi^{*}(-)\otimes \sO_{\bP\sE_{X}}(-k)$ from
$\aD_{\parf}(X)$ to $\aA_{k}$ admits a refinement to a DG-functor
$\DGparf\to \DGparf[\bP\sE_{X}]$, which we can lift to a spectral
functor $\DSparf\to \aA'_{k}$.  Viewed as map of localization pairs
$(\DSparf,0)\to (\aA'_{k},\aA'_{k+1})$, the induced map of
triangulated quotients $\aD_{\parf}(X)\to \aA_{k}/\aA_{k+1}$ is an
equivalence \cite[\S 2.7]{ThomasonBlowUp}.  Theorem~\ref{thmgentwo}
then shows that the induced map $THH(X)\to
\CTHH(\aA'_{k},\aA'_{k+1})$ is a weak equivalence.  In particular, we
obtain split cofibration sequences
\[
THH(\aA^{\prime}_{k+1}) \to THH(\aA^{\prime}_k) \to
THH(X), 
\]
and hence weak equivalences
\[
THH(\aA^{\prime}_k) \htp THH(\aA^{\prime}_{k+1}) \times
THH(X).
\]
for $k=0,\dotsc,r-1$. Using Corollary~\ref{corgentwo}, we get the
corresponding results for $TC$. This completes the proof of 
Theorem~\ref{intproj}.
\end{proof}

\section{The cyclotomic trace from $K^{B}$}
\label{secnoncon}

In this section we show that the cyclotomic trace from $K$-theory to
$THH$ and $TC$ factors through Thomason-Trobaugh's construction of
Bass' non-connective $K$-theory $K^{B}$ \cite[\S6]{ThomasonTrobaugh}.
Using a version of Bass' fundamental theorem for $THH$, we factor
the cyclotomic trace map from connective $K$-theory on affine schemes
(commutative rings) 
through $K^{B}$.  This factorization holds more generally for maps from
$K$-theory to any theory satisfying the appropriate analogue of Bass'
fundamental theorem, and is natural for such functors to a (strict
point-set) category of spectra.  Since the trace map admits
a model in which it is a map of presheaves restricted to affine covers, we
obtain the factorization $K^{B}\to THH$ on the level of presheaves,
which we show lifts to a map of presheaves $K^{B}\to TC$.  For
quasi-compact semi-separated schemes, $K^{B}$ is equivalent to the
\v{C}ech 
hypercohomology spectrum of its presheaf by
\cite[8.4]{ThomasonTrobaugh}. 
The work of the previous section (and \cite{GeisserHesselholt}) shows
that for such schemes $THH$ and $TC$ are each equivalent to both the
hypercohomology spectrum and \v{C}ech hypercohomology spectrum of
their respective 
presheaves.  This then constructs the trace map $K^{B}\to TC$ for
all quasi-compact, semi-separated schemes.

We begin by discussing the analogue of Bass' fundamental theorem that
we need.  For the purposes of this section, we say that a covariant
functor $F$ from commutative rings to some (point-set) category of
spectra is a \term{Bass functor} when it comes with a natural
transformation $\tau\colon \Sigma F(R)\to
F(R[t,t^{-1}])$ and satisfies:
\begin{enumerate}
\item For any $R$, $F(R)$ is connective ($\pi_{n}F(R)=0$ for $n<0$),
\item For any $R$ and any $n\geq 0$, the sequence
\[
0\to \pi_{n}F(R)\to \pi_{n}F(R[t])\oplus \pi_{n}F(R[t^{-1}])\to
\pi_{n}F(R[t,t^{-1}])
\]
induced by the inclusion maps is exact, and
\item For any $R$ and any $n>0$, the composite map 
\begin{multline*}
\qquad \pi_{n-1}F(R)=\pi_{n}\Sigma F(R)\to \pi_{n}F(R[t,t^{-1}])\\
\to \Coker\big(\pi_{n}F(R[t])\oplus \pi_{n}F(R[t^{-1}])\to
\pi_{n}F(R[t,t^{-1}])\big)
\end{multline*}
induced by $\tau$ is an isomorphism.
\end{enumerate}
A map of Bass functors is a natural transformation $F\to G$ that
commutes with the maps $\tau$.  The key fact we need to apply this
theory is the following (well-known) theorem, whose proof we review at
the end of this section. 

\begin{thm}\label{thmtracebass}
The $K$-theory functor and the $THH$ functor admit models that are
Bass functors with the trace map a map of Bass functors.
\end{thm}

As an immediate consequence of the definition, a Bass functor in
particular comes 
with a natural $4$-term exact sequences
\[
0\to \pi_{n}F(R)\rightarrow \pi_{n}F(R[t])\oplus \pi_{n}F(R[t^{-1}])\rightarrow
\pi_{n}F(R[t,t^{-1}])\rightarrow \pi_{n-1}F(R)\to 0
\]
for $n>0$ with the map
$\pi_{n}F(R[t,t^{-1}])\rightarrow \pi_{n-1}F(R)$
naturally split.  This exact sequence
and splitting are functorial in maps of Bass functors.  Bass'
construction extends these sequences to all $n$:

\begin{defn}[Bass' Construction]
For a Bass functor $F$, let $\beta_{n}F=\pi_{n}F$ and let
\[
\tau_{n}\colon \beta_{n}F(R)\to \beta_{n+1}F(R[t,t^{-1}])
\]
be the map induced by $\tau$ for $n\geq 0$.  Inductively, for $n\leq 0$,
define 
\[
\beta_{n-1}F(R)=\Coker\big(\beta_{n}F(R[x])\oplus \beta_{n}F(R[x^{-1}])\to
\beta_{n}F(R[x,x^{-1}])\big)
\]
and $\tau_{n}\colon \beta_{n-1}F(R)
\to \beta_{n}F(R[t,t^{-1}])$ to be the induced map on cokernels
\[
\Coker\left(\vcenter{\xymatrix@R=2pt{%
\beta_{n}F(R[x])\\
\oplus \\
\beta_{n}F(R[x^{-1}])\ar[dd]\\\relax\mathstrut\\
\beta_{n}F(R[x,x^{-1}])}}\right)
\to
\Coker\left(\vcenter{\xymatrix@R=2pt{%
\beta_{n+1}F(R[x,t,t^{-1}])\\
\oplus\\
\beta_{n+1}F(R[x^{-1},t,t^{-1}])\ar[dd]\\\relax\mathstrut\\
\beta_{n+1}F(R[x,x^{-1},t,t^{-1}])}}
\right)
\]
\end{defn}

Applied to the $K$-theory functor, Bass' construction defines Bass' negative
$K$-groups.  Applied to the $THH$ functor, $\beta_{n}THH=0$ for $n<0$
since the map 
\[
\pi_{0}THH(R[x])\oplus \pi_{0}THH(R[x^{-1}])\to \pi_{0}THH(R[x,x^{-1}])
\]
is surjective.  (It is the map $R[x]\oplus R[x^{-1}]\to R[x,x^{-1}]$.)
Thomason and Trobaugh extended Bass'
construction to a construction on spectra suitable for application to
general Bass functors as defined above.  The following is essentially
a simplification of \cite[6.3]{ThomasonTrobaugh}. 

\begin{lem}\label{lemTT}
Let $F$ be a Bass functor.  There exists a functor $F^{B}$ from
commutative rings to spectra and a natural transformation 
$F\to F^{B}$ that is an isomorphism on $\pi_{n}$ for
$n\geq 0$ and induces (as indicated in \cite[6.3]{ThomasonTrobaugh})
a natural isomorphism $\beta_{n}F\to \pi_{n}F^{B}$ for $n<0$.  The
functor and natural transformation are functorial in maps of Bass functors.
\end{lem}

We need a few of the details of the construction.  Thomason and Trobaugh
construct $F^{B}$ as the homotopy colimit of a sequence of functors 
\[
F=F_{0}=F'_{0}\to F'_{-1}\to F'_{-2}\to \dotsb F'_{-k}\to \dotsb.
\]
The functor $F'_{-k-1}$ is formed inductively as a homotopy pushout 
\[
\xymatrix{
\Omega^{k} F_{-k} \ar[r] \ar[d] & F'_{-k} \\
\Omega^{k+1} F_{-k-1}
}
\]
for functors $F_{-k}$ which come with natural transformations $\Sigma
F_{-k}\to F_{-k-1}$.  The functor $F_{-k-1}$ is defined inductively
as the homotopy cofiber of the natural map
\[
F_{-k}(R[x])\cuph_{F_{-k}(R)} F_{-k}(R[x^{-1}])\to F_{-k}(R[x,x^{-1}])
\]
(where ``$\cuph$'' denotes the homotopy pushout).
The map $\Sigma F_{-k}\to F_{-k-1}$ comes from the canonical map
$\Sigma F_{-k}(R)\to F_{-k}(R[t,t^{-1}])$, constructed just as in Bass'
construction in algebra, as the induced map on cofibers coming from the
natural commutative diagram
\[
\xymatrix{%
\Sigma F_{-k}(R[x])\cuph_{\Sigma F_{-k}(R)} \Sigma F_{-k}(R[x^{-1}])\ar[r]\ar[d] 
&\Sigma F_{-k}(R[x,x^{-1}])\ar[d]\\
F_{-k}(R[t,t^{-1},x])\cuph_{F_{-k}(R[t,t^{-1}])} F_{-k}(R[t,t^{-1},x^{-1}])\ar[r]
&F_{-k}(R[t,t^{-1},x,x^{-1}]).
}
\]
Our notation differs slightly from that of
\cite[6.3]{ThomasonTrobaugh}; our $F_{-k}$ is their $\Sigma^{k}F^{-k}$.
As a consequence of this construction, we get the following observation.

\begin{prop}\label{propbasscyc}
If $F$ is a Bass functor and factors through cyclotomic spectra (with
$\tau$ a natural map of cyclotomic spectra), then
the functor $F^{B}$ factors through cyclotomic spectra and the natural transformation 
$F\to F^{B}$ is a natural transformation of cyclotomic spectra.
\end{prop}

Combining the Thomason-Trobaugh lemma with
Theorem~\ref{thmtracebass}, we get a natural transformation of functors
\[
K^{B} \to THH^{B} \htp THH.
\]
The model of $THH$ constructed below and referred to in
Theorem~\ref{thmtracebass} satisfies the hypotheses of
Proposition~\ref{propbasscyc}.  Then $THH^{B}$ is a cyclotomic
spectrum and we can form a functor $TC$ as the appropriate limit (or
pro-object).  We obtain the following commutative
diagram of functors.
\[
\xymatrix{%
K^{B}\ar[r]&TC^{B}\ar[d]&TC\ar[l]_{\htp}\ar[d]\\
&THH^{B}&THH\ar[l]^{\htp}
}
\]
This extends the
trace to non-connective $K$-theory.

The remainder of the section proves Theorem~\ref{thmtracebass}.
We begin by observing that $THH$ algebraically satisfies the analogue
of Bass' fundamental theorem.  For a commutative ring $R$, the
Eilenberg Mac\,Lane spectrum $HR$ is a commutative ring spectrum in any of the
modern categories of spectra.  We have a weak equivalence of
associative ring spectra $HR \sma TS\to HR[t]$, where $TS$ is the free
associative ring spectrum on the sphere
spectrum (or a cofibrant model of it).  We then get a weak equivalences
of $THH(R)$-modules
\begin{multline*}
THH(R[t]) \htp THH(R \sma TS) \htp THH(R) \sma THH(TS)\\ 
\iso THH(R)\sma_{HR}(HR \sma THH(TS)).
\end{multline*}
We also have the identifications 
\[
HR\sma THH(TS) \iso THH^{HR}(HR \sma TS)\htp THH^{HR}(R[t]),
\]
where $THH^{HR}(R[t])$ is as defined in \cite[\S IX.1.7]{EKMM} and is 
essentially the spectrum whose homotopy groups are
$HH^{R}_{*}(R[t])$.  Since this Hochschild homology is a free module
over $R$, we obtain the computation 
\[
\pi_{*}THH(R[t]) \iso \pi_{*}THH(R) \otimes_{R} HH^{R}_{*}(R[t])
\iso \pi_{*}THH(R) \otimes_{R} R[t] \langle 1,\sigma_{t}\rangle,
\]
where $1$ is in degree zero and $\sigma_{t}$ is in degree one.
This is an isomorphism of $\pi_{*}THH(R)$-modules, and is natural in $R$ and
$TS$, though not obviously in $R[t]$.

Writing $TS[t^{-1}]$ for the localization of $TS$ under multiplication
by the generator of $\pi_{0}S$ (which we are thinking of as $t$), we
have a weak equivalence of associative ring spectra $HR\sma
TS[t^{-1}]\to HR[t,t^{-1}]$, and as above, we get the weak equivalences
\[
THH(R[t,t^{-1}]) \htp THH(R \sma TS[t^{-1}]) \htp 
THH(R)\sma_{HR} THH^{HR}(HR[t,t^{-1}])
\] 
and the computation
\begin{multline*}
\pi_{*}THH(R[t,t^{-1}]) \iso \pi_{*}THH(R) \otimes_{R} HH^{R}_{*}(R[t,t^{-1}])\\
\iso \pi_{*}THH(R) \otimes_{R} R[t,t^{-1}] \langle 1,\sigma_{t}\rangle,
\end{multline*}
Again, this is an isomorphism of $\pi_{*}THH(R)$-modules, and is
natural in $R$ and $TS[t^{-1}]$, though not obviously in
$R[t,t^{-1}]$.  

The map $R[x]\to R[t,t^{-1}]$ sending $x$ to $t^{-1}$ is induced by a
map of associative ring spectra $TS\to TS[t^{-1}]$, namely, the map
induced by the map $S\to TS[t^{-1}]$ sending the generator of
$\pi_{0}S$ to $t^{-1}$ in $\pi_{0}TS[t^{-1}]\iso \bZ[t,t^{-1}]$.
Thus, we can compute the maps in Bass' sequence for $\pi_{*}THH$ in
terms of Hochschild homology.  This then becomes an easy computation
with resolutions: the inclusion $R[t^{-1}]\to R[t,t^{-1}]$ induces the
map of $\pi_{*}THH(R)\otimes_{R}R[t^{-1}]$-modules that sends $1$ to
$1$ and $\sigma_{t^{-1}}$ to $-t^{-2}\sigma_{t}$.  It follows that the
sequence of graded abelian groups
\[
0 \to \pi_{*}THH(R) \to \pi_{*}THH(R[t])\oplus \pi_{*}THH(R[t^{-1}])
\to \pi_{*}THH(R[t,t^{-1}])
\]
is exact and the map $\pi_{*-1}THH(R)\to \pi_{*}THH(R[t,t^{-1}])$
induced by the inclusion and multiplication by $t^{-1}\sigma_{t}$
induces an isomorphism from $\pi_{*-1}THH(R)$ onto the cokernel of the
last map above.

Thomason and Trobaugh \cite[\S6]{ThomasonTrobaugh} prove an analogous
formulation of Bass' fundamental theorem for $K$-theory: The three
term sequence is exact and the map $K_{n-1}R\to K_{n}R[t,t^{-1}]$
induced by the inclusion and multiplication by $t$ (viewed as an
element of 
$K_{1}R[t,t^{-1}]$) induces an isomorphism onto the cokernel for
$n>1$.  Since the Dennis trace map takes the element $t$ in
$K_{1}(\bZ[t,t^{-1}])$ to the element $t\sigma_{t^{-1}} =
-t^{-1}\sigma_{t}$ in $HH_{1}(\bZ[t,t^{-1}])$, multiplication by the
image of $t$ under the trace to $THH$ also provides an isomorphism
from $\pi_{*-1}THH(R)$ to the cokernel for $\pi_{*}THH$.  We now have
what we need to prove Theorem~\ref{thmtracebass}.

\begin{proof}[Proof of Theorem~\ref{thmtracebass}]
We give an argument that uses minimal details of the construction of
the cyclotomic trace.  The next section contains a review of the
cyclotomic trace and a more direct construction of the extension to
non-connective $K$-theory.  
For this argument, we use the model of the trace map described in
\cite[\S2.1.6]{DundasMcCarthy}, with some of the modifications in
\cite[\S6.3]{GeisserHesselholt} that involve the multiplicative
structure.  We regard $K$ as a functor from exact categories to
symmetric spectra, and we regard $TH$ (in the notation of
\cite[\S6.3]{GeisserHesselholt}) as a functor from exact categories to
symmetric spectra of cyclotomic spectra.  The
point-set category of spectra we work in is the
category of symmetric spectra of orthogonal spectra, and we use the
free functor in the orthogonal spectrum direction
$F_{0}^{\orth}$ (and geometric realization) to convert symmetric
spectra to symmetric spectra of orthogonal spectra.

Since $K$ and $TH$ are functors on exact categories, to get functors
on commutative rings, we need a model of the exact category of finitely
generated projective modules that is strictly functorial in maps of
commutative rings.  For this, consider the category $\aP(R)$ whose
objects are pairs $(P,m)$ where $P$ is a projective submodule of
$R^{m}$, and whose maps $(P,m)\to (Q,n)$ are the $R$-module maps $P\to
Q$. This is an exact category in the evident way.  A map of
rings $R\to R'$ induces a map $\aP(R)\to \aP(R')$ by
extension of scalars and the canonical identification
$R'\otimes_{R}R^{m}\iso R^{\prime m}$; this makes $\aP$ a functor from
commutative rings to exact categories.  Defining
$K(R)=K(\aP(R))$ and $THH(R)=TH(\aP(R))$ (in the notation of
\cite{GeisserHesselholt}), we obtain functors from commutative rings
to symmetric spectra and from commutative rings to symmetric spectra
of cyclotomic spectra, respectively.  We also obtain a natural
transformation of symmetric spectra of orthogonal spectra
$F_{0}^{\orth}F_{0}K(R)\to THH(R)$.

We have a bi-exact strictly associative tensor product on $\aP(R)$
defined by the usual tensor product and the (lexicographical order)
identification $R^{m}\otimes_{R} R^{n}\iso R^{mn}$.  As observed in
\cite[\S6.3]{GeisserHesselholt}, it follows that $K(R)$ is naturally
an associative ring symmetric spectrum.  Moreover, using the
associative smash product pairing of a space and a cyclotomic spectrum,
$THH(R)$ becomes a module over $K(R)$ in the category of symmetric
spectra of cyclotomic spectra.  If it were possible, what we
would like to do as a next step is choose a map $\bar t$ from $S^{1}$
to the zeroth space of $K(\bZ[t,t^{-1}])$ representing the element $t$
of $K_{1}(\bZ[t,t^{-1}])$.  We would then get natural point-set maps
of symmetric spectra $\tau \colon \Sigma K(R)\to K(R[t,t^{-1}])$ using
the inclusions and multiplication 
\[
K(R)\sma K(\bZ[t,t^{-1}])\to K(R[t,t^{-1}])\sma K(R[t,t^{-1}])
\to K(R[t,t^{-1}]).
\]
However, although $t$ is represented in the first space of
$K(\bZ[t,t^{-1}])$, it is not represented in the zeroth space.  To fix
this, we use a variant of the functor $M$ of \cite[3.2.1]{ShipleyD}.

For a symmetric spectrum $X$, let $MX$ be the symmetric spectrum of
topological spaces
\[
MX=\hocolim_{\bn\in \aI}\Omega^{n}|sh_{n}X|
\]
where $sh_{n}X$ is the symmetric spectrum $sh_{n}X_{k}=X_{k+n}$
\cite[2.2.12]{HSS}.  The inclusion of $X$ as $\Omega^{0}sh_{0}X$
induces a natural transformation $X\to MX$.  When $X$ is a positive
$\Omega$-spectrum like $K(R)$, $MX$ is an $\Omega$-spectrum and $X\to
MX$ is a positive level equivalence, cf.~\cite[p.~168]{ShipleyD}.
Moreover, $M$ is a monoidal functor with
\begin{multline*}
MX\sma MY \iso \hocolim_{(\bm,\bn)\in \aI^{2}}
(\Omega^{m}|sh_{m}X| \sma \Omega^{n}|sh_{n}Y|)\\
\to \hocolim_{\bp\in \aI}\Omega^{p}(|sh_{p}(X\sma Y)|)
= M(X\sma Y)
\end{multline*}
induced by disjoint union functor $\aI\times \aI\to \aI$, the map
\[
\Omega^{m}|X_{m+k}| \sma \Omega^{n}|Y_{n+\ell}|
\to \Omega^{m+n}|X_{m+k}\sma Y_{n+\ell}|
\to \Omega^{m+n}|(X\sma Y)_{m+k+n+\ell}|,
\]
and the appropriate permutation $(X\sma Y)_{m+k+n+\ell}\iso (X\sma
Y)_{m+n+k+\ell}$ on the $(m+n+k+\ell)$-th space of $X\sma Y$.
We have the analogous endofunctor $M$ and natural transformation $X\to
MX$ on the category of symmetric spectra
of cyclotomic spectra.  The cyclotomic spectra $THH(R)$ are
$\Omega$-spectra in the symmetric spectrum direction, and so the map
$THH(R)\to MTHH(R)$ is a level equivalence in the symmetric spectrum
direction. 

Now we choose a point-set map $\bar t$ from $S^{1}$ to the
zeroth space of $MK(\bZ[t,t^{-1}])$.
We define $\tau \colon \Sigma MK(R)\to MK(R[t,t^{-1}])$ to be the
natural transformation of symmetric spectra induced by multiplication
with the point-set 
representative $\bar t$ of $t$.  Likewise, we define $\tau \colon
\Sigma MTHH(R)\to MTHH(R[t,t^{-1}])$ to be the natural transformation
of symmetric spectra of cyclotomic spectra
induced by multiplication by $\bar t$.  Using $F_{0}^{\orth}MK(R)$
and $MTHH(R)$ as our models for $K(R)$ and $THH(R)$, this constructs
$K$ and $THH$ as Bass functors, with $THH$ a Bass functor in
a point-set category of cyclotomic spectra, and the cyclotomic trace a
natural transformation of Bass functors.
\end{proof}


\section{The cyclotomic trace for DG-Waldhausen categories}
\label{apptrace}

In Section~\ref{secappl}, we implicitly constructed the cyclotomic trace
connecting the $K$-theory of a scheme to the $TC$ and $THH$ of the
associated spectral derived category via the comparison to the
Geisser-Hesselholt definition of these spectra in terms of
hypercohomology.  This streamlined approach allowed us to avoid the
lengthy technical development necessary for a more intrinsic
construction of the cyclotomic trace, and was sufficient for our
applications.  In this section, we complete the theory of
$TC$ and $THH$ of spectral derived categories by describing an
intrinsic construction of the cyclotomic trace.

Our construction of the cyclotomic trace follows the
perspective of \cite[\S 2.1.6]{DundasMcCarthy} that the trace should be
regarded as ``the inclusion of 
the objects'' from a Waldhausen category to a model of $THH$ which
mixes the cyclic bar construction and Waldhausen's $\Sdot$
construction.  In order to enable this mixing, we work with a class of
Waldhausen categories equipped with a DG-enrichment that is compatible
with the Waldhausen structure.  We call these DG-Waldhausen
categories; they are in particular complicial biWaldhausen categories
\cite[1.2.11]{ThomasonTrobaugh}. 

\begin{defn}\label{defdgwald}
A \term{DG-Waldhausen category} consists of 
a small full subcategory $\aC$ of the category of complexes of an abelian
category $Ab_{\aC}$ (which is part of the structure), and a
subcategory $w\aC$ of $\aC$ called the weak equivalences, satisfying
the following properties.
\begin{enumerate}
\item $\aC$ contains zero.
\item $\aC$ is closed under pushouts along degreewise split
monomorphisms and pullbacks along degreewise split epimorphisms.
\item $\aC$ is closed under cones and cocones.
\item The weak equivalences contain the quasi-isomorphisms of
complexes, are preserved by pushout along degreewise-split
monomorphisms and pullback along degreewise-split epimorphisms, and
satisfy Waldhausen's saturation and extension properties. 
\end{enumerate}
A \term{DG-exact functor} from $(\aC,Ab_{\aC},w\aC)$ to
$(\aC',Ab_{\aC'},w\aC')$ is an
additive functor $Ab_{\aC} \to Ab_{\aC'}$ that takes $\aC$ into $\aC'$
and $w\aC$ into $w\aC'$.
\end{defn}

By abuse of language, we usually call $\aC$ the DG-Waldhausen
category. In the definition, the cone $CX$ and cocone $C'X$ of a
complex $X$ are the usual contractible complexes that fit into the
short exact sequences
\begin{gather*} 
0\to X\to CX\to X[1] \to 0\\
0\to X[-1]\to C'X\to X \to 0.
\end{gather*}
Waldhausen's saturation property on the weak equivalences means that
$w\aC$ satisfies ``two-out-of-three'': for composable maps $f,g$ in $\aC$,
if any two of the maps $f$, $g$, and $g\circ f$ are in $w\aC$ then so
is the third.  Waldhausen's extension property means that when 
\[
\xymatrix{%
0\ar[r]&X\ar[d]_{\htp}\ar[r]&Y\ar[d]\ar[r]&Z\ar[d]_{\htp}\ar[r]&0\\
0\ar[r]&X'\ar[r]&Y'\ar[r]&Z'\ar[r]&0
}
\]
is a commutative diagram of degree-wise split short exact sequences
with the maps $X\to X'$ and $Z\to Z'$ in $w\aC$, then the map $Y\to
Y'$ is in $w\aC$.  As a consequence, the subcategory $\aC^{w}$ of
$w\aC$-acyclic objects (those objects weakly equivalent to 0) is
closed under extensions; the extension property is equivalent
to this closure condition.

A DG-Waldhausen category obtains the structure of a pretriangulated
DG-category with the usual mapping complexes and also the structure of
a Waldhausen category (in fact a complicial biWaldhausen category)
with the cofibrations the degreewise-split monomorphisms.  Therefore
we can construct both its algebraic $K$-theory (using the Waldhausen
category structure), as well as its $THH$ and $TC$ (lifting the
DG-category structure to an associated spectral category structure).
The weak equivalences of the Waldhausen category structure specify
additional homotopical data beyond that in the mapping spectra: The
natural homotopy category of the Waldhausen category structure is the
localization of the homotopy category associated to the DG-category
obtained by localizing with respect to the weak equivalences.  In the
terminology of 
Section~\ref{secgenloc}, this homotopy category is the triangulated
quotient of $\aC$ by the subcategory $\aC^{w}$ of $w\aC$-acyclics.
Thus, the proper notion of $THH$ and $TC$ are the $THH$ and $TC$ of
the localization pair, $\CTHH(\aC/\aC^{w})$ and $CTC(\aC/\aC^{w})$.

We now review the construction of algebraic $K$-theory in preparation
for constructing the trace map.
Recall Waldhausen's $\Sdot$ construction produces a simplicial
Waldhausen category from a Waldhausen category.  In the case of a
DG-Waldhausen category $\aC$, the $\Sdot$ construction produces a
simplicial DG-Waldhausen category.
Let $\Ar[n]$ denote the
category with objects $(i,j)$ for $0\leq i\leq j\leq n$ and a unique
map $(i,j)\to (i',j')$ for $i\leq i'$ and $j\leq j'$.  $\Sdot[n]\aC$
is defined to be the full subcategory of the category of functors
$A\colon \Ar[n]\to \aC$ such that:
\begin{itemize}
\item $A_{i,i}=0$ for all $i$, 
\item The map $A_{i,j}\to A_{i,k}$ is a cofibration (degreewise-split
monomorphism) for all $i \leq 
j \leq k$, and
\item The diagram
\[  \xymatrix@-1pc{%
A_{i,j}\ar[r]\ar[d]&A_{i,k}\ar[d]\\A_{j,j}\ar[r]&A_{j,k}
} \]
is a pushout square for all $i \leq j \leq k$, 
\end{itemize}
where we write $A_{i,j}$ for $A(i,j)$.  The last two conditions can be
simplified to the hypothesis that each map $A_{0,j}\to A_{0,j+1}$ is a
cofibration and the induced maps $A_{0,j}/A_{0,i}\to A_{i,j}$ are
isomorphisms.  This becomes a DG-Waldhausen category by taking the
abelian category to be the category of functors $\Ar[n]\to Ab_{\aC}$
and defining a map $A\to B$ to be a weak equivalence when each
$A_{i,j}\to B_{i,j}$ is a weak equivalence in $\aC$.  Note that $A\to
B$ is a degreewise-split monomorphism when each $A_{i,j}\to B_{i,j}$
and each induced map $A_{i,k}\cup_{A_{i,j}}B_{i,j}\to B_{i,k}$ is a
degreewise-split monomorphism. 
An ordered map $\{ 1,\dotsc,m \}\to \{ 1,\dotsc,n \}$ induces a
functor $\Ar[m]\to \Ar[n]$ and hence a DG-exact functor $\Sdot[n]\aC\to
\Sdot[m]\aC$, making $\Sdot\aC$ a simplicial DG-Waldhausen category.
Because each $\Sdot[p]\aC$ is itself a DG-Waldhausen category, the
$\Sdot$ construction can be iterated to form multisimplicial
DG-Waldhausen categories.

For any DG-Waldhausen category $\aD$, let $w_{q}\aD$ denote the
DG-category whose objects consist of a sequence of $q$ composable weak
equivalences in $\aD$ (with $w_{0}\aD = \aD$).  Using this
construction and iterating the $\Sdot$ construction, we obtain
multi-simplicial DG-categories $w\subdot\Sdot^{(n)}\aC$.
The inclusion of $\aD$ as $\Sdot[1]\aD$ induces an $(n+2)$-simplicial
map
\[
\ob(w\subdot\Sdot^{(n)}\aC)\sma S^{1}\subdot \to
\ob(w\subdot\Sdot^{(n+1)}\aC),
\]
where $S^{1}\subdot$ denotes the standard simplicial model of the
circle (with one non-de\-gen\-er\-ate vertex and one non-degenerate 1-simplex).
These structure maps together with the natural $\Sigma_{n}$ action on
the categories $w_{q}\Sdot[p]^{(n)}\aC$ give the collection of simplicial sets
\[
\{ \diag\ob(w\subdot\Sdot^{(n)}\aC) \mid n\geq 0 \}
\]
the structure of a symmetric spectrum.  Waldhausen showed that the adjoint
attaching maps  
\[
|\diag\ob(w\subdot\Sdot^{(n)}\aC)| \to \Omega|\diag\ob(w\subdot\Sdot^{(n+1)}\aC)|
\]
are weak equivalences for $n>0$; i.e., the geometric realization is a
positive fibrant symmetric spectrum of topological spaces. 

\begin{defn}[Waldhausen]
$K\aC$ is the symmetric spectrum 
\[K\aC(n)=\diag \ob(w\subdot\Sdot^{(n)}\aC).\]
\end{defn}

To mix the $\Sdot$ construction with the cyclic bar construction, we
use the more convenient DG-categories $\barw_{q}\aD$ in place of the
DG-categories $w_{q}\aD$.  For a DG-Waldhausen category $\aD$, let
$\barw_{q}\aD$ be the DG-category that is the full subcategory of
$w_{q}\aD$ consisting of those objects where each weak
equivalence in the sequence is also a cofibration
(degreewise-split monomorphism).  The advantage of $\barw_{q}\aD$ over
$w_{q}\aD$ is that the limit defining its mapping complexes is a
homotopy limit.  Waldhausen also used
this construction; the following is 
a special case of Lemma~1.6.3 of \cite{WaldhausenKT}.

\begin{prop}
For a DG-Waldhausen category $\aD$, the inclusion of
$\ob(\barw\subdot\aD)$ in $\ob(w\subdot\aD)$ is a
weak equivalence.  
\end{prop}

For the construction of the cyclotomic trace, we use the associated
spectral category functor of Definition~\ref{defassoc} to lift the
multi-simplicial DG-categories $\barw\subdot\Sdot^{(n)}\aC$ to
multi-simplicial spectral categories, which by abuse, we denote with
the same notation.  For any space $X$, the spaces 
\[
|THH(\barw\subdot\Sdot^{(n)}\aC)(X)|
\]
then fit together into a symmetric spectrum (indexed on $n$) of
topological spaces.   For each orthogonal $S^{1}$-representation $V$
(q.v., Notation~\ref{notHM2})  let
$\WTHH(\aC)(V)$ be the symmetric spectrum defined by
\[
\WTHH(\aC)(V)(n)=|THH(\barw\subdot\Sdot^{(n)}\aC)(S^{V})|.
\]
As we let $V$ and $n$ vary, $\WTHH(\aC)$ has the structure of a
symmetric spectrum in the category of cyclotomic spectra.  Let $Q$ be
a $\Omega$-spectrum replacement functor in the category of cyclotomic
spectra that is enriched in based spaces as an endofunctor
(topologizing the mapping spaces in the category of cyclotomic spectra
as subspaces of the mapping spaces in the category of orthogonal
$S^{1}$-spectra).  Then applying $Q$ objectwise in symmetric spectra
to the cyclotomic spectra in $\WTHH(\aC)$, we obtain a weakly
equivalent symmetric spectrum of cyclotomic spectra where each
constituent orthogonal $S^{1}$-spectrum is an $\Omega$-spectrum.  Let
$\WT(\aC)=Q\WTHH(\aC)$.

The spectrum $\WT(\aC)\simeq \WTHH(\aC)$ lies between $K\aC$ and
$\CTHH(\aC/\aC^{w})$ 
in the stable category.
Write $\bar K\aC$ for the symmetric spectrum of topological
spaces $\bar K\aC(n)=|\ob\barw\subdot\Sdot^{(n)}\aC|$.  We then get a
symmetric spectrum of (non-equivariant) orthogonal spectra 
$F^{\orth}_{0}\bar K\aC(n)$ representing the same object in the stable
category using the free functor from spaces to orthogonal spectra. The
inclusion of objects (via the 
identity) induces a map of symmetric spectra of orthogonal spectra
\[
K\aC\htp F^{\orth}_{0}\bar K\aC \to \WTHH(\aC)\to \WT(\aC),
\]
natural in DG-exact functors.  Likewise, using the free functor from
spaces to symmetric spectra, we obtain a
map of symmetric spectra of orthogonal spectra
\[
THH(\aC)\htp F^{\symm}_{0}THH(\aC) \to \WTHH(\aC)\to \WT(\aC),
\]
natural in DG-exact functors, induced by the identification of
$w_{0}\Sdot^{(0)}\aC$ as $\aC$.  Finally, writing $CT(\aC/\aC^{w})$
for the cofiber
\[
CT(\aC/\aC^{w})=C(T(\aC^{w})\rightarrow T(\aC))
\]
(where, as in Definition~\ref{defTC}, $T(\aC)=QTHH(\aC)$),
we obtain the comparison map
\[
\CTHH(\aC/\aC^{w})\htp F^{\symm}_{0}CT(\aC/\aC^{w})\to \WT(\aC),
\]
as follows: The functor
$\aC^{w}\to \barw_{1}\aC$ that sends a $w\aC$-acyclic object $a$ to the
weak equivalence $0\to a$ induces a map from the cone on
$T(\aC^{w})$ to $|T(w\subdot(\aC))|$ that restricts on the
face $T(\aC^{w})$ to the inclusion of $T(\aC^{w})$ in
$T(\aC)=T(\barw_{0}\aC)$.  This then extends to the map from the
cofiber $CT(\aC/\aC^{w})$ above.  Similar observations
construct the 
symmetric spectrum of orthogonal spectra
\[
\WTC(\aC)=|TC(w\subdot\Sdot^{(-)}\aC)|
\]
and maps
\[
K\aC\htp F^{\orth}_{0}\bar K\aC\to \WTC(\aC)\from 
F^{\symm}_{0}CTC(\aC/\aC^{w}) \htp CTC(\aC/\aC^{w}).
\]
The following is the main theorem of this section and is proved below.

\begin{thm}\label{thmmainappb}
For a DG-Waldhausen category $\aC$, the maps 
\[
F^{\symm}_{0}CTHH(\aC/\aC^{w})\to \WTHH(\aC)
\qquad \text{and}\qquad 
F^{\symm}_{0}CTC(\aC/\aC^{w})\to \WTC(\aC)
\]
are level equivalences of symmetric spectra of orthogonal spectra.
\end{thm}

We can now define the trace.

\begin{defn}
The cyclotomic trace maps from $K$-theory to $TC$ and from $K$-theory
to $THH$ are the zigzags
\[
\xymatrix@R-1pc@C-1.125pc{%
&\WTC(\aC)\ar[dd]
&F^{\symm}_{0}CTC(\aC/\aC^{w})\htp CTC(\aC/\aC^{w})\hspace{9em}
\ar[l]_-{\htp}\ar[dd]<-7em>\ar[dd]<0em>\\
K\aC\htp F^{\orth}_{0}\bar K\aC\ar[ur]\ar[dr]\\
&\WT(\aC)
&F^{\symm}_{0}CT(\aC/\aC^{w})\htp CT(\aC/\aC^{w})\htp CTHH(\aC/\aC^{w})\ar[l]^-{\htp}\\}
\]
Every map in the diagram is natural in DG-exact functors.
\end{defn}

When we restrict to appropriate categories of schemes or pairs of
schemes as in \cite[\S6]{ThomasonTrobaugh}, we can factor the trace
above through non-connective $K$-theory.  Essentially, we take $\bar
K$ and $\WT$ as our model functors to spectra (which here would be
the point-set category of symmetric spectra of orthogonal spectra)
applied to the appropriate DG-Waldhausen category model for perfect
complexes (as in \cite[\S3]{ThomasonTrobaugh}), depending on the kind
of naturality required for the maps of schemes.  For any of these
models, we get natural pairings
\[
\bar K(X \on (X-U)) \sma K_{f}(\bZ[t,t^{-1}]) \to \bar K(X[t,t^{-1}] \on
(X[t,t^{-1}] -U[t,t^{-1}]))
\]
and
\begin{multline*}
\WTHH(X \on (X-U)) \sma K_{f}(\bZ[t,t^{-1}])\\ 
\to \WTHH(X[t,t^{-1}] \on (X[t,t^{-1}] -U[t,t^{-1}])),
\end{multline*}
where $K_{f}(\bZ[t,t^{-1}])$ denotes the Waldhausen $K$-theory
symmetric ring spectrum of the exact category with objects the
canonical free modules 
\[
0, \bZ[t,t^{-1}], (\bZ[t,t^{-1}])^{2}, (\bZ[t,t^{-1}])^{3} \dotsc .
\]
The arguments presented in Section~\ref{secnoncon} extend to this
context to construct the non-connective cyclotomic trace.

The remainder of the section proves Theorem~\ref{thmmainappb}.  A
version of the Additivity Theorem, as always, provides the key lemma.
Given DG-Waldhausen categories $\aA$,$\aB$,$\aC$ and DG-exact functors
$\phi\colon \aA\to\aB$, 
$\psi \colon \aC\to \aB$, let $\aE(\aA,\aB,\aC)$ be the 
DG-Waldhausen category where an object consists of:
\begin{enumerate}
\item A tuple $(a,b,c)$ of objects $a\in\aA$, $b\in\aB$, and $c\in\aC$, and
\item A degreewise-split short exact sequence in $\aB$,
\[
0\to \phi a\to b\to \psi c \to 0.
\]
\end{enumerate}
The mapping complex in $\aE(\aA,\aB,\aC)$ from $(a,b,c)$
to $(a',b',c')$ is 
\begin{equation}\label{eqmapcxa}
(\aA(a,a')\times_{\aB(\phi a,b')}\aB(b,b'))\times_{\aB(\psi c,\psi c')}\aC(c,c'),
\end{equation}
which is isomorphic to
\begin{equation}\label{eqmapcxb}
\aA(a,a')\times_{\aB(\phi a,\phi a')}(\aB(b,b')\times_{\aB(b,\psi c')}\aC(c,c')).
\end{equation}
Note that each of the maps
\begin{gather*}
\aB(b,b')\to \aB(\phi a,b'), \\
\aB(b,b')\to \aB(b,\psi c'), \\
\aA(a,a')\times_{\aB(\phi a,b')}\aB(b,b')\to \aB(\psi c,\psi c'),\\
\aB(b,b')\times_{\aB(b,\psi c')}\aC(c,c')\to \aB(\phi a,\phi a')
\end{gather*}
is a degreewise-split epimorphism, and so the limits
in~\eqref{eqmapcxa} and~\eqref{eqmapcxb} are homotopy limits.

We have DG-exact functors
\begin{gather*}
\alpha \colon \aE(\aA,\aB,\aC)\to \aA\\
\beta \colon \aE(\aA,\aB,\aC)\to \aB\\
\gamma \colon \aE(\aA,\aB,\aC)\to \aC
\end{gather*}
induced by the forgetful functor and a DG-exact functor
\[
\sigma \colon \aA\times \aC\to \aE(\aA,\aB,\aC)
\]
induced by $\sigma(a,c)=(a,\phi a\oplus \psi c,c)$ (and the split
short exact sequence).  
The version of the additivity theorem we prove compares the maps
induced on $THH$ by $\sigma$ and $\alpha \vee \gamma$.

\begin{thm}[Additivity Theorem]
The functors 
\[
THH(\aA)\vee THH(\aC)\to THH(\aE(\aA,\aB,\aC))\to THH(\aA)\times  THH(\aC)
\]
induced by $\sigma$ and $\alpha \times \gamma$ are inverse weak equivalences.
\end{thm}

\begin{proof}
Consider the DG-exact functor $\phi'\colon \aA\to \aE(\aA,\aB,\aC)$
that takes $a$ in $\aA$ to $(a,\phi a,0)$.  By~\eqref{eqmapcxa}, we see
that this is a DK-embedding.  Now by Theorem~\ref{thmgenone}, it
suffices to show that the functor $\psi'\colon \aB\to
\aE(\aA,\aB,\aC)$ (sending $c$ to $(0,\psi c,c)$) induces an
equivalence from homotopy category $\pi_{0}\aB$ to the triangulated
quotient $\pi_{0}\aE(\aA,\aB,\aC)/\pi_{0}\aA$.  This is a
straightforward calculation from~\eqref{eqmapcxb}.
\end{proof}

We can apply this to understand the effect both of $\barw_{q}$ and
$\Sdot[p]$ on $THH$.  An element of $\barw_{q}$ of $\aC$ is a sequence
of degreewise-split maps 
\[
c_{0}\to \dotsb \to c_{q}
\]
such that each quotient $c_{i}/c_{i-1}$ is in $\aC^{w}$.  Choosing
quotients, we get a DK-equivalent DG-category $\bar W_{q}$ that is a
DG-Waldhausen category.  Furthermore, we can
identify $\bar W_{q+1}\aC$ as $\aE(\bar W_{q}\aC,\aC,\aC^{w})$, for the
functor $\phi \colon \bar W_{q}\aC\to \aC$ that sends the sequence
pictured above to $c_{q}$.  As a consequence we get the following corollary.

\begin{cor}\label{corappbone}
For all $q$, the map 
\[
\underbrace{THH(\aC^{w})\vee \dotsb \vee THH(\aC^{w})}_{\text{$q$ factors}}
\vee THH(\aC)
\to THH(\barw_{q}\aC)
\]
induced by the map that sends $(a_{1},\dotsc,a_{q},c)$ to 
\[
c\to c\oplus a_{1}\to \dotsb \to c\oplus a_{1}\oplus \dotsb \oplus a_{q}
\]
is a weak equivalence.
\end{cor}

Similarly, for any DG-Waldhausen category $\aD$, the DG-Waldhausen category
$\Sdot[r]\aD$ is DK-equivalent (via a DG-exact functor) to 
the DG-category $\aE(\Sdot[r-1]\aD, \aD, \aD)$ for the
functor $\phi \colon S_{r}\aD\to \aD$ that takes $\{A_{i,j}\}$ to
$A_{0,r-1}$.  We use this observation to prove the following
corollary.

\begin{cor}\label{corappbtwo}
For each $n$ and $q$, the map
\[
\Sigma |THH(\barw_{q}\Sdot^{(n)}\aC)|\to |THH(\barw_{q}\Sdot^{(n+1)}\aC)|
\]
is a weak equivalence.
\end{cor}

\begin{proof}
We can write
$\Sigma |THH(\barw_{q}\Sdot^{(n)}\aC)|$ as the geometric realization
of a multi-sim\-pli\-cial object with one more simplicial direction,
$THH(\barw_{q}\Sdot^{(n)}\aC)\sma S^{1}\subdot$,
where $S^{1}\subdot$ denotes the standard simplicial model of the circle.
The map in the statement is induced by the map on geometric realizations of
the map of multi-simplical objects
\[
THH(\barw_{q}\Sdot^{(n)}\aC)\sma S^{1}\subdot
\to
THH(\barw_{q}\Sdot^{(n+1)}\aC).
\]
Using the standard isomorphisms 
\[
\Sdot[r]\barw_{q}\iso
\barw_{q}\Sdot[r], \qquad 
\Sdot[r]\Sdot[p]\iso \Sdot[p]\Sdot[r],
\]
and writing $\aD=\barw_{q}\Sdot[p_{1}]\dotsb \Sdot[p_{n}]\aC$, we are
looking at maps of the form
\[
\bigvee_{r}THH(\aD)\to THH(\Sdot[r]\aD).
\]
Using the relationship of $\Sdot[r]\aD$ and
$\aE(\Sdot[r-1]\aD,\aD,\aD)$ as above, we see by induction that this
map is a weak equivalence.
\end{proof}

Combining these two corollaries, we prove Theorem~\ref{thmmainappb}.

\begin{proof}[Proof of Theorem~\ref{thmmainappb}]
We can identify the map $\CTHH(\aC/\aC^{w})\to |THH(\barw\subdot
\aC^{w})|$ above as the induced map on geometric realization
of the map of simplicial objects
\[
\underbrace{THH(\aC^{w})\vee \dotsb \vee THH(\aC^{w})}_{\text{$\ssdot$ factors}}
\vee THH(\aC)
\to THH(\barw\subdot\aC),
\]
and is a weak equivalence by Corollary~\ref{corappbone}.  For $n>0$,
the $n$-th level of the symmetric spectrum of orthogonal spectra
$F^{\symm}_{0}\CTHH(\aC/\aC^{w})$ is $\Sigma^{n}\CTHH(\aC/\aC^{w})$.  It
now follows from Corollary~\ref{corappbtwo}, that the map 
\[
F^{\symm}_{0}\CTHH(\aC/\aC^{w})\to \WTHH(\aC)
\]
is a level equivalence of symmetric spectra of orthogonal spectra and
Theorem~\ref{thmmainappb} follows. 
\end{proof}

\section{$THH$ and $TC$ of small spectral model categories}
\label{appcofiber}

Our treatment of $THH$ and $TC$ of spectral categories in the rest of
the paper took the perspective that all the homotopy information
is encoded in the mapping spectra.  In the context of closed model
categories enriched over symmetric spectra, the weak equivalences
encode an additional localization.  We can extract a spectral category
satisfying the hypotheses of the main discussion of the paper from
such a model category by restricting to the full spectral subcategory
of cofibrant-fibrant objects.  However, this subcategory is not
usually preserved by naturally-occurring functors between model
categories, which tend to preserve only cofibrant or only fibrant
objects.

In this section, we present a construction of $THH$ of a small
spectral model category in terms of either the full subcategory of
cofibrant or fibrant objects.  (Here we must use the original
convention of Quillen that closed model categories are closed under
finite limits and colimits rather than the modern convention that they
are closed under all small limits and colimits.)
The construction is in terms of a ``cofiber $THH$''
description, exactly as in the $THH$ of localization pairs constructed
in Section~\ref{secgenloc}.  Since the quotient of the subcategory
of cofibrants by the acyclic cofibrants is the homotopy category of
the model category, we can regard this pair as analogous to a
localization pair, although it may not satisfy the hypotheses of the
definition.  Nevertheless, a similar (but easier) proof applies to
compare the $THH$ of this pair to the $THH$ of the cofibrant-fibrants.
The main theorem of this section is the following.

\begin{thm}
Let $\aM$ be a small closed model category enriched over symmetric
spectra, satisfying the symmetric spectrum version of SM7.  Write
$\aA$ for the subcategory of acyclic objects (objects weakly
equivalent to the zero object), and subscripts $c$ and $f$ for the
subcategories of cofibrant and fibrant objects, respectively, of $\aM$
and likewise of $\aA$.  In the following diagram, the
vertical map is always a weak equivalence, the left-hand map is a weak
equivalence if $\aM$ is left proper, and the right-hand map is a weak
equivalence if $\aM$ is right proper.
\[
\xymatrix{%
&THH(\aM_{cf})\ar[d]\\
\CTHH(\aM_{f}/\aA_{f})&\CTHH(\aM_{cf}/\aA_{cf})\ar[r]\ar[l]&
\CTHH(\aM_{c}/\aA_{c})
}
\]
\end{thm}

Since for any pair of objects in $\aA_{cf}$, the symmetric spectrum of maps
is homotopically trivial, $THH(\aA_{cf})$ is homotopically trivial,
and it then follows that the vertical map is a weak equivalence.

Of the remaining statements in the theorem, we treat the case of the
right horizontal map in detail; the case of the left horizontal map is
similar (and in fact follows by considering the opposite category).
Let $\aM'=\aM^{\Cell}$ be the pointwise cofibrant spectral category
weakly equivalent to $\aM$ from Proposition~\ref{propcofrep}, and let
$\aM'_{c}$, $\aA'$, and $\aA'_{c}$ be the appropriate subcategories.
As in Section~\ref{secgenloc}, we define the
$(\aM'_{c},\aM'_{c})$-bimodule $\fLM$ by
\[
\fLM(x,y)=\TB(\aM'_{c}(-,y);\aA'_{c};\aM'_{c}(x,-))
\]
and $\fQM$ as the cofiber of the composition map $\fLM\to\aM'_{c}$. 
The following lemma lists the properties of $\fQM$ we need in the
proof of the theorem.

\begin{lem}
Let $x$ be an object of $\aM_{c}$.
\begin{enumerate}
\item For $y$ in $\aM_{cf}$, the map $\aM'_{c}(x,y)\to \fQM(x,y)$ is a
weak equivalence.
\item If $\aM$ is right proper, then $\fQM(x,-)$ preserves weak
equivalences.
\end{enumerate}
\end{lem}

\begin{proof}
Since the mapping spectrum from a cofibrant acyclic
object to a fibrant object is homotopically trivial, for any object $y$ in
$\aM_{cf}$, $\fLM(x,y)$ is homotopically trivial, and the map
$\aM'_{c}(x,y)\to \fQM(x,y)$ is a weak equivalence. This proves~(i).
To prove~(ii), it suffices to show that for any object $y$ and
any fibrant replacement $y\to y'$, the map $\fQM(x,y)\to \fQM(x,y')$
is a weak equivalence.  Factor the initial map $*\to y'$ as an acyclic
cofibration followed by a fibration $a'\to y'$, and let $a$ be a
cofibrant replacement of the
pullback $y\times_{y'}a'$.
\[
\xymatrix{%
a\ar@{->>}[r]^(.4){\htp}
&y\times_{y'}a'\ar[r]\ar@{->>}[d]&a'\ar@{->>}[d]\\
&y\ar[r]_{\htp}&y'
}
\]
We obtain from this fibration pullback square (and the symmetric
spectrum version of SM7) the homotopy (co)cartesian
square of $\aM_{c}$-modules on the left below, and from this, the homotopy
cocartesian square of $\aM'_{c}$-modules on the right below.
\[
\xymatrix{%
\aM_{c}(-,a)\ar[r]\ar@{->>}[d]&\aM_{c}(-,a')\ar@{->>}[d]\\
\aM_{c}(-,y)\ar[r]&\aM_{c}(-,y')
}\qquad
\xymatrix{%
\aM'_{c}(-,a)\ar[r]\ar[d]&\aM'_{c}(-,a')\ar[d]\\
\aM'_{c}(-,y)\ar[r]&\aM'_{c}(-,y')
}
\]
Looking at the construction of $\fLM$, the bar construction preserves
homotopy cocartesian squares in either variable, and so we see that
the square on the left below is homotopy cocartesian; it follows that
the square on the right below is homotopy cocartesian.
\[
\xymatrix{%
\fLM(x,a)\ar[r]\ar[d]&\fLM(x,a')\ar[d]\\
\fLM(x,y)\ar[r]&\fLM(x,y')
}\qquad 
\xymatrix{%
\fQM(x,a)\ar[r]\ar[d]&\fQM(x,a')\ar[d]\\
\fQM(x,y)\ar[r]&\fQM(x,y')
}
\]
The hypothesis that $\aM$
is right proper implies that the map $a\to a'$ is a weak equivalence
and therefore that $a$ is in $\aA_{c}$.   It follows that $\fQM(x,a)$
and $\fQM(x,a')$ are homotopically trivial, and that
$\fQM(x,y)\to\fQM(x,y')$ is a 
weak equivalence.
\end{proof}

As in Section~\ref{secgenloc}, we have a natural weak equivalence
relating 
\[
\CTHH(\aM_{c}/\aA_{c})\simeq
\CTHH(\aM'_{c}/\aA'_{c})
\]
with $THH(\aM_{c};\fQM)$,
compatible with the map from $THH(\aM'_{c})\simeq THH(\aM_{c})$.
Applying part~(i) of 
the lemma, to complete the proof of the theorem, it suffices to show
that the map
\[
THH(\aM'_{cf};\fQM)\to THH(\aM'_{c};\fQM)
\]
is a weak equivalence.

Our strategy as in Section~\ref{secgenloc} is to apply the Dennis-Waldhausen
Morita argument~\ref{propcoremorita} to reduce to proving an
objectwise statement.  For this, we use the weak equivalence of
$(\aM'_{c},\aM'_{c})$-bimodules 
$\TB(\aM'_{c};\aM'_{c};\fQM)\to \fQM$.  Then the Dennis-Waldhausen
Morita argument gives us a weak equivalence 
\[
THH(\aM'_{cf};\TB(\aM'_{c};\aM'_{c};\fQM))\htp
THH(\aM'_{c};\TB(\fQM;\aM'_{cf};\aM'_{c}))
\]
and likewise the analogous weak equivalence with $\aM'_{cf}$ replaced
by $\aM'_{c}$.  Now it suffices to show that the map
\[
THH(\aM'_{c};\TB(\fQM;\aM'_{cf};\aM'_{c})) \to
THH(\aM'_{c};\TB(\fQM;\aM'_{c};\aM'_{c}))
\]
is a weak equivalence.  This map is induced by the map of
$(\aM'_{c},\aM'_{c})$-bimodules 
\[
\TB(\fQM;\aM'_{cf};\aM'_{c}) \to \TB(\fQM;\aM'_{c};\aM'_{c})
\]
and so it suffices to show that the map
\[
\TB(\fQM(-,y);\aM'_{cf};\aM'_{c}(x,-))\to \fQM(x,y)
\]
is a weak equivalence for each pair of objects $x$,$y$ in $\aM_{c}$.
This is clear from the Two-Sided Bar Construction Lemma and part~(i)
of the lemma when $y$ is in
$\aM_{cf}$.  It then follows for 
arbitrary $y$ by part~(ii) of the lemma.



\end{document}